\documentclass[a4paper,11pt,twoside]{article}

\usepackage[dvips]{graphics}
\usepackage{latexsym}
\usepackage{amsmath}
\usepackage{amsfonts}
\usepackage{amssymb}
\usepackage{amscd}
\usepackage[mathscr]{eucal}
\usepackage{epsfig}

\def\Dim{\emph{Proof : }}
\def\cvd{\nopagebreak\par\rightline{$_\blacksquare$}}

\def\Kn#1#2{\mathbb{#1}^{#2}}                 
\def\E#1#2{\left\langle #1,#2 \right\rangle}  
\def\fut{\mathrm{I}^+}                         
\def\pass{\mathrm{I}^-}                         
\def\cfut{\mathrm{J}^+}                         
\def\cpass{\mathrm{J}^-}   

\def\OO{\textrm{O}}
\def\SOO{\textrm{SO}}

\def\ISO{\textrm{Iso}}

\def\Diffeo{\textrm{Diffeo}}

\def\ort#1{#1^\perp}                          
\def\eps{\varepsilon}                        

\def\Ftilde{\widetilde{F}_\tau}                    
\def\dom{\mathcal D_\tau}             
\def\dompass{\mathcal D^-_\tau}       

\def\tr{\mathrm{tr}}                          
\def\stab{\mathrm{stab}}                      
\def\graph{\mathrm{graph}}

\def\coom{\textrm{H}^}                         

\def\der#1#2{\frac{\partial#2}{\partial y_#1}}
\def\supp{\mathrm{supp}}

\newtheorem{defi}{Def.}[section]
\newtheorem{teo}{Theorem}[section]
\newtheorem{prop}[teo]{Proposition}
\newtheorem{lemma}[teo]{Lemma}
\newtheorem{corollario}[teo]{Corollary}
\newtheorem{oss}[teo]{Remark}
\newtheorem{es}{Example}[section]

\begin{document}

\title{Flat  Spacetimes  with Compact \\ Hyperbolic Cauchy Surfaces}
\author{Francesco Bonsante}
\maketitle


\section{Introduction}
We study the flat $(n+1)$-spacetimes $Y$ admitting a Cauchy surface
diffeomorphic to a compact hyperbolic $n$-manifold $M$. Roughly speaking, we
show how to construct a canonical future complete one, $Y_\rho$, among all
such spacetimes sharing a same holonomy $\rho$.
We study the geometry of $Y_\rho$ in terms of its \emph{canonical cosmological
  time} (CT). In particular we study the asymptotic behaviour of the level
surfaces of the cosmological time.\par
The present work generalizes the case $n=2$ treated in \cite{Mess} taking from
\cite{Ben-Guad} the emphasis on the fundamental r\^ole played by the canonical
cosmological time. 
 In particular Mess showed that if $F$ is a closed surface of genus
$g\geq 2$ then the linear holonomies of the Lorentzian flat structures on $\mathbb
R\times F$ such that $\{0\}\times F$ is a spacelike surface
are faithful and discrete representations of
$\pi_1(F)$ in $\SOO^+(2,1)$.\par
Moreover he proved that every representation  $\varphi: \pi_1(F)\rightarrow\ISO(\Kn{M}{3})$ 
whose linear part is faithful and discrete is a holonomy for some Lorentzian
structure on $\mathbb R\times F$ such that $\{0\}\times F$ is a spacelike
surface.
In particular he showed that there exists a unique maximal future complete
convex domain of $\Kn{M}{3}$,
called \emph{domain of dependence},  which is $\varphi(\pi_1(F))$-invariant  such that the quotient is
a globally hyperbolic manifold homeomorphic to $\mathbb R\times F$ with
regular cosmological time.\par
If we fix the linear holonomy $f:\pi_1(F)\rightarrow\SOO^+(2,1)$ 
these domains (and so the affine deformations of the representation $f$) 
are parametrized by the measured geodesic laminations on
$\Kn{H}{2}/f(\pi_1(F))$. The link between domains of dependence and
measured geodesic laminations is the \emph{Gauss map} of the CT-level surfaces.\par
On the other hand Benedetti and Guadagnini noticed in \cite{Ben-Guad} 
the \emph{singularity in the past} of a domain of dependence is a real tree
which is dual to the lamination. Moreover they argued that the action of
$\pi_1(F)$ on the CT-level
surface $\widetilde S_a=T^{-1}(a)$ converges in the Gromov sense to the action of
$\pi_1(F)$ on the singularity for
$a\rightarrow 0$ and to the action of $\pi_1(F)$ on the hyperbolic plane $\Kn{H}{2}$ for
$a\rightarrow\infty$.
Thus the asymptotic states of the cosmological time materialize the duality
between geometric real trees (realized by the singularity in the past) and the
measured geodesic laminations in the hyperbolic surface $F$, according to
Skora's theorem \cite{Skora}.\par
In this paper we try to generalize this approach in higher dimension. 
By extending Mess' method to any dimension $n$ 
we associate to each $Y_\rho$  a $\Gamma$-invariant \emph{geodesic
  stratification} of $\Kn{H}{n}$ (see section \ref{CT-section} for the
definition), and we discuss the duality between geodesic stratifications and
singularities in the past.
In particular we recover this duality at least for an interesting class of so
called spacetimes with \emph{simplicial singularity in the past}.\\

We briefly describe the contents of the paper in section
\ref{preliminaries-section}. We first remind some basic facts about Lorentzian
spacetimes and hyperbolic manifolds, then we give a quite articulated
statement of the main results. 
Some assertions will be fully described and proved in the following sections.
In the last section we shall discuss some related questions and open problems.
\section{Preliminaries and Statement of the Main Theorem}\label{preliminaries-section}

In this section we recall few basic facts about Lorentzian geometry, the geometry of the Minkowski space, and
hyperbolic space. An exaustive trateament about Lorentzian geometry, including
a careful analysis of global causality questions,  can be
found in \cite{Hawking} or in \cite{Beem}. 
For an introduction to  hyperbolic space, see \cite{Benedetti}.
In the last part of the section we state the main theorem which we shall prove
in the following sections.
\paragraph{Spacetimes.}
A \emph{Lorentzian $(n+1)$-manifold} $(M,\eta)$ is given by a smooth
$(n+1)$-manifold $M$ (this includes the topological assumption that $M$ is
metrizable and with countable basis) and
a symmetric non-degenerate $2$-form $\eta$ with signature equal to $(n,1)$. A basis
$(e_0,\ldots,e_n)$ of a tangent space $T_pM$ is \emph{ orthonormal} if the matrix
of $\eta_p$ with respect to this basis is $diag(-1,1,\ldots,1)$. A
tangent vector $v$ is  \emph{spacelike}
(resp. \emph{timelike}, \emph{null}, \emph{non-spacelike}) if $\eta(v,v)$ is
positive (resp. negative, null, non-positive). 
A $\mathrm{C}^1$-curve in $M$ is  \emph{chronological}
(resp. \emph{causal}) if the  speed vector is timelike (resp. non-spacelike).\par
Let $M$ be a Lorentzian connected $(n+1)$-manifold. Consider the subset $\mathcal{C}$ of
the tangent bundle $TM$ formed by the timelike tangent vectors: it turns out
that either $\mathcal{C}$ is connected or it has two connected
components. We say that $M$ is \emph{time-orientable} if $\mathcal{C}$ has two
connected components. A time-orientation is a choice of one of these
components. A \textbf{spacetime} is a Lorentzian connected time-orientable
manifold provided with a time-orientation. Let  $M$ be a space-time and
$\mathcal{C}_+$ be the chosen component. A non-spacelike 
tangent vector is  \emph{future-directed} (resp. \emph{past
  directed}) if it is not the null vector and it lies (resp. it does not
lie) in the closure of $\mathcal{C}_+$.
A causal curve is  \emph{future-directed} if its speed vector is future
directed.
Let $p\in M$, \emph{the future of $p$} (resp. \emph{the past}) is the subset
$\fut(p)$ (risp. $\pass(p)$) of $M$ formed by the points of $M$ which are
terminal end-points of future-directed (resp. past-directed) chronological
curves which start in $p$.
If we replace the chronological curves by causal curves we obtain the
\emph{causal future} $\cfut(p)$ (resp. \emph{causal past} $\cpass(p)$ ) of
$p$.
Let $\gamma:I\rightarrow M$ be a causal curve. The Lorentzian lenght of
$\gamma$ is
\[
      \ell(\gamma):=\int_I\sqrt{-\eta(\gamma'(t),\gamma'(t))}\mathrm{d}t.
\] 
Given $p\in M$ and  $q\in\cpass(p)$ the Lorentzian distance between $p$ and $q$ is 
\[
      d(p, q):=\{\sup\ell(\gamma)|\gamma\textrm{ is a causal curve whose
        endpoints are } p \textrm{ and } q\}.
\]
For every $p\in M$ we can take
\[
   \tau(p):=\sup\{d(p, q)|q\in\cpass(p)\}.
\]
This define a function
\[
  \tau:M\rightarrow\mathbb R\cup\{+\infty\}
\]
which can be very degenerate (for istance $\tau\equiv+\infty$ on the Minkowski space
$\Kn{M}{n+1}$).
We are interested in the spacetimes such that the function $\tau$ is a
\textbf{canonical cosmological time} (CT): this means that $\tau$ is finite and
increasing on every directed causal curves (i.e. $\tau$ is a \emph{time}) and
it is \emph{regular}, that is $\tau$ tends to $0$ over every inextendable past
directed causal curve. Spacetimes with regular cosmological time have been
pointed out and studied in \cite{Andersson}. We recall the following general result.
\begin{teo}
Suppose that $M$ is a spacetime with a regular cosmological time $\tau$.
Then $\tau$ is twice differentiable  almost everywhere. Moreovere for every
$p\in M$ there is an inextendable in the past timelike  geodesic
$\gamma:(0,\tau(p)]\rightarrow M$ which has terminal point in $p$ and such
that $\tau(\gamma(t))=t$. The
level surfaces $S_a=\tau^{-1}(a)$ are future Cauchy surfaces.
\end{teo}
 We recall that a \textbf{Cauchy surface} is an embedded hypersurface $S$ of $M$ such that every
inextendable causal curve in $M$ intersects $S$ in a unique point.\par
Finally we recall that a spacetime $M$ is \textbf{globally hyperbolic} if for every
$p,q\in M$ the set $\cfut(p)\cap\cpass(q)$ is compact. It is the strongest
global causality assumption and implies strict constraints on the topology of
$M$. In particular in \cite{Geroch} it is shown that $M$ is globally hyperbolic
if and only if there is a Cauchy surface $S$ in $M$, and in this case $M$ is
homeomorphic to $\mathbb R\times S$.
\paragraph{Minkowski space.}
The \emph{Minkowski $(n+1)$-space-time} $\Kn{M}{n+1}$ is the flat simply
connected complete Lorentzian $(n+1)$-manifold (it is unique up to isometry).
 Let $(x_0,\ldots,x_n)$ be the natural coordinates
on $\Kn{R}{n+1}$, then a concrete model for $\Kn{M}{n+1}$ is  $\Kn{R}{n+1}$
provided with the Lorentz form
\[
     \eta=-\mathrm{d}x_0^2+\mathrm{d}x_1^2+\ldots+\mathrm{d}x_n^2.
\]
In what follows we shall always use this model.
Notice that the frame $\left(e_i=\frac{\partial}{\partial x_i}\right)_{i=0,\ldots,n}$ is
parallel and orthonormal. Thus we can identify in a standard way the tangent
space $(T_x\Kn{M}{n+1},\eta_x)$ with $\Kn{R}{n+1}$ provided with the scalar
product $\E{\cdot}{\cdot}$ defined by the rule
 $\E{v}{w}=-v_0w_0+v_1w_1+\ldots+v_nw_n$. \par 
Minkowski space is an orientable and  time-orientable  Lorentz manifold. Let us put on
it the standard orientation (such that the canonical basis $(e_0,\ldots,e_n)$
is positive) and the standard time-orientation (a timelike tangent vector $v$
is future directed if $\E{v}{e_0}<0$).
By \emph{orthonormal affine coordinates} we mean a set $(y_0,\ldots,y_n)$ of
affine coordinates on $\Kn{M}{n+1}$ such that the frame
$\{\frac{\partial}{\partial y_i}\}$ is \emph{orthonormal and positive} and the vector 
$\frac{\partial}{\partial y_0}$ is \emph{future directed}.\par
Consider the isometry group of $\Kn{M}{n+1}$. 
It is easy to see that
$f$ is an isometry of $\Kn{M}{n+1}$ if and only if
it is affine and $\mathrm{d}f(0)$ belongs to the group $\OO(n,1)$ 
of the
linear trasformations of $\Kn{R}{n+1}$ which preserve the scalar product $\E{\cdot}{\cdot}$.  
It follows that the group of isometries of $\Kn{M}{n+1}$ is generated by
$\OO(n,1)$ and the group of traslations $\Kn{R}{n+1}$. Furthermore
$\Kn{R}{n+1}$ is a normal subgroup of $\ISO(\Kn{M}{n+1})$ (in fact it is the
kernel of the map $\ISO(\Kn{M}{n+1})\ni f\mapsto\mathrm{d}f(0)\in\OO(n,1)$) so
that $\ISO(\Kn{M}{n+1})$ is isomorphic to $\Kn{R}{n+1}\rtimes\OO(n,1)$.\par
 Notice that $\OO(n,1)$ is the isotropy group of
$0$ in $\ISO(\Kn{M}{n+1})$. It is a semisimple Lie group and it has four connected
components.
The connected component of the identity $\SOO^+(n,1)$ is the group of linear trasformations
which preserve  orientation  and time-orientation. It is the 
\emph{Lorentz group}.  There are two proper subgroups which contain
$\SOO^+(n,1)$: the group $\SOO(n,1)$ of linear isometries which preserve
orientation of $\Kn{M}{n+1}$ and the group $\OO^+(n,1)$ of linear isometries which preserve
time-orientation of $\Kn{M}{n+1}$. In each of these groups the index of $\SOO^+(n,1)$ is
$2$.  \par
Geodesics in $\Kn{M}{n+1}$ are  straight lines. 
There are three types of geodesics up to isometry: spacelike, timelike and
null. Notice that they are classified by the restriction of the
form $\eta$ on them. The totally geodesic $k$-planes are the affine $k$-planes
in $\Kn{M}{n+1}$.
Also the $k$-planes are classified up to isometry by the restriction
of the Lorentz form on them. Hence we say that a $k$-plane $P$ is \emph{spacelike} if 
 $\eta|_P$ is a flat Riemannian form, we say that it is \emph{timelike} if
 $\eta|_P$ is a flat Lorentz form, finally we say that $P$ is \emph{null} if
 $\eta|_P$ is a degenerated form.
\paragraph{Hyperbolic space.}
Let $\Kn{H}{n}$ be the set of the points in the future of $0$ which have
Lorentzian distance from $0$ equal to $1$. If we identify $\Kn{M}{n+1}$ with
the tangent space $T_0\Kn{M}{n+1}$ via the exponential map we get
\[
     \Kn{H}{n}=\{x\in\Kn{M}{n+1}|\E{x}{x}=-1,\quad x_0>0\}.
\]
It follows that the tangent space $T_x\Kn{H}{n}$ is the space $\ort x$. 
Since $x$ is timelike $T_x\Kn{H}{n}$ is spacelike so
that $\Kn{H}{n}$ has a natural Riemannian structure. An easy calculation shows that $\Kn{H}{n}$ is the
simply connected complete Riemannian manifold with constant sectional curvature
equal to $-1$.\par
A geodesic in $\Kn{H}{n}$ is the intersection of $\Kn{H}{n}$ with a timelike $2$-plane which
passes through $0$.
More generally a totally geodesic $k$-submanifold ($k$-plane) of $\Kn{H}{n}$ is the
intersection of $\Kn{H}{n}$ with a timelike $(k+1)$-plane which passes through $0$  (notice that such
intersection is  transverse).
Thus it follows that a subset $C$ of $\Kn{H}{n}$ is convex if and
only if it is the intersection of $\Kn{H}{n}$ with  a convex cone with apex in $0$.\par
Clearly $\Kn{H}{n}$ is invariant for the group $\OO^+(n,1)$ and
furthermore this group acts by isometries on $\Kn{H}{n}$. It can be shown that
$\OO^+(n,1)$ is in fact the isometry group of $\Kn{H}{n}$.
The group of orientation-preserving isometries of $\Kn{H}{n}$  is
identified with $\SOO^+(n,1)$.\par
Let $\Kn{P}{n}$ be the set of the lines which passes through $0$ and 
$\pi:\Kn{M}{n+1}\rightarrow\Kn{P}{n}$ the natural projection. Then
$\pi|_{\Kn{H}{n}}$ is a diffeomorphism of $\Kn{H}{n}$ onto the set of time-like lines. 
The closure of this set is a closed ball and its boundary is formed by the set
of the null lines. Let $\partial\Kn{H}{n}$ be the set of null lines and put on
$\overline{\Kn{H}{n}}:=\Kn{H}{n}\cup\partial\Kn{H}{n}$ the topology which makes
the natural map $\pi:\overline{\Kn{H}{n}}\rightarrow\Kn{P}{n}$ a homeomorphism
onto its image. Notice that every $g\in\OO^+(n,1)$ extends uniquely to a homeomorphism
of $\overline{\Kn{H}{n}}$.\par
Now we can  classify the elements of $\OO^+(n,1)$.
We say that $g\in\OO^+(n,1)$ is \emph{elliptic} if $g$ has a timelike
eigenvector (and in this case the respective eigenvalue is $1$).
We say that $g$ is \emph{parabolic} if it is not elliptic and it has a unique
null eigenvector (and in this case the relative eigenvalue is $1$).
Finally we say that $g$ is \emph{hyperbolic} if it is not elliptic and it has two
null eigenvectors (in this case there exists $\lambda>1$ such that the
respective eigenvalues are $\lambda$ and $\lambda^{-1}$).
Notice that if $g$ is hyperbolic there exists a unique geodesic $\gamma$ in
$\Kn{H}{n}$ which is invariant for $g$. In this case $\gamma$ is called the
\emph{axis} of $g$.
\paragraph{Geometric structures.} 
We shall consider only \emph{oriented}
manifolds or spacetimes. We shall be concerned with \emph{hyperbolic
  $n$-manifolds} (i.e. Riemannian $n$-manifolds locally isometric to $\Kn{H}{n}$) and with \emph{flat}
$(n+1)$-\emph{ spacetimes} (i.e. spacetimes locally isometric to $\Kn{M}{n+1}$).
By using the convenient setting of $(X,G)$-manifolds (see for istance chap. B
of \cite{Benedetti}) we can say that hyperbolic manifolds and flat spacetimes
are, by definition, $(X,G)$-manifolds where $(X,G)$ is respectively
$(\Kn{H}{n},\SOO^+(n,1))$ and $(\Kn{M}{n+1},\ISO(\Kn{M}{n+1})$).\par
We summarize few basic facts about such $(X,G)$-manifolds. Let $N$ be a
$(X,G)$-manifold and fix an universal covering $\pi:\widetilde N\rightarrow N$. Then
the $(X,G)$-structure on $N$ lifts to a $(X,G)$-structure on $\widetilde N$ such
that:
\begin{enumerate}
\item
the covering map $\pi$ is a local isometry;
\item
the group $\pi_1(N, x_0)$ acts by isometries on $\widetilde N$ in such way that
$N=\widetilde N/\pi_1(N, x_0)$ and $\pi$ is identified with the quotient map.
\end{enumerate}
Let us summarize by saying that $\pi:\widetilde N\rightarrow N$ is a
\emph{$(X,G)$-universal covering}.
\begin{prop}
Given a $(X,G)$-universal covering $\pi:\widetilde N\rightarrow N$ there exists a
pair $(D,\rho)$ such that:
\begin{enumerate}
\item
 $D:\widetilde N\rightarrow X$ is a local isometry;
\item
 $\rho:\pi_1(N, x_0)\rightarrow G$ is a representation;
\item
 the map $D$ is $\pi_1(N, x_0)$-equivariant in the following sense
\[
      D(\gamma(x))=\rho(\gamma)D(x)\qquad\textrm{for all
      }\gamma\in\pi_1(N,x_0)\textrm{ and }x\in\widetilde N.
\]
\end{enumerate}
Moreover given two such pairs $(D,\rho)$ and $(D',\rho')$ there exists a unique $g\in
G$ such that
\[
   D'=gD\qquad\textrm{and}\qquad \rho'=g\rho g^{-1}.
\]
\end{prop}
\cvd
\begin{defi}
With the notation of the proposition $D$ is called a \textbf{developing map} of
$N$ and $h$ is the holonomy representative compatible with $D$. The
coniugacy class of $h$ is called the \textbf{holonomy} of the $(X,G)$-manifold
$N$.
\end{defi}
\begin{oss}
\begin{enumerate}
\item\emph{
Generally $D$ is only a local isometry neither injective nor surjective.}
\item\emph{
If $D$ is a global isometry between $\widetilde N$ and $X$, we say that the
$(X,G)$-manifold $N$ is} complete. \emph{A hyperbolic manifold is} complete
\emph{ as a Riemannian manifold iff it is $(\Kn{H}{n},\SOO^+(n,1))$-complete.
If $N$ is complete then $\rho$ is a faithful representation and the image
$\Gamma$ acts freely and properly discontinuously on $X$. The isometry $D$
induces an isometry $\hat D:\widetilde N/\pi_1(N,x_0)\rightarrow X/\Gamma$.}
\end{enumerate}
\end{oss}
Let $M:=\Kn{H}{n}/\Gamma$ be a complete hyperbolic $n$-manifold. 
 Notice that $\Gamma$ acts freely and properly discontinuously on the whole
$\fut(0)$. The \emph{future complete Minkowskian cone} on $M$ is the flat Lorentz spacetime 
$\mathcal{C}^+(M):=\fut(0)/\Gamma$. Notice that $\fut(0)$ has regular cosmological time $\widetilde
 T:\fut(0)\rightarrow\mathbb R_+$  which is in fact a real analytic
 submersion with level surfaces
 $\mathbb{H}_a=\{x\in\fut(0)|-x_0^2+x_1^2+\ldots+x_n^2 =-a^2\}$.
For every $p\in\fut(0)$ we have $\widetilde T(p)=d(p,0)$ and the origin is the unique
 point with this property.
Every
$\mathbb{H}_a$ is a Cauchy surface of $\fut(0)$ so that it is globally
 hyperbolic.\par
Since $\widetilde T$ is $\Gamma$-invariant it induces the cosmological time
 $T:\mathcal C^+(M)\rightarrow\mathbb R_+$ with level surface $\widetilde S_a=\mathbb
 H_a/\Gamma$.
Notice that $M=\widetilde S_1$ so that $C^+(M)$ is diffeomorphic to $\mathbb R_+\times M$.
\par
We are interested in studying  the  globally hyperbolic flat spacetimes $Y$
which admit a Cauchy surface diffeomorphic to $M$
(hence $Y$ is diffeomorphic to $\mathbb R_+\times M$). We shall provide a
complete discussion of this problem under the assumption that $M$ is
\emph{compact}. So from now on $M:=\Kn{H}{n}/\Gamma$ is a compact hyperbolic manifold.\par 
The set of globally hyperbolic flat spacetime structures on $\mathbb R_+\times
M$ has a natural topology (induced by the $\mathrm C^\infty$-topology on
symmetric forms).
Let us denote by $Lor(M)$  this space.
We know that $\Diffeo(\mathbb{R}_+\times M)$ acts continuously on
$Lor(M)$. The quotient  $\mathcal{M}_{Lor}(M)$ is called the \emph{moduli
  space}, whereas the \emph{Teichm\"uller space} $\mathcal{T}_{Lor}(M)$ is the
quotient of the group by the action of $Lor(M)$  of homotopically trivial
diffeomorphisms.
Notice that two structures which differ by a homotopically trivial
diffeomorphism give the same holonomy (up to coniugacy), so that the holonomy
depends only on the class of the structure in the Teichm\"uller space. \par
For every group $G$ denote by $\mathcal R_G$ the set of representations 
\[
   \pi_1(M,x_0)\rightarrow G
\]
up to coniugacy.
As $\pi_1(M,x_0)\cong\pi_1(\mathbb R_+\times M,x_0)$ we have a continuous
holonomy map 
\[
   \rho:\mathcal T_{Lor}(M)\rightarrow\mathcal R_{\ISO(\Kn{M}{n+1})}
\]
with linear part
\[
   \mathrm d\rho: \mathcal T_{Lor}(M)\rightarrow\mathcal R_{\SOO^+(n,1)}
\]
In a recent paper \cite{Andersson2} it is shown that every linear holonomy is
faithful with discrete image (for $n=2$ this was deduced in \cite{Mess} as a
corollary of a theorem of Goldman).
So we shall often confuse the linear holonomy with its image subgroup of
$\SOO^+(n,1)$ (up to coniugacy). If $n\geq 3$ Mostow rigidity theorem implies
that the linear holonomy group coincides with $\Gamma$ (up to coniugacy).
Thus if $n\geq 3$  the image of the holonomy map $h:\mathcal
T_{Lor}(M)\rightarrow\mathcal R_{\ISO(\Kn{M}{n+1})}$ is contained in
 \[
   \mathcal{R}(\Gamma)=\{[\rho]\in \mathcal
   R_{\ISO(\Kn{M}{n+1})}|\,
   \mathrm{d}(\rho(\gamma))(0)=\gamma\,\textrm{ for
   all}\,\gamma\in\Gamma\}.
\]
When $n=2$ we have to vary the hyperbolic structure on $M$ (i.e. the group
 $\Gamma$) which is now a closed surface of genus $g\geq 2$.
Anyway $\mathcal R(\Gamma)$ is the key object to be understood.\par
Let $\rho$ a representation of $\Gamma$ into $\ISO(\Kn{M}{n+1})$ whose linear part
is the identity. Thus $\rho(\gamma)=\gamma+t_\gamma$ where $t_\gamma=\rho(\gamma)(0)$ is the traslation
part.
By imposing the homomorphism condition we obtain 
\[
    t_{\alpha\beta}=t_\alpha+\alpha t_\beta\qquad\forall\alpha,\beta\in\Gamma.
\]
So that $(t_\gamma)_{\gamma\in\Gamma}$ is a cocycle in $Z^1(\Gamma,\Kn{R}{n+1})$.
Conversely if $(t_\gamma)_{\gamma\in\Gamma}$ is a cocyle then the map
$\Gamma\ni \gamma\mapsto \gamma+t_\gamma\in\ISO(\Kn{M}{n+1})$ is a homomorphism.
Hence the homomorphisms of $\Gamma$ into $\ISO(\Kn{M}{n+1})$ whose linear part
is the identity are parametrized by the cocycles in $Z^1(\Gamma,\Kn{R}{n+1})$.\par
Take two such representations $\rho$ and $\rho'$ and let
$(t_\gamma)_{\gamma\in\Gamma}$ and $(t'_\gamma)_{\gamma\in\Gamma}$ be the respective traslation parts.
Suppose now that $\rho$ and $\rho'$ are conjugated by
some element $f\in\ISO(\Kn{M}{n+1})$. Then we have that the linear part of $f$
commutes with the elements of $\Gamma$. Since the centralizer of
$\Gamma$ in $\SOO^+(n,1)$ is trivial  $f$ is a pure traslation of a
vector $v=f(0)$. Now by imposing the condition $\rho'(\gamma)=f\rho(\gamma) f^{-1}$ we
obtain that $t_\gamma-t'_\gamma=\gamma v-v$ so that $t_\gamma$ and $t'_\gamma$
differ by a coboundary.
Conversely if $(t_\gamma)_{\gamma\in\Gamma}$ and $(t'_\gamma)_{\gamma\in\Gamma}$ are cocycles
which differ by a coboundary then they induce representations which are
conjugated.
Hence there is a natural identification between $\mathcal{R}(\Gamma)$ and 
the cohomology group $\coom1(\Gamma,\Kn{R}{n+1})$. 
In what follows we use this identification without mentioning it. In
particular for a cocycle $\tau$ we denote by $\rho_\tau$ and $\Gamma_\tau$
respectively  the
homomorphism corrisponding to $\tau$ and its image. 
\paragraph{Main results.}
We can state now the main results of this paper. 
\begin{teo}\label{main}
For every $[\tau]\in\coom1(\Gamma,\Kn{R}{n+1})$ there is a unique $[Y_\tau]\in\mathcal
T_{Lor}(M)$ represented by a maximal globally hyperbolic future complete
spacetime $Y_\tau$ that admits a pair $(D,\rho)$ of compatible developing map
\[
      D:\widetilde Y_\tau\rightarrow\Kn{M}{n+1}
\]
and holonomy representative
\[
    \rho:\pi_1(Y_\tau)\cong\pi_1(M)\rightarrow\ISO(\Kn{M}{n+1})
\]
such that
\begin{enumerate}
\item
$\rho=\rho_\tau$.
\item
$D$ is injective and so it is an isometry onto its image $\mathcal D_\tau$
which is a proper convex domain of $\Kn{M}{n+1}$. Moreover it is a
\emph{future set} (i.e. $\mathcal D_\tau=\fut(\mathcal D_\tau)$).
\item
The action of $\pi_1(M)$ on $\mathcal D_\tau$ via $\rho$ is free and properly
discontinuous so that the developing map $D$ induces an isometry between
$Y_\tau$ and $\mathcal D_\tau/\pi_1(M)$. 
\item
The spacetime $\mathcal D_\tau$ has a \textbf{canonical cosmological time} $\widetilde
T:\mathcal D_\tau\rightarrow\mathbb R_+$ which is a $\mathrm C^1$-submersion.
Every level surface $\widetilde S_a$ is the graph of a proper $\mathrm C^1$-convex
function defined over the horizontal hyperplane $\{x_0=0\}$.
\item
The map $\widetilde T$ is $\pi_1(M)$-equivariant and induces the canonical cosmological
time $T:Y_\tau\rightarrow\mathbb R_+$; this is a proper $\mathrm
C^1$-submersion and every level surface $S_a=\widetilde S_a/\pi_1(M)$ is $\mathrm
C^1$-diffeomorphic to $M$.
\item
For every $p\in\mathcal D_\tau$ there exists a unique
$r(p)\in\pass(p)\cap\partial\mathcal D_\tau$ such that $\widetilde
T(p)=d(p,r(p))$. The map $r:\mathcal D_\tau\rightarrow\partial \mathcal D_\tau$ is
continuous. The image $\Sigma_\tau:=r(\mathcal D_\tau)$ is said the
\textbf{singularity in the past}. $\Sigma_\tau$ is spacelike-arc connected and
contractile and $\pi_1(M)$-invariant. Moreover the map $r$ is $\pi_1(M)$-equivariant.
\end{enumerate}
The map
\[
   \mathcal R(\Gamma)\ni[\rho_\tau]\mapsto[Y_\tau]\in\mathcal T_{Lor}(M)
\]
is continuous section of the holonomy map.\par
The same statement holds by replacing ``future'' with ``past''. Let us call
$Y^-_\tau$ and $\mathcal D^-_\tau$ the corresponding spaces.\par
Every globally hyperbolic flat spacetime with compact spacelike Cauchy surface
and holonomy group equal to $\rho_\tau(\pi_1(M))$ is
diffeomorphic to $M\times\mathbb R_+$ and  embeds isometrically either into $Y_\tau$
or into $Y^-_\tau$.
\end{teo} 
\cvd
In the section \ref{Convergence-sec} we shall study the metric properties of
the surfaces $\widetilde S_a$. We look at the asymptotic behaviour
of the metrics properties of the action of $\Gamma$ on $\widetilde S_a$ for
$a\rightarrow+\infty$ and for $a\rightarrow 0$. In particular we focus
on the \textbf{Gromov convergence} when $a\rightarrow+\infty$ and on the
convergence of the \textbf{marked lenght spectrum} when $a\rightarrow 0$ (the
definition of these concepts are given in section \ref{Convergence-sec}). 
The principal result that we get is the following.
\begin{teo}\label{main2}
Let $\tau\in Z^1(\Gamma, \Kn{R}{n+1})$ and $\dom\subset\Kn{M}{n+1}$ be the 
universal cover of $Y_\tau$.
Let $\widetilde S_a$ be the CT level surface
$\widetilde T^{-1}(a)$ of $\dom$ and let $d_a$ be the natural distance on $\widetilde
S_a$. We have that $\widetilde S_a$ is a $\pi_1(M)$-invariant spacelike surface
and $\pi_1(M)$ acts by isometries on it.\par
When $a\rightarrow +\infty$ the $\pi_1(M)$-action
on the rescaled surface $(\widetilde S_a, \frac{d_a}{a})$ converges in the Gromov sense
to the action of $\,\pi_1(M)$ on $\Kn{H}{n}$.\par
When $a\rightarrow 0$ the\emph{ marked lenght spectrum} of the
$\pi_1(M)$-action on $(\widetilde S_a, d_a)$ converges to the spectrum of the
$\pi_1(M)$-action on the singularity in the past $\Sigma$.
\end{teo}
\cvd
Now we point out some comments and corollaries.\\
The following statement is an immediate consequence of theorem \ref{main}.
\begin{corollario}
Let $F$ be a $n$-manifold and suppose that there exists a Lorentzian flat
structure on $\mathbb R\times F$ such that $\{0\}\times F$ is a spacelike
surface.
Suppose that the holonomy group for such structure is $\Gamma_\tau$ for some
$\tau\in Z^1(\Gamma,\Kn{R}{n+1})$. Then $F$ is diffeomorphic to
$M=\Kn{H}{n}/\Gamma$.
\end{corollario}
\cvd
By studying  the action of $\Gamma_\tau=\rho_\tau(\pi_1(M))$ on the boundary $\partial \mathcal
D_\tau$ we shall show that $\Gamma_\tau$ \emph{does not act freely and properly
discontinuously on the whole $\Kn{M}{n+1}$}.\\

On the domain $\mathcal D_\tau$ there is a natural field $-N$ which is the
Lorentzian gradient of the cosmological time $\widetilde T$. We have that $N$ is a timelike
vector, furthermore $N$ is future directed by the choice of the sign. Notice that $N(x)$ is the normal
vector to $\widetilde S_{\widetilde T(x)}$ at $x$, so that we call it the \textbf{normal
  field}. Under the identification of $T_x\Kn{M}{n+1}$ with $\Kn{M}{n+1}$ it
results that $N(x)\in\Kn{H}{n}$ (in fact it is also called the \textbf{Gauss
  map} of the surfaces $\widetilde S_a$). 
The restriction $N|_{\widetilde S_a}$ is a surjective and proper map.
The map $N$ is $\pi_1(M)$-equivariant so
that it induces a map $\overline N:\mathcal
D_\tau/\Gamma_\tau\rightarrow\Kn{H}{n}/\Gamma$. For all $a>0$ the restriction the
of the map $\overline N|_{\widetilde S_a}$ has degree $1$.\\

When $n=2$ it turns out that the singularity $\Sigma_\tau$ is a \emph{real
  tree}. Moreover Mess showed that the $N$-images of the fibers of the
  retraction $r$ produce  a $\Gamma$-invariant \emph{geodesic lamination} $L$ of
  $\Kn{H}{2}$. According to Skora's theorem \cite{Skora} $L$ is the geodesic
  lamination dual to the real tree $\Sigma_\tau$. More precisely, the complete
  duality is realized by suitably equipped $L$ with a transverse measure
  $\mu$. Finally the triple $(M,L,\mu)$ determines the spacetime $Y_\tau$.\\

In section \ref{CT-section} we shall see that for $n\geq 2$ the $N$-images of
the fibers of the retraction $r$ determine a \emph{geodesic stratification} of
$\Kn{H}{n}$ (we introduce this notion in section \ref{CT-section} and we prove
that in dimension $n=2$ the geodesic stratifications are the geodesic
laminations).
Moreover in section \ref{measure-section} we introduce the notion of
\emph{measured geodesic stratification} and we show that every measured
geodesic stratification enables us to construct a spacetime $Y_\tau$. For some
technical reasons we are not able, for the moment, to show that this
corrispondence is bijective. However we present an interesting class of
spacetimes: the spacetimes with \emph{simplicial singularity}. In
dimension $n=2$ these spacetimes have singularity which is a
\emph{simplicial tree} and the corresponding geodesic lamination is a
multicurve.
We shall show that for spacetimes with simplicial singularity the complete
duality between singularity and geodesic stratification is realized in very
explicit way.


\section{Construction of $\mathcal D_\tau$}\label{construction-sec}

Let $\Gamma$ be a free-torsion co-compact discrete subgroup of $\SOO^+(n,1)$
and $M:=\Kn{H}{n}/\Gamma$.
Fix $[\tau]\in\coom1(\Gamma,\Kn{R}{n+1})$ and let $\Gamma_\tau$ be the image of
the homomorphism associated with $\tau$. Moreover for every $\gamma\in\Gamma$
we shall denote by $\gamma_\tau$ the affine trasformation $x\mapsto \gamma(x)+\tau_\gamma$.
In this section we construct a $\Gamma_\tau$-invariant future complete convex
domain of $\Kn{M}{n+1}$. Moreover we show that
the action of $\Gamma_\tau$ on this domain is free and properly discontinuous
and the quotient is diffeomorphic to $\mathbb R_+\times M$.\\

First let us show that there is a $\mathrm C^\infty$-embedded hypersurface $\Ftilde$ of $\Kn{M}{n+1}$ which
is spacelike (i.e. $T_p\Ftilde$  is a spacelike subspace of $T_p\Kn{M}{n+1}$)
and $\Gamma_\tau$-invariant such that the quotient $\Ftilde/\Gamma_\tau$ is
diffeomorphic to $M$. 
We start with an easy and useful lemma (see \cite{Mess}).
\begin{lemma}\label{spacelike manifolds are graph}
Let $S$ be a manifold and $f:S\rightarrow\Kn{M}{n+1}$ be an $\mathrm
C^r$-immersion $(r\geq 1)$
such that $f^*\eta$ is a complete Riemannian  metric on $S$ . 
Then $f$ is an embedding.
Moreover fix orthonormal affine coordinates $(y_0,\ldots,y_n)$.
 Then $f(S)$ is a graph over the horizontal plane $\{y_0=0\}$.
\end{lemma}
\Dim
In coordinates put $f(s)=(i_0(s),\ldots,i_n(s))$.
Let $\pi:S\rightarrow\{y_0=0\}$ be the canonical projection ($\pi(s)=(0,i_1(s),\ldots,i_n(s))$) 
We have to show that $\pi$ is a $\mathrm{C}^r$-diffeomorphism.\par
Notice that $\pi$ is a $\mathrm{C}^r$-map between Riemannian manifold. We
claim that $\pi$ is distance increasing i.e.
\begin{equation}\label{distance increasing}
     \E{\mathrm d\pi(x)[v]}{\mathrm d\pi(x)[v]}\geq \left(f^*\eta\right)(v,v).
\end{equation}   
The lemma follows easily from the claim: in fact the equation (\ref{distance
  increasing}) implies that $\pi$ is a local
  $\mathrm{C}^r$-diffeomorphism. Furthermore a classical argument shows that
  $\pi$ is path-lifting and so $\pi$ is a covering map. Since the horizontal
  plane is simply connected it follows that $\pi$ is a
  $\mathrm{C}^r$-diffeomorphism.\par
Let us  prove tha claim. Let $v\in T_sS$ and let  $\mathrm d
  f(s)[v]=(v_0,\ldots,v_n)$ then we have 
 $\mathrm d\pi(s)[v]=(0,v_1,\ldots,v_n)$. Thus
\[
     f^*\eta(v,v)=\E{\mathrm d\pi(x)[v]}{\mathrm d\pi(x)[v]}-v_0^2.
\]
\cvd
\begin{oss}\emph{
Let $S$ be a $\Gamma_\tau$-invariant spacelike hypersurface in $\Kn{M}{n+1}$
such that the action of $\Gamma_\tau$ on it is free and properly discontinuous.
Suppose that $S/\Gamma_\tau$ is compact. By the Hopf-Rinow theorem we know that $S$ is
complete and so the previous lemma applies.
}\end{oss} 
Now we  want to construct a
$\Gamma_\tau$-invariant spacelike hypersurface $\Ftilde$. 
In fact we shall construct $\Ftilde$ in a particular class of spacelike
hypersurfaces.\par
\begin{defi}
A closed connected spacelike  hypersurface $S$  divides $\Kn{M}{n+1}$ in two
components, the future and the past of $S$.
We say that $S$ is \emph{future convex} (resp. \emph{past convex}) if $\fut(S)$
(resp. $\pass(S)$) is a convex set and $S=\partial\fut(S)$
(resp. $S=\partial\pass(S)$). 
Moreover $S$ is \emph{future strictly convex} 
(resp. \emph{past strictly convex}) if $\fut(S)$ (resp. $\pass(S)$) is
strictly convex.
\end{defi}
\begin{oss}\emph{
The hyperbolic space $\Kn{H}{n}\subset\Kn{M}{n+1}$ is an example of spacelike
future strictly convex hypersurface in $\Kn{M}{n+1}$.
}\end{oss}
Let $N_0$ be the flat Lorentz structure on $[1/2,3/2]\times M$ given by the
standard inclusion $[1/2,3/2]\times M\subset\mathcal{C}(M)$ (where
$\mathcal{C}(M))$ is the Minkowskian cone on $M$). Let 
$\widetilde N_0=\{x\in\Kn{M}{n+1}| x\in\fut(0)\textrm{ and } d(0,x)\in[1/2,3/2]\}$
be the universal covering of $N_0$.
The following theorem was stated by Mess (\cite{Mess}) for the case $n=2$. However his proof 
runs in all dimensions. We relate it here for sake of completeness.
\begin{teo}\label{existence of structure}
 If $U$ is a bounded neighborhood of $0$ in $Z^1(\Gamma, \Kn{R}{n+1})$ there exists
$K>0$ and a $\mathrm C^\infty$-map 
\[
     dev:U\times\widetilde N_0\rightarrow\Kn{M}{n+1}
\]
such that\\
1. for every $\sigma\in U$ the map 
\[
   dev_\sigma:\widetilde N_0\ni x\mapsto dev(\sigma,x)\in\Kn{M}{n+1}
\]
 is a developing map whose holonomy is the representation associated  with $\sigma$;\\
2. $dev_0$ is the multiplication by $K$;\\
3. $dev_\sigma(\Kn{H}{n})$ is a strictly convex spacelike
hypersurface. 
\end{teo}
\Dim
For a Thurston theorem (see \cite{Epstein}) there exists a neighborhood $U_0$ of $0$ in
$Z^1(\Gamma, \Kn{R}{n+1})$ and a $\mathrm C^\infty$-map
\[
   dev': U_0\times \widetilde N_0\rightarrow\Kn{M}{n+1}
\]
such that\\ 
1. for every $\sigma\in U_0$ the map $dev'_\sigma:\widetilde
N_0\rightarrow\Kn{M}{n+1}$ is a developing map whose holonomy
 is the representation associated  with $\sigma$;\\
2. $dev'_0$ is the identity.\par
By using the compactness of $M$ it is easy 
to show that if $U_0$ is choosen
sufficiently small then $dev_\sigma(\Kn{H}{n})$ is a spacelike future convex 
hypersurface.\par
Now fix $K>0$ so that $K\cdot U_0\supset U$. Let us define
$dev: U\times\widetilde N_0\rightarrow\Kn{M}{n+1}$ by the rule
\begin{equation}\label{developing}
    dev(\sigma,x):=Kdev'(\sigma/K,x)
\end{equation}
It is staightforward to see that $dev_\sigma$ is a developing map whose
holonomy is the representation associated with $\sigma$. Clearly $dev_0$ is the
multiplication by $K$ and
$dev_\sigma(\Kn{H}{n})=K\cdot dev'_{\sigma/K}(\Kn{H}{n})$ is a future convex
spacelike surface invariant for $\Gamma_\sigma$.
\cvd
Now fix a bounded neighborhood of $0$ in $Z^1(\Gamma, \Kn{R}{n+1})$ which
contains the cocycle $\tau$. Consider the map $dev$ of the previous theorem
and let $\Ftilde$ be the hypersurface $dev_\tau(\Kn{H}{n})$. Then $\Ftilde$
is a $\Gamma_\tau$-invariant future stricly convex spacelike hypersurface 
such that the $\Gamma_\tau$-action on it is \emph{free and
properly discontinuous}  and $\Ftilde/\Gamma_\tau\cong M$.
Clearly in the same way we can obtain a $\Gamma_\tau$-invariant spacelike
hypersurface  $\Ftilde^-$ which is\emph{ past} strictly convex such that
$\Ftilde^-/\Gamma_\tau\cong M$. \\

Now given a hypersurface $\widetilde F$ we construct a
natural domain $D(\widetilde F)$ which includes it. Furthermore we show that if
$\widetilde F$ is $\Gamma_\tau$-invariant and the action on $\widetilde F$ is free and
properly discontinuous then the same holds for $D(\widetilde F)$.
\begin{defi}
Given a spacelike hypersurface $\widetilde F$ the \textbf{domain of dependence} of $\widetilde
F$ is the set $D(\widetilde F)$ of the points
$p\in\Kn{M}{n+1}$ such that all inextendable causal curves which pass through $p$
 intersect $\widetilde F$.
\end{defi}
 The following is a well known result (for istance see \cite{Geroch}).
\begin{prop}\label{description of domain of dependence}
The domain of dependence $D(\widetilde F)$ is open.
Moreover if $\widetilde F$ is complete (i.e. the natural Riemannian structure on
it is complete) then a point $p\in\Kn{M}{n+1}$ lies in $D(\widetilde  F)$ if and only if
each null line which passes through $p$ intersects $\widetilde F$.
\end{prop}
\cvd
\begin{prop}\label{points not in domain}
Let $\widetilde F$ be a complete spacelike $\mathrm C^1$-hypersurface.
Let $p\notin D(\widetilde F)$ and $v$ be a null vector such that the line 
$p+\mathbb R v$  does not intersect $\widetilde F$. Then the null plane 
$P=p+\ort{v}$ does not intersect $\widetilde F$.
\end{prop}
\Dim
Suppose that $S:=\widetilde F\cap P$ is not empty. Since this intersection is
transverse it follows that it is a closed $(n-1)$-submanifold of $\widetilde F$ and so
it is complete.\par
Fix a set of orthonormal affine coordinates $(y_0,\ldots,y_n)$ such that $p$
is the origin and $P=\{y_0=y_1\}$ (i.e. $v=(1,1,0,\ldots,0)$). Consider the map
\[
    \pi:S\ni(y_0,y_1,\ldots,y_n)\rightarrow (0,0,y_2,\ldots,y_n)\in\{y_0=y_1=0\}.
\]  
As well as in the proof of lemma \ref{spacelike manifolds are graph} we argue
that $\pi$ is a diffeomorphism. Thus there exists $s\in\Kn{R}{}$ such that 
$q=(s,s,0,\ldots,0)\in\widetilde F$. But $q$ lies on the line $p+\mathbb R v$ and
this is a contradiction.
\cvd
\begin{corollario}\label{convexity of domain of dependence}
Let $\widetilde F$ be a complete spacelike hypersurface.
The domain of dependence $D(\widetilde F)$  is a convex set. Moreover for every $p\notin
D(\widetilde F)$ there exists a null  plane through $p$ which is a support plane for
$D(\widetilde F)$.\par 
Suppose that $\widetilde F$ is \emph{$\Gamma_\tau$-invariant} and  $D(\widetilde F)$
is not the whole $\Kn{M}{n+1}$ 
then either
\begin{eqnarray*}
      D(\widetilde F)=\bigcap_{\begin{array}{l}
                          {\scriptstyle P\textrm{ null  plane}}\\
                          {\scriptstyle P\cap\widetilde F=\varnothing}
                          \end{array}} \fut(P)\qquad  \textrm{  or}\\ 
      D(\widetilde F)=\bigcap_{\begin{array}{l}
                          {\scriptstyle P\textrm{ null  plane}}\\
                          {\scriptstyle P\cap\widetilde F=\varnothing}
                          \end{array}} \pass(P).\qquad\quad
\end{eqnarray*}  
 Thus $D(\widetilde F)$ is a future or past set (i.e. $D(\widetilde
 F)=\fut(D(\widetilde F))$ or $D(\widetilde F)=\pass(D(\widetilde F))$).
 \end{corollario}
\Dim
\begin{figure}
\begin{center}
\input{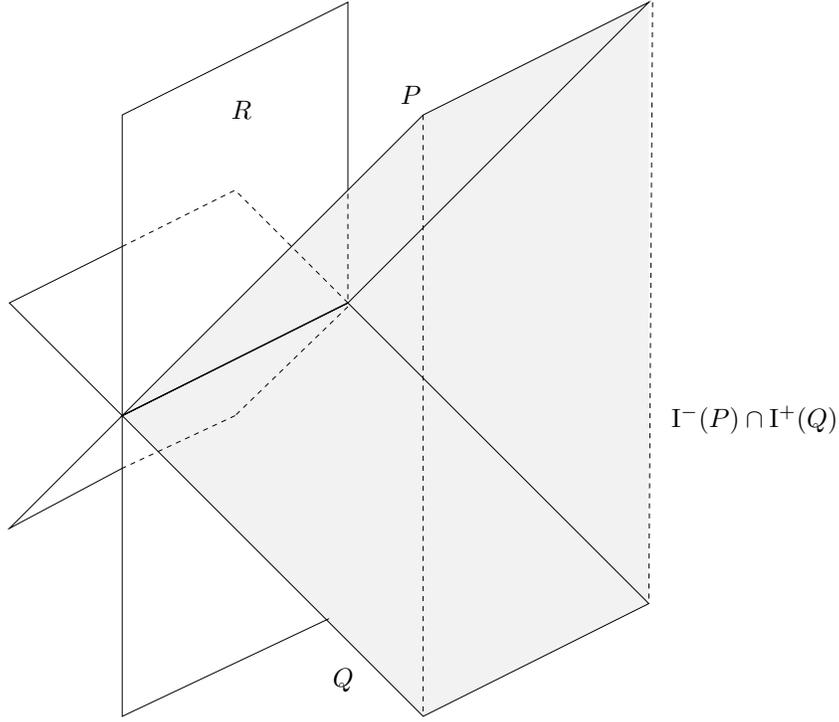}
\caption{{\small If $P$ and $Q$ are not parallel there exists a timelike
    support plane}}\label{nullsupp.fig}
\end{center}
\end{figure}
The proposition \ref{points not in domain} implies that for every $p\notin D(\widetilde
F)$ there exists a null plane $P$ which passes through $p$ and does not intersect
$D(\widetilde F)$. Thus $D(\widetilde F)$ is a convex set.\par
Suppose now that $\widetilde F$ is $\Gamma_\tau$-invariant and suppose that
$D(\widetilde F)$ is not the whole $\Kn{M}{n+1}$. We have to show that either
$D(\widetilde F)$ is contained in the future of the null support planes or it is
contained in the past of the null support planes.
Suppose by contradiction that there exist null support planes $P$
and $Q$ such that $D(\Ftilde)\subset\pass(P)\cap\fut(Q)$. 
First suppose that $P$ and $Q$ are not parallel. 
Then there exists a timelike support plane $R$ (see fig.\ref{nullsupp.fig}). 
Fix now an affine coordinates system $(y_0,\ldots,y_n)$ such that $R=\{y_n=0\}$. By
lemma \ref{spacelike manifolds are graph} we know that $\widetilde F$ is a graph of a function
defined on the horizontal plane $\{y_0=0\}$. Then $\widetilde F\cap
R\neq\varnothing$ and this is a contradiction.\par
Suppose now that we cannot choose non-parallel $P$ and $Q$. Then it follows
that null tangent planes are all parallel. Thus let $v$ be the null vector
orthogonal to all null tangent planes and $[v]$ the corresponding point on
$\partial\Kn{H}{n}$. 
Since we have that $\Gamma_\tau$ acts
on $\widetilde F$ we have that $\Gamma_\tau$ permutes the null support planes of
$D(\widetilde F)$. It follows that $\Gamma\cdot[v]=[v]$. But $\Gamma$
is a discrete co-compact group and so it does not fix any point in
$\overline{\mathbb H}^n$. 
\cvd 
\begin{oss}\emph{
The completeness of $\widetilde F$ is an essential hypothesis. For istance the
domain of dependence of $\widetilde F-{p}$ is not convex. 
}\end{oss}
\begin{oss}\emph{
If $\widetilde F$ is future (past) convex then $D(\widetilde F)$ is future (past)
complete.
}\end{oss}
\begin{prop}\label{action on domain is properly discontinuous}
Let $\widetilde F$ be a $\Gamma_\tau$-invariant spacelike hypersurface such that the
$\Gamma_\tau$-action on it is free and properly discontinuous. 
Then  $\Gamma_\tau$ acts freely and properly discontinuously on the whole
$D(\widetilde F)$. Moreover $D(\widetilde F)/\Gamma_\tau$ is diffeomorphic to 
$\mathbb R_+\times\widetilde F/\Gamma_\tau$.
\end{prop}
\Dim
Since $\Gamma_\tau$ is torsion free it is sufficient to show that the action
is properly discontinuous.\par
Let $K,H\subset D(\widetilde F)$ be compact sets. We have to show that the set
$\Gamma(K,H)=\{\gamma\in\Gamma|\gamma_\tau(K)\cap H\neq\varnothing\}$ is
finite.
By using proposition \ref{description of domain of dependence} we get that
the sets $C=(\cfut(K)\cup\cpass(K))\cap\widetilde F$ and
$D=(\cfut(H)\cup\cpass(H))\cap\widetilde F$ are compact. Furthermore
$\gamma_\tau (C)=(\cfut(\gamma_\tau(K))\cup\cpass(\gamma_\tau(K)))\cap\widetilde F$.
Thus $\Gamma(K,H)$ is contained in $\Gamma(C,D)$. Since the action of
$\Gamma_\tau$ on $\widetilde F$ is  properly discontinuous  we have that
$\Gamma(C,D)$ is finite.\par
Notice that $\widetilde F/\Gamma_\tau$ is a Cauchy surface in $D(\widetilde
F)/\Gamma_\tau$ and thus $D(\widetilde F)/\Gamma_\tau\cong\mathbb R_+\times\widetilde F/\Gamma_\tau$.
\cvd
We have constructed a $\Gamma_\tau$-invariant future convex hypersurface
$\Ftilde$ such that $\Ftilde/\Gamma_\tau\cong M$. Now consider $D(\Ftilde)$.
We know that it is a $\Gamma_\tau$-invariant future complete convex set and
$D(\Ftilde)/\Gamma_\tau\cong \mathbb R_+\times M$. \emph{From now on we denote
$D(\Ftilde)$ by $\dom$}.\par
Notice we can also consider the domain of dependence of $\Ftilde^-$. In the
same way we can 
show that $D(\Ftilde^-)$ is $\Gamma_\tau$-invariant, past complete convex
domain of $\Kn{M}{n+1}$ and $D(\Ftilde^-)/\Gamma_\tau\cong \mathbb R_+\times
M$. \emph{We shall denote it by $\dompass$}.

In the remaining part of this section we prove that $\dom$ is not the whole
$\Kn{M}{n+1}$. This is a necessary condition for $\dom$ to have regular
cosmological time. In the next section we shall see that this condition is in
fact sufficient.
In order to prove that $\dom$ is a proper subset we need some geometric properties of the
$\Gamma_\tau$-invariant  future convex sets.
\begin{lemma}
Let $\Omega$ be a proper convex set of $\Kn{M}{n+1}$. Fix a set of orthonormal
affine coordinates $(y_0,\ldots,y_n)$ then $\Omega$  is a future convex set if and only
if $\partial\Omega$ is the graph on the horizontal plane $\{y_0=0\}$ of a
$1$-Lipschitz convex function.
\end{lemma}
\Dim
The \emph{if} part is quite evident. Hence suppose that $\Omega$ is a proper
future convex set. First let us show that the projection on the horizontal
plane $\pi:\partial\Omega\rightarrow\{y_0=0\}$ is a homeomorphism. Since $\partial\Omega$ is a
topological manifold it is sufficient to show that the projection is
bijective.
Since $\Omega$ is a future set we have that two points on $\partial\Omega$ are
not chronological related and so the projection is injective.\par
It remains to show that given $(a_1,\ldots,a_n)$ there exists $a_0$ such that
$(a_0,a_1,\ldots,a_n)\in\partial\Omega$.
Fix $p\in\partial\Omega$.
It is easy to see that there exist $a_+$ and $a_-$ such that
$(a_+,a_1,\ldots,a_n)\in\fut(p)$ and $(a_-,a_1,\ldots,a_n)\in\pass(p)$.  Since $\fut(p)\subset\Omega$ and
$\pass(p)\cap\Omega=\varnothing$
there exists $a_0$ such that $(a_0,\ldots,a_n)\in\partial\Omega$.\par
It follows that $\partial\Omega$ is  the graph of a function $f$. Since $\Omega$
is future convex we have that $f$ is convex. Since two points on
$\partial\Omega$ are not chronological related we have that $f$ is
$1$-Lipschitz.
\cvd
\begin{lemma}\label{support planes}
Let $\Omega$ be a $\Gamma_\tau$-invariant proper future convex set. Then for every
$u\in\Kn{H}{n}$ there exists a plane $P=p+\ort{u}$ such that $\Omega\subset
\fut(P)$.
\end{lemma}
\Dim
Consider the set $K$ of the vectors $v$ in $\Kn{M}{n+1}$ which are orthogonal
to some support plane of $\Omega$. Clarly $K$ is a convex cone with apex in $0$. Since $\Omega$ is
future complete it is easy to see that the vectors in $K$ are not
spacelike. So that the projection $\mathbb PK$ of $K$ in $\Kn{P}{n}$ is a convex subset of
$\overline{\mathbb H}^n$. Since $\Omega$ is $\Gamma_\tau$-invariant $K$ is
$\Gamma$-invariant. 
Then $\mathbb P K$ is a $\Gamma$-invariant convex set of $\overline{\mathbb H}^n$. Since it is not
empty (there exists at least a support plane) and $\Gamma$ is co-compact then $K$
contains the whole $\Kn{H}{n}$ and the lemma follows.
\cvd
\begin{lemma}\label{timelike coordinate is proper}
Let $\Omega$ be as in the last lemma. For each timelike vector $v$ the
function
\[
    \partial\Omega\ni x\mapsto \E{x}{v}\in\mathbb R
\]
is proper. If $v$ is future directed
$\lim_{x\in\partial\Omega\,\,x\rightarrow\infty}\E{x}{v}=-\infty$ (we mean by
$\infty$ the point of the Alexandroff compactification of $\partial\Omega$).\par
Moreover there exists a unique support plane $P_v$ of $\Omega$ such that it is
orthogonal to $v$  and $P_v\cap\partial\Omega\neq\varnothing$.
\end{lemma}
\Dim
Fix a timelike vector $v$. Clearly we can suppose that $v$ is future directed and $\E{v}{v}=-1$.
Fix a set of orthonormal affine coordinates $(y_0,\ldots,y_n)$ with the origin
in $0$ and such that $\der{0}{\,}=v$. Notice that $\E{x}{v}=-y_0(x)$.
Now since $\Omega$ is future complete then the boundary $\partial\Omega$ is
graph of a convex function $f:\{y_0=0\}\rightarrow\mathbb R$.\par
 We have to show
that $f$ is proper and $\lim_{x\rightarrow\infty}f(x)=+\infty$.
Thus it is sufficient to show that the set $K_C=\{x|f(x)\leq C\}$ is compact
for every $C\in\mathbb R$. Since $f$ is convex $K_C$ is a closed convex subset of
$\{y_0=0\}$. Suppose by contradiction that it is not compact. It is easy to
see that there exists $\overline x\in\{y_0=0\}$ and a horizontal vector $w$
such that the ray 
$\{\overline x+tw|t\geq 0\}$ is
contained in $K_C$.\par
We can suppose  $\E{w}{w}=1$ so that the vector  $u=\sqrt{2}v+w$ is timelike
$\E{u}{u}=-1$.
Lemma \ref{support planes} implies that there
exists $M\in\mathbb R$ such that $\E{p}{u}\leq M$ for all
$p\in\partial\Omega$.
On the other hand let $p_t=(\overline x+tw)+f(\overline x+tw)v$. We have $p_t\in\partial\Omega$ and
\[
    \E{p_t}{u}=-\sqrt{2}f(q+tw)+\E{q+tw}{q+tw}\geq-\sqrt{2}C+\E{q+tw}{q+tw}.
\]
Since $\E{q+tw}{q+tw}\rightarrow+\infty$ we get a contradiction.
\cvd
\begin{prop}\label{domain of dependence is proper}
Let $\Omega$ be a $\Gamma_\tau$-invariant future complete convex proper subset of
$\Kn{M}{n+1}$. Then there exists a null support plane of $\Omega$.
\end{prop}
\Dim
Take $q\in\partial\Omega$ and $v\in\Kn{H}{n}$ such that $P=q+\ort{v}$ is a
support plane in $x$. Fix $\gamma\in\Gamma$ and consider the sequence of
support planes $P_k:=\gamma_\tau^k(P)$. If this sequence does not escape to the
infinity there is a subsequence which converges to a support plane $Q$. The
normal direction of $Q$ is the limit of the normal directions of $P_k$. On the
other hand the normal direction of $P_k$ is the direction of
$\gamma^k(v)$. Since in the projective space $[\gamma^k(v)]$ tends to a null
direction we have that $Q$ is a null support plane.\par
Thus we have to prove that $P_k$ does not escape to the infinity.  
Let $v_k=|\E{v}{\gamma^k v}|^{-1}\gamma^k v$. We know that $v_k$ converges to an
attractor eigenvector of $\gamma$ in $\Kn{M}{n+1}$. On the other hand we have
\[
P_k=\{x\in\Kn{M}{n+1}|\E{x}{v_k}\leq\E{\gamma_\tau^k q}{v_k}\}.
\]
Thus the
sequence $P_k$ does not escape to the infinity if and only if the coefficients
$C_k=\E{\gamma_\tau^k q}{v_k}$ are bounded. Since $\{v_k\}$ is relatively compact
 in $\Kn{M}{n+1}$ it is  sufficient to show that the
coefficients
\[
C'_k:=\E{\gamma_\tau^k q-q}{v_k}
\]
are bounded.
For $\alpha\in\Gamma$ let $z(\alpha)=\alpha_\tau(q)-q$. It is easy to see that
$z$ is a cocycle (and in fact the difference $z(\alpha)-\tau_\alpha=\alpha
q-q$ is a coboundary). Thus we have
\begin{eqnarray}\label{domain of dependence is proper-formula1}
    C'_k= \Big|\frac{\E{z(\gamma^k)}{\gamma^k v}}{\E{\gamma^k v}{v}}\Big|= 
          \Big| \frac{\E{\gamma^{-k}z(\gamma^k)}{v}}{\E{\gamma^k v}{v}}\Big|=\nonumber\\
       = \Big|\frac{\E{z(\gamma^{-k})}{v}}{\E{\gamma^k v}{v}}\Big|.
\end{eqnarray}
Now let $\lambda>1$ by the maximum eigenvalue of $\gamma$. Denote by $\|\cdot\|$ the euclidean
norm of $\Kn{R}{n+1}$. Then we have that $\|\gamma^{-1}(x)\|\leq\lambda\|x\|$
for every $x\in\Kn{R}{n+1}$.
Since $z(\gamma^{-k})=-\sum_{i=1}^{k}\gamma^{-i}(z(\gamma)$ it follows that
$\|z(\gamma^{-k})\|\leq K\lambda^k$ for some $K>0$. Thus we have
that $|\E{z(\gamma^{-k})}{v}|\leq K'\lambda^k$
(in fact we have $|\E{x}{y}|\leq\|x\|\|y\|$).\par
On the other hand let $v=x^++x^-+x^0$ where $x^+$ is eigenvector for
$\lambda$ and $x^-$ is eigenvector for $\lambda^{-1}$ and $x^0$ is in
the orthogonal of $\mathrm{Span}(x^+,x^-)$. Since $v$ is future directed
timelike vector it turns out that $x^+$ and $x^-$ are future directed. Thus we have
\[
    \E{\gamma^k
    v}{v}=(\lambda^k+\lambda^{-k})\E{x^+}{x^-}\,+\,\E{x^0}{\gamma^k x^0}.
\]
Now notice that $\ort{\mathrm{Span}(x^+,x^-)}$ is spacelike and
$\gamma$-invariant. We deduce that $\E{x^0}{\gamma^k x^0}\leq\E{x_0}{x_0}$ so
that there exists $M>0$ such that
$|\E{\gamma^kv}{v}|>M\lambda^k$.  Thus $|C'_k|\leq K'/M$ and this concludes
the proof.
\cvd
Now we can easily prove that $\dom$ is not the whole $\Kn{M}{n+1}$.
\begin{corollario}
Let $\widetilde F$ be a $\Gamma_\tau$-invariant future convex spacelike
hypersurface. Then there is a null support plane which does not intersect
$\widetilde F$. Hence $D(\widetilde F)\neq\Kn{M}{n+1}$.
\end{corollario}
\Dim
Take $\Omega=\fut(\widetilde F)$ and use  proposition \ref{domain of
  dependence is proper}.
\cvd 
In dimension $n+1=4$ there is an easier argument to prove that $\dom\neq\Kn{M}{n+1}$.
Notice that if $1$ is not an eigenvalue for some $\gamma\in\Gamma$ then
the trasformation $\gamma_\tau$ has a fixed point, namely
$z:=(\gamma-1)^{-1}(\tau_\gamma)$. Generally we say that  $\gamma\in\SOO^+(3,1)$ is
\emph{loxodromic} if $\gamma-1$ is invertible. So that if $\Gamma$ contains a
loxodromic element then $\Gamma_\tau$ does not acts freely on
$\Kn{M}{3+1}$ (and in particular $\Kn{M}{3+1}$ does not coincides with $\dom$). 
\begin{lemma}
Let $\Gamma$ be a discrete cocompact subgroup of
$\SOO^+(3,1)$.
There exists a loxodromic element $\gamma\in\Gamma$.
\end{lemma}
\Dim  
We use the identification $\SOO^+(3,1)\cong\ISO(\Kn{H}{3})\cong PSL(2,\mathbb C)$. A
hyperbolic element $\gamma\in PSL(2,\mathbb C)$ is not loxodromic if and only
if $\tr\gamma\in\mathbb R$. Hence suppose by contradiction that every hyperbolic
$\gamma\in\Gamma$ has real trace.\par
Fix  a hyperbolic element $\gamma_0\in\Gamma$. Up to coniugacy we can suppose
\[
  \gamma_0=\left[\begin{array}{cc} \lambda & 0\\ 0 &
  \lambda^{-1}\end{array}\right]
\]
with $\lambda\in\mathbb R_+$. Moreover we can suppose that $\gamma_0$ is a generator of
the stabilizer of the axis $l_0$ with endpoints $0$ and $\infty$.
Now let $\alpha$ be a generic element of $\Gamma$
\[
 \alpha=\left[\begin{array}{cc} a & b\\ c &
   d\end{array}\right].
\]
By general facts on Klenian groups we know that either $\alpha$ fix the
geodesic $l_0$ or it does not fix $0$ nor $\infty$. We deduce that either
$b=c=0$ or $bc\neq 0$. Suppose $\alpha\notin\stab(l_0)$.
By imposing $\tr\alpha\in\mathbb R$ and $\tr\left(\gamma_0\alpha\right)\in\mathbb R$ we
obtain
\begin{eqnarray*}
    a+d\in\mathbb R;\\
   \lambda a+\lambda^{-1} d\in\mathbb R.
\end{eqnarray*}
Thus $a,d\in\mathbb R$. Since $ad-bc=1$ we can write
\[
 \alpha=\left[\begin{array}{cc} A & B e^{i\theta}\\ C e^{-i\theta} &
   D\end{array}\right]
\] 
with $A,B\in\mathbb R$, $C,D\in\mathbb{R}-\{0\}$ and $\theta\in[0,\pi)$.\par
Now let $\beta\in\Gamma-\stab(l_0)$ 
\[
\beta=\left[\begin{array}{cc} A' & B' e^{i\theta'}\\ C e^{-i\theta'} &
   D'\end{array}\right].
\]
The first entry of $\alpha\beta$ is $AA'+BC'e^{i(\theta-\theta')}$. If we
impose that it is real we deduce that  $\theta=\theta'$ (notice that $BC'\neq0$).
So there exists a $\theta_0$ such that
for every $\gamma\in\Gamma-\stab(l_0)$
\[
  \gamma=\left[\begin{array}{cc} A & B e^{i\theta_0}\\ C e^{-i\theta_0} &
   D\end{array}\right]
\]
with $A,B,C,D\in\mathbb R$. Hence the rotation $R_{-\theta_0}$ conjugates
$\Gamma$ in $PSL(2,\mathbb R)$ and so $\Gamma$ is Fuchsian. But this is a contradiction.
\cvd

\section{The Cosmological Time and the  Singularity in the Past}\label{CT-section}

In this section we shall see that the cosmological time on
$\dom/\Gamma_\tau$ is regular and the level surfaces are homeomorphic to
$M$. Furthermore we study the boundary of $\dom$ and we shall
see that it determines a \emph{geodesic stratification} in $M$. If $n=2$ this
stratification is in fact the geodesic lamination which Mess associated with $\tau$.\\

We study the geometry of a general class of domain of $\Kn{M}{n+1}$, the
\textbf{regular convex domain}. We shall see that $\dom$ is
$\Gamma_\tau$-invariant regular convex domain. The most results of this
section are quite general and we do not use the action of the group
$\Gamma_\tau$. We shall see that every regular domain $\Omega$ is provided with a
regular \textbf{cosmological time} $T$, a \textbf{retraction} $r$ on the \textbf{ singularity
  in the past}, and a \textbf{normal field} $N:\Omega\rightarrow\Kn{H}{n}$ which is (up
to the sign) 
the Lorentzian gradient of $T$. Moreover if the domain is $\Gamma_\tau$-invariant then all
these objects are $\Gamma_\tau$-invariant. Finally we see that if the normal
field is surjective (this is the case for istance if the domain is $\Gamma_\tau$-invariant)
then the images in $\Kn{H}{n}$ of the fibers of the retraction give a
\emph{geodesic stratification}. This stratification is
$\Gamma$-invariant if $\Omega$ is $\Gamma_\tau$-invariant.
\begin{defi}
Let $\Omega\subset\Kn{M}{n+1}$ be a nonempty convex open set. We say that $\Omega$ is a
future complete (resp. past complete) regular convex domain if there exists a
family of null support planes $\mathcal F$ which contains at least $2$
non-parallel planes and such that
\[
       \Omega=\bigcap_{P\in\mathcal
       F}\fut(P)\qquad(\textrm{resp. }\Omega=\bigcap_{P\in\mathcal
       F}\pass(P)).
\]
\end{defi}
\begin{oss}\emph{
The condition that there exists at least $2$ non-parallel planes escludes that
$\Omega$ is the whole $\Kn{M}{n+1}$ or that $\Omega$ is the future of a null
plane. These domains have not regular cosmological times. On the other hand if
there exists at least $2$ non-parallel null support planes there exists a
spacelike support plane and we shall see that this condition assure the existence of a regular
cosmological time.
} \end{oss}
\begin{oss}\emph{
Let $\widetilde F$ be a  $\Gamma_\tau$-invariant future convex complete  spacelike
hypersurface. By corollary \ref{convexity of domain of dependence} we know that
 $D(\widetilde F)$ is the intersection of either the future or the past of its
 null support planes. By \ref{domain of dependence is proper} it results that
 $D(\widetilde F)$ is not the whole $\Kn{M}{n+1}$. Finally proposition \ref{support
   planes} assure us that there exist spacelike support planes. It follows
 that $D(\widetilde F)$ is a future complete regular domain. In particular $\dom$
 is a future complete regular domain}\par
\emph{
On the other hand we shall see that if $\Omega$ is a $\Gamma_\tau$-invariant regular
domain , then the cosmological time $T_\Omega$ is regular and
$\Omega=D(\widetilde S_a)$ where $\widetilde S_a$ is the  level surface $T_\Omega^{-1}(a)$. Moreover
we shall prove  that $\widetilde S_a$ is a $\Gamma_\tau$-invariant spacelike complete hypersurface. Thus we
have that  $\Gamma_\tau$-invariant regular domains are domains of
dependence of some $\Gamma_\tau$-invariant future convex complete spacelike hypersurface.
}
\end{oss} 

We want to describe the cosmological time on $\dom$ and in general on a future
 complete regular domain.
 First we show that every future complete convex set \emph{which has at
least a spacelike support plane} has a regular cosmological time which is a
$\mathrm C^1$-function.
\begin{prop}\label{CT and retraction on future convex set}
Let $A$ be a future complete convex  subset of $\Kn{M}{n+1}$ and $S=\partial
A$. Suppose that there exists a spacelike support plane. 
Then for every $p\in A$ there exists a unique $r(p)\in S$ which maximizes
the Lorentzian distance in $\overline A\cap\cpass(p)$. Moreover the map $p\mapsto r(p)$ is continuous.\par
The point $r=r(p)$ is the unique point in $S$ such that the plane
$r+\ort{(p-r)}$ is a support plane for $A$.\par
The cosmological time of $A$ is expressed by the formula
\[
      T(p)=\sqrt{-\E{p-r(p)}{p-r(p)}}.
\]
It is $\mathrm{C}^1$ and $-T$ is convex. The Lorentzian gradient of $T$ is given
by the formula
\[
      \nabla_LT(p)=-\frac{1}{T(p)}\left(p-r(p)\right).
\]
\end{prop}
\Dim
Since $A$ is convex the Lorentzian distance in $A$ is the restriction of the
Lorentzian distance in $\Kn{M}{n+1}$, that is 
\[
              d(p,q)=\sqrt{-\E{p-q}{p-q}}\quad\textrm{ for every }q\in
              A\textrm{ and }p\in\fut(A).
\]
Fix $p\in A$ and let $P$ be a spacelike support plane of $A$. 
Notice that $\cpass(p)\cap\cfut(P)$ is compact and $\cpass(p)\cap A\subset\cpass(p)\cap\cfut(P)$.
Thus there exists a point $r\in\overline A\cpass(p)$ which maximize the Lorentzian distance
from $p$. Clearly $r$ lies into the boundary $S$.\par
\begin{figure}
\begin{center}
\input{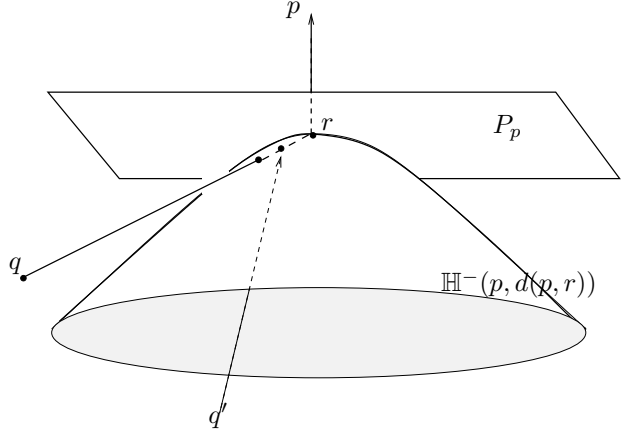}
\caption{$P_p$ is a support plane for $A$.}\label{support}
\end{center}
\end{figure}
Now we have to show that $r$ is unique. Let us define 
\[
\mathbb H^-(p,\alpha)=\{x\in\pass(p)|\,
d(p,x)=\alpha\}.
\] 
We know that $\mathbb H^-(p,d(p,r))$ is a past convex spacelike surface. Let
$r'\in \mathbb H^-\big(p,d(p,r)\big)\cap S-\{r\}$, the segment $(r,r')$ is contained
in $\pass(\mathbb H^-(p,d(p,r)))$ so that for
every $s\in (r,r')$ it turns out that $d(p,s)>d(p,r)$. On the other
hand we have that $(r,r')\in \overline A$ and this contradicts the choice of $r$.\par 
We have to see that the map $p\mapsto r(p)$ is continuous. Let $p_k\in A$ such
that $p_k\rightarrow p\in A$ and let $r_k=r(p_k)$. First let us show that
$\{r_k\}_{k\in\mathbb N}$ is bounded. Let $q\in\fut(p)$, notice that there
exists $k_0$ such that $p_k\in\cpass(q)$ for every $k\geq k_0$, so that 
$r_k\in\cpass(q)\cap S$ for $k\geq k_0$. Since $\cpass(q)\cap S$ is compact we
deduce that $\{r_k\}$ is contained in a compact subset of $S$.
Hence it is sufficent to prove that if $r_k\rightarrow r$ then
$r=r(p)$.
Let $q\in A$ then we have $\E{p_k-r_k}{p_k-r_k}\leq\E{p_k-q}{p_k-q}$. By
passing to the limit we obtain that $r$ maximizes the Lorentz distance.\par
Let $p\in A$ and $P_p$ be the plane $r(p)+\ort{(p-r(p))}$. We claim that $P_p$
is a support plane for $A$. In fact notice that $P_p$ is the tangent plane of
the set $\mathbb H^-(p,d(p,r))$ at the point $r(p)$.
Suppose that there exists 
$q\in A\cap\pass(P_p)$: we have  $(q,r)\subset A$. On the other hand we have
that there exists $q'\in(q,r)\cap\pass\left(\mathbb
  H^-(p,d(p,r))\right)$ (see fig.\ref{support}). Then $d(p,q')>d(p,r)$ and this is a contradiction.
Conversely let $s\in A$ such that $s+\ort{(p-s)}$ is a support plane for $A$.
An analogous argument shows that $s$ is in the past of $p$ and maximizes the
Lorentzian distance.\par
Now we can prove that the cosmological time $T$ is $\mathrm C^1$. We
shall use the following elementary fact:\\

\emph{ Let $\Omega\subset\mathbb R^N$ be an open set, and
  $f:\Omega\rightarrow\mathbb R$. Suppose that there exists
  $f_1,f_2:\Omega\rightarrow\mathbb R$ such that:}\\
\emph{1. $f_1\leq f\leq f_2$;}\\
\emph{2. $f_1(x_0)=f_2(x_0)=f(x_0)$;}\\
\emph{2. $f_1$ and $f_2$ are $\mathrm C^1$  and $\mathrm
  df_1(x_0)=\mathrm df_2(x_0)$.}\\
\emph{ Then $f$ is differentiable in $x_0$ and $\mathrm df(x_0)=\mathrm
  df_1(x_0)$.}\\

Let $p\in A$ and $r=r(p)$; fix a set of orthonormal affine coordinates
$(y_0,\ldots,y_n)$ such that the origin is $r(p)$ and $P_p$ is the plane 
$\{y_0=0\}$ (so that the point $p$ has coordinates $(\mu,0,\ldots,0)$).
Consider the functions
$f_1:A\ni y\mapsto+y_0^2-\sum_{i=1}^n y_i^2\in\mathbb R$ and $f_2: A\ni
y\mapsto y_0^2\in\mathbb R$. It is straightforward to recognize that
$f_1\leq T^2\leq f_2$ and $f_1(p)=T^2(p)=f_2(p)$. Moreover we have
$\nabla_Lf_1(p)=-2\frac{\partial}{\partial y_0}=\nabla_Lf_2(p)$. It
turns out that $T^2$ is differentiable in $p$ and
$\nabla_L\left(T^2\right)=-2(p-r)$. Thus $T$ is  differentiable in $p$
and$ \nabla_LT(p)=-\frac{1}{T(p)}\left(p-r\right)$.\par
Finally let us show that $-T$ is convex. Let $\varphi(p):=-T^2(p)=\E{p-r(p)}{p-r(p)}$.
 We have to show that 
\begin{equation}\label{CT is concave}
    -\varphi\left(tp+(1-t)q\right)\geq\left(t\sqrt{-\varphi(p)}+(1-t)\sqrt{-\varphi(q)}\right)^2
\qquad\textrm{ for all } p,q\in A\textrm{ and }t\in[0,1].
\end{equation}
Since $t\, r(p)+(1-t)r(q)\in A$ we have by the definition of $r$ that
\begin{eqnarray*}
   -\varphi(tp+(1-t)q)\geq\qquad\qquad\qquad\qquad\qquad\qquad
         \qquad\qquad\qquad\qquad\qquad\qquad\qquad\qquad\quad\quad \\
     \geq-\E{\left(tp+(1-t)q\right)\,-\,\left(t\,r(p)+(1-t)r(q)\right)\,}
       {\,\left(tp+(1-t)q\right)\,-\,\left(t\,r(p)+(1-t)r(q)\right)}=\\
      =-\E{t\left(p-r(p)\right)\,+\,(1-t)\left(q-r(q)\right)\,} 
{\,t\left(p-r(p)\right)\,+\,(1-t)\left(q-r(q)\right)}=\qquad\qquad\qquad\quad\\
      =-\Big(t^2\varphi(p)+(1-t)^2\varphi(q) +2t(1-t)\E{p-r(p)}{q-r(q)}\Big).
\end{eqnarray*}
Since $p-r(p)$ and $q-r(q)$ are future directed timelike vectors we have that
$\E{p-r(p)}{q-r(q)}\leq-\sqrt{\varphi(p)\varphi(q)}$ so that we have
\begin{eqnarray*}
  -\varphi(tp+(1-t)q)\geq \left(t^2 (-\varphi(p))+(1-t)^2(-\varphi(q))
  +2t(1-t)\sqrt{\varphi(p)\varphi(q)}\right)=\\
  \left(t\sqrt{-\varphi(p)}-(1-t)\sqrt{-\varphi(q)}\right)^2.
\end{eqnarray*}
\cvd
\begin{corollario}\label{CT limit}
We use  the notation of proposition \ref{CT and retraction on future convex set}.
For all $(p_k)_{k_\in\mathbb N}\subset A$ such that $p_k\rightarrow\overline
p\in\partial A$ we have
\[
    \lim_{k\rightarrow+\infty}T(p_k)=0.
\]
\end{corollario}
\Dim
Let $q\in\fut(\overline p)$, it is easy to see that $p_k\in\fut(q)$ for all
$k>>0$. By arguing as in proposition \ref{CT and retraction on future convex
  set} we see that $\{r(p_k\}$ is a bounded set. Up to passing to a
subsequence we can suppose that $r(p_k)\rightarrow\overline r$. Since
$p_k-r(p_k)$ is a timelike vector, $\overline p-\overline r$ is a
non-spacelike vector. On the other hand since $S$ is achronal set we get that
$\overline p-\overline r$ is null vector. Since
$T^2(p_k)=-\E{p_k-r(p_k)}{p_k-r(p_k)}$ we get
$\lim_{k\rightarrow+\infty} T^2(p_k)=-\E{\overline p-\overline r}{\overline p-\overline r}=0$.
\cvd
\begin{corollario}
With the above notation  let $\widetilde S_a=T^{-1}(a)$ for $a>0$. Then we
have that $\widetilde S_a$ is a future convex spacelike hypersurface and 
$T_p\widetilde S_a=\ort{(p-r(p))}$ for all $p\in \widetilde S_a$.\par
We have that $\fut(\widetilde S_a)=\bigcup_{b>a}\widetilde S_b$. Moreover let $r_a:\fut(\widetilde S_a)\rightarrow
\widetilde S_a$ be the projection and $T_a$ the CT of $\fut(\widetilde S_a)$ then we have:
\begin{eqnarray*}
r_a(p)=\widetilde S_a\cap[p,r(p)]\\
T_a(p)=T(p)-a.
\end{eqnarray*}
\end{corollario}
\cvd
Let $A$ be as in proposition \ref{CT and retraction on future convex set}, $T$
the CT of $A$ and $r:A\rightarrow\partial A$ the retraction. 
The \textbf{normal field} on $A$ is the map $N:A\rightarrow \Kn{H}{n}$
defined by the rule $N(p):=\frac{1}{T(p)}(p-r(p))$, it coincides up to the sign
with the lorentzian gradient of $T$ on $A$ (we have defined $N=-\nabla_L T$ instead
of $N=\nabla_L T$ because we want that $N$ is future directed). Let $\widetilde S_a=T^{-1}(a)$ then
$N|_{\widetilde S_a}$ is the normal field on $\widetilde S_a$.
Notice that the following identity holds
\[
   p= r(p)\,+\,T(p)N(p)\qquad\textrm{ for all } p\in A.
\]
Thus every point in $A$ is decomposed in \emph{singularity part} $r(p)$ and a
\emph{hyperbolic part} $T(p)N(p)$. We shall see that such decomposition plays
an important r\^ole to recover the announced duality. The following inequalities
are consequence of the fact that $r(p)+\ort{N(p)}$ is a support plane of $A$.
\begin{corollario}
With the above notation we have that
\begin{eqnarray}\label{fundamental inequality}
  \E{q}{p-r(p)}<\E{r(p)}{p-r(p)}&\nonumber\\
  \E{T(p)N(p)-T(q)N(q)}{r(p)-r(q)}\geq 0&\quad\textrm{ for all } p,q\in A.
\end{eqnarray}
\end{corollario}
\cvd
We denote by $\Sigma_A$ the image of the retraction $r:A\rightarrow\partial A$
and we refer to it as the \textbf{ singularity in the past}. Notice that if $r_0=r(p)$
then by (\ref{fundamental inequality}) the plane $r_0+\ort{(p-r_0)}$ is a
\emph{spacelike} support plane for $A$. Conversely let
$r_0\in\partial A$ and suppose that there exists a future directed timelike vector $v$
such that the spacelike plane $r_0+\ort{v}$ is a support plane. We have  that
$p_\lambda=r_0+\lambda v\in A$ for 
$\lambda>0$ and by proposition \ref{CT and retraction on future convex set}
we have $r(r_0+\lambda v)=r_0$.
 \begin{corollario}\label{fiber of retraction}
Let $A$ a future convex set wich has a spacelike support plane. Then
$r_0\in\Sigma_A$ if and only if there exists a timelike vector $v$ such that the
plane $r_0+\ort{v}$ is a support plane. Moreover 
\[
       r^{-1}(r_0)=\{r_0+v|\,r_0+\ort{v}\textrm{ is a support plane}\}.
\]
\end{corollario}
\cvd
\begin{oss}\emph{
Notice that the map $r: A\rightarrow \Sigma_A$  continuously extends to a retraction
$r:A\cup\Sigma_A\rightarrow \Sigma_A$. This map is  a  deformation
retraction (in fact the maps $r_t(p)=t(p-r(p))+r(p)$ give the
homotopy). Therefore $\Sigma_A$ is} contractile.
\end{oss} 
Now let $\Omega$ be a \textbf{future complete  regular domain}. Let us use
this notation:
\begin{itemize}
\item
$T$ is the cosmological time on $\Omega$ and $\widetilde S_a=T^{-1}(a)$;
\item
$r:\Omega\rightarrow\partial\Omega$ is the retraction and
$N:\Omega\rightarrow\Kn{H}{n}$ the normal field;
\item
$\Sigma=r(\Omega)$ is the singularity in the past.
\end{itemize}
\begin{lemma}
$\widetilde S_a$ is a Cauchy surface for $\Omega$. Moreover
$\Omega$ is the domain of dependence of $\widetilde S_a$.
\end{lemma}
\Dim
Since $\Omega$ is a regular domain we have that $D(\widetilde S_a)\subset\Omega$. Now let
$p\in\Omega$ and let $v$ be a future directed non-spacelike vector. There
exists $\lambda>0$ such that $p+\lambda v\in\fut(\widetilde S_a)$ so that $T(p+\lambda
v)>a$. On the other hand there exists $\mu<0$ such that $p+\mu
v\in\partial\Omega$, and by corollary \ref{CT limit} we have  that $\lim_{t\rightarrow \mu}T(p+ t v)=0$. So
that there exists a $\lambda'\in\mathbb R$ such that $T(p+\lambda' v)=a$.
Thus $\Omega\subset D(\widetilde S_a)$ and the proof is complete.
\cvd
\begin{oss}\label{quotient of CT level surface}
\emph{
If $\Omega$ is a} $\Gamma_\tau$-invariant \emph{regular domain then $T$ is a
$\Gamma_\tau$-invariant function. It follows that the CT level surface $\widetilde S_a$ are
$\Gamma_\tau$-invariant future convex spacelike hypersurface. Moreover we
have  $r\circ\gamma_\tau=\gamma_\tau\circ r$ and
$N\circ\gamma_\tau=\gamma\circ N$. Thus $\Sigma$ is $\Gamma_\tau$-invariant
subset of $\partial\Omega$.}\par
\emph{
From the last lemma it follows that  $\widetilde S_a/\Gamma_\tau$ is a Cauchy surface of
$\Omega/\Gamma_\tau$. In particular if we take $\Omega=\dom=D(\Ftilde)$ (see
section \ref{construction-sec}) we have that
$\widetilde S_a/\Gamma_\tau$ is homeomorphic to $\Ftilde/\Gamma_\tau\cong M$ (in fact two
Cauchy surface are always homeomorphic).
}
\end{oss} 
We want to give a more precise description of the map $r$ for a regular
domain $\Omega$. In paticular we want to describe the singularity and the
fiber of a point on the singularity in terms of geometric properties of the
boundary $\partial \Omega$.
Let us start with a simple remark.
\begin{lemma}\label{null support plane for domain of dependence}
For every $p\in \Omega$ there exists a future directed null vector $v$ such
that the ray $p+\mathbb R_+v$ is contained in $\partial\Omega$.
Furhermore we have
\[
     \Omega=\bigcap\{\fut(p+\ort v)|\, p\in\partial\Omega \textrm{ and } v
                         \textrm{ is a null vector such that }
                          p+\mathbb R_+v\subset\partial\Omega\}.
\]
\end{lemma}
\Dim
Let $p\in\partial \Omega$, then 
there exists a null future directed vector $v$ such that the 
ray $p+\ort v$ is a support plane for $\Omega$. Since $\fut(p)\subset\Omega$
we have that the ray $p+\mathbb R_+ v$  is contained in $\partial\Omega$.\par
Conversely suppose that a ray $R=p+\mathbb R_+v$ is contained in
$\partial\Omega$. 
By Han-Banach theorem there exists a hyperplane $P$ such that
$\Omega$ and $R$ are contained in the opposite (closed)
semispaces bounded by $P$. Since $\Omega$ is future complete and $P$ is a
support plane of $\Omega$ it results that $P$ is not timelike.
Since $P$ does not intersect $R$ transversally it follows that $v$ is parallel
to $P$ so that $P=p+\ort v$.
\cvd
\begin{prop}\label{singularity for domain of dependence}
A point $p\in\partial\Omega$ lies in $\Sigma$ if and only if there exist at
least $2$ future directed null rays which are contained in $\partial\Omega$ and
pass through $p$.\par
Moreover let $p\in\Sigma$ then $r^{-1}(p)$ is the intersection of $\Omega$
with the convex hull of the null rays which are contained in $\partial\Omega$
and pass through $p$. 
\end{prop}
\Dim
We use the following elementary fact about convex sets:\\

\emph{Let $V$ be a vector space and $G\subset V^*$. Consider the convex 
   $K=\{v\in V| g(v)\leq C_g\textrm{ for }g\in G\}$.
    Suppose that the following property holds:
    if $g_n\rightarrow g$ and $C_{g_n}\rightarrow C$  then $g\in G$ and
   $C_g\leq C$. Then for all $v\in\partial K$
   the set $G_v=\{g\in G|C_g=g(v)\}$ is non empty. Moreover  the plane $v+P$ is
   a support plane of $K$ if and only if there exists $h$ in the convex hull
   of $G_v$ such that $P=\ker h$.}\\

Consider the family $L$ of the null future directed vectors which are orthogonal
to some null support plane. For every $v\in L$ let 
$C_v=\sup_{r\in\Omega}\E{v}{r}$. By 
lemma \ref{null support plane  for domain of dependence} we get
\[
\Omega=\{x\in\Kn{M}{n+1}|\E{x}{v}\leq C_v\,\forall v\in L\}.
\] 
Now fix $p\in\partial\Omega$ and let $L(p)$ be the set of the null future
directed vectors $v$ such that $p+\ort v$ is a support plane of $\Omega$.  
We can apply the remark on convex sets stated above to the
family $L$ (in fact the scalar product $\E{\cdot}{\cdot}$ gives an
identification of $\Kn{R}{n+1}$ with its dual) and we obtain that $L(p)$ is non-empty. Moreover fix a future
directed  non-spacelike vector $v$, then $p+\ort v$ is a support plane if and
only if $v$ belongs to the convex hull of $L(p)$. By corollary \ref{fiber of
  retraction} we get that $r^{-1}(p)$ is the intersection of $\Omega$ with
the convex hull of $p+L(p)$.\par
Notice that a null future directed vector $v$ lies in $L(p)$ if and only if
$p+\mathbb R_+v$ is contained in $\partial\Omega$.
\cvd
Now for $p\in\Sigma$ let us define a subset of $\Kn{H}{n}$
\[
    \mathcal{F}(p):=N(r^{-1}(p)).
\]
In the following corollary we point out that $\mathcal F(p)$ is an ideal convex set of
$\Kn{H}{n}$. We recall that a convex set $C$  of $\Kn{H}{n}$ is \emph{ideal} if it is
the convex hull of boundary points. 
\begin{corollario}
Fix $p\in\Sigma$ and let $L(p)$ as in the last proposition. Denote by $\hat
L(p)$ the set of the points on $\partial\mathbb H^n$ which correspond to points in
$L(p)$.
Then $\mathcal F(p)=N(r^{-1}(p))$ is the convex hull in $\Kn{H}{n}$ of $\hat L(p)$
\end{corollario}
\cvd
Thus we  see that each point $p$ in the singularity $\Sigma$ corresponds to
an ideal convex set $\mathcal F(p)$.
Now we study as the convex sets $\{\mathcal F(p)\}_{p\in\Sigma}$ stay in $\Kn{H}{n}$.
We recall that given two convex sets $C,C'\subset\Kn{H}{n}$ we say that a
hyperplane $P$ \emph{separes} $C$ from $C'$ if $C$ and $C'$ are contained in  the
opposite closed semispaces bounded by $P$.
\begin{prop}\label{domain gives stratification}
Let $\Omega$ be a future complete regular domain. 
For every $p,q\in\Sigma$  the plane $\ort{(p-q)}$ separes $\mathcal F(p)$ from
$\mathcal F(q)$.
The segment $[p,q]$ is contained in
$\Sigma$ if and only if  $\mathcal F(p)\cap \mathcal F(q)\neq\varnothing$. In this case for all $r\in (p,q)$ we have
\[
   \mathcal F(r)=\mathcal F(p)\cap\ort{(p-q)}=\mathcal F(q)\cap\ort{(p-q)}=\mathcal F(p)\cap \mathcal F(q).
\]
\end{prop} 
\Dim
The inequality (\ref{fundamental inequality}) implies that 
$\E{tv}{p-q}\leq\E{sw}{p-q}$
for every $v\in \mathcal F(q)$, $w\in \mathcal F(p)$ and $t,s\in\mathbb R_+$. This inequality
can be satisfied if and only if $\E{v}{p-q}\leq 0$ and $\E{w}{p-q}\geq 0$ for
all $v\in \mathcal F(q)$ and $w\in \mathcal F(p)$. This show that $\ort{(p-q)}$
separes $\mathcal F(p)$  from $\mathcal F(q)$.\par
Suppose now that $\mathcal F(p)\cap \mathcal F(q)\neq\varnothing$. We have that $\mathcal F(p)\cap
\mathcal F(q)$ is contained in $\ort{(p-q)}$. Let $v\in \mathcal F(p)\cap \mathcal F(q)$ and let $P_v$ be
the unique support plane which is orthogonal to $v$  and intersects
$\partial\Omega$.
We know that $P_v$ passes through $p$ and $q$ so that the segment $[p,q]$ is
contained in $\partial\Omega$. Finally since $P_v$ is a spacelike support plane
which passes through every $r\in(p,q)$ we have $[p,q]\subset\Sigma$ and
moreover $\mathcal F(p)\cap \mathcal F(q)\subset \mathcal F(r)$.\par
Conversely suppose that $[p,q]$ is contained in $\Sigma$. Let $r\in (p,q)$ and
$v\in \mathcal F(r)$. Then we have that $\E{v}{p-r}\leq 0$ and $\E{v}{r-q}\geq
0$. Since $p-r$ and $r-q$ have the same direction we argue that
$\E{v}{r-q}=0$ and $\E{v}{r-p}=0$ so that $v\in \mathcal F(p)\cap \mathcal F(q)$.\par
In order to conclude the proof we have to show that $\mathcal F(r)\supset
\mathcal F(p)\cap\ort{(p-q)}$. We know that $\mathcal F(p)\cap\ort{(p-q)}$ is the convex hull
of
$\hat L(p)\cap\ort{(p-q)}$. Thus it is sufficient to show that $L(r)\supset
L(p)\cap\ort{(p-q)}$.
Now fix $v\in L(p)\cap\ort{(p-q)}$ and consider the plane $P=p+\ort v$. The
intersection of this plane with $\overline\Omega$ includes the ray $p+\mathbb R_+v$ and the segment $[p,q]$.
Since this intersection is a convex we have that $r+\mathbb R_+v$ is a subset
of $P\cap\overline\Omega$  and thus $v\in L(r)$.
\cvd
Let us give a general definition.
\begin{defi}
A \textbf{geodesic stratification} of $\Kn{H}{n}$ is a family $\mathcal C=\{C_i\}_{i\in I}$ such
that
\begin{enumerate}
\item
$C_i$ is an ideal convex set of $\Kn{H}{n}$;
\item
$\Kn{H}{n}=\bigcup_{i\in I}C_i$;
\item
For every $i,j\in I$ (with $i\neq j$) there exists a support plane $P_{i,j}$
which separes $C_i$ from $C_j$. Furthermore if $C_i\cap C_j\neq\varnothing$
then $C_i\cap C_j=C_i\cap P_{i,j}=C_j\cap P_{i,j}$ .
\end{enumerate}
Every $C_i$ is called  \emph{piece} of the stratification.\par
We say that the stratification is $\Gamma$-invariant if
$\gamma(C_i)\in\mathcal C$ for all
$\gamma\in\Gamma$ and for all $C_i\in\mathcal C$  
\end{defi}
If $\Omega$ is a future complete regular domain of $\Kn{M}{n+1}$ such that \emph{the
normal fied $N$ is surjective} we have that $\{\mathcal F(p)\}_{p\in\Sigma}$ is a
geodesic stratification.\par 
If $\Omega$ is a \emph{$\Gamma_\tau$-invariant} future complete regular domain
 by lemma \ref{support planes}  the normal field $N$ is
surjective. In this case it is evident that the stratification
$\{\mathcal F(p)\}_{p\in\Sigma}$ is $\Gamma$-invariant.\par
Let $C$ be a convex set, we say that a point $p$ is \emph{internal} if all support
planes which pass through $p$ contains $C$. Let us denote by $bC$ the set of
the point of $C$ which are not internal. Notice that $bC$ is not the topological
boundary. If $\dim C=k$ then $bC$ have a natural decomposition in
convex pieces which ideal convexes $C_i$ with $\dim C_i<k$ (see
\cite{Epstein}).\par
Now if $\mathcal C$ is a  geodesic stratification of $\Kn{H}{n}$ we can add to
it the convexes $C$ which are pieces of the decomposition of $bC_i$
for some $C_i\in\mathcal C$. It is easy to see that so we obtain a new
geodesic stratification $\overline{\mathcal C}$ which we call the \emph{completeness}
of $\mathcal C$. 
Notice that $\overline{\overline C}=\overline C$. A geodesic startification
wich coincides with its completeness is called \emph{complete}.\par
Now we can define the $k$-stratum of $\mathcal C$ (for $1\leq
k\leq n-1$) as the set
\[
       X_{(k)}=\bigcup\{F\in\overline{\mathcal C}|\dim F\leq k\}.
\]
Notice that if $\mathcal C$ is $\Gamma$-invariant then also
$\overline{\mathcal C}$ is so. Moreover in this case the strata are $\Gamma$-invariant subsets.\par
It is easy to see that $X_{(n-1)}$ is a closed set (in fact $\Kn{H}{n}-X_{(n-1)}$
is the union of the interior of the $n$-dimensional convex of $\mathcal C$).
Notice that if $n=2$  we have only the $1$-stratum which in fact is a geodesic
lamination of $\Kn{H}{2}$. Conversely if $L$ is a $\Gamma$-invariat geodesic lamination there is
a unique complete geodesic stratification $\mathcal C$ such that $L$ is the
$1$-stratum of $\mathcal C$.
For $n=2$ we know that the stratification is continuous in the sense that if
$r_k\in C_k$ and $r_k\rightarrow r\in\Kn{H}{n}$ then there exists a piece
$C\in\overline {\mathcal C}$ such that $C_k\rightarrow C$ with respect of the Hausdorff
topology.\par
Unfortunately in dimension $n>2$ the geodesic stratification are more complicated: the
strata $X_{(k)}$ are not closed for $k\neq n-1$ and we have not the
continuity. However, as we are going to see,  the stratifications which arise
from a particular class of future complete regular domains are weakly continuous in the following sense:
\begin{defi}
A geodesic stratification $\mathcal C$ is \emph{ weakly continuous} if the
following property holds. Suppose $x_k$ be a convergent sequence of
$\Kn{H}{n}$ and put $x=\lim x_k$. Let $F_k$ be a pieces which
contains $x_k$ and suppose that $F_k\rightarrow F$ in the Hausdorff
topology. Then there exists a piece $G\in\overline{\mathcal C}$ such that $F\subset G$.
\end{defi}
\begin{prop}
Let $\Omega$ be a future complete regular domain  such that the
normal field $N$ is surjective and the restriction $N|_{\widetilde S_1}$ is a proper map
(where $\widetilde S_1$ is the CT-level surface $T^{-1}(1)$). Then       
the geodesic stratification $\mathcal C$ associated with it
is weakly continuous.
\end{prop}
\Dim
Let $(x_k)\subset\Kn{H}{n}$ be a convergent sequence $x_k\rightarrow x$. Let
$F_k$ be a piece which contains $x_k$ and suppose that  $F_k\rightarrow F$. We have to see
that $F$ is contained in a piece $G$.\par
Let $r_k\in \Sigma$ such that $F_k=\mathcal F(r_k)$ and let $p_k=r_k+x_k\in 
\widetilde S_1$.
Since $N|_{\widetilde S_1}$ is a proper map there exists a convergent subsequence $p_{k(j)}$.
Let $p=\lim p_{k(j)}$ and $r=r(p)$. We want to show that $F$ is contained in
$\mathcal F(r)$. Since $F$ is the convex hull of $\hat L_F=F\cap\partial\mathbb H^n$ 
it is sufficient to show that $\hat L_F\subset
\hat L(r)$. Now let $[v]\in\hat L_F$, we know that there exist a sequence
$[v_n]\in\hat L(r_n)$ such that $[v_n]\rightarrow[v]$ in $\partial\Kn{H}{n}$.
We have that $r_n+\mathbb R_+ v_n\subset\partial\Omega$. Since
$\partial\Omega$ is closed this implies  that $r+\mathbb R_+
v\subset\partial\Omega$. Thus we conclude $[v]\in F_r$.
\cvd
Consider the regular domain $\dom$. Notice that the map $N:\widetilde S_1\rightarrow
\Kn{H}{n}$ induces to the quotient a map $\overline
N:\widetilde S_1/\Gamma_\tau\rightarrow \Kn{H}{n}/\Gamma = M$. Since $\widetilde S_1/\Gamma_\tau$ is
compact (in fact it is homeomorphic to $M$) it is easy to see that $N$ \emph{is a
proper map}. Thus the last proposition applies.
\begin{corollario}
Let $\mathcal C_\tau$ be the stratification associated with $\dom$. Then it is
weakly continuous.
\end{corollario}
\cvd
We postpone a more careful discussion about $\Gamma$-invariant geodesic stratification to the
last sections. In the last part of this section we consider the future
complete regular domain $\dom$ which we have constructed in section
\ref{construction-sec}. 
We prove that
$\Gamma_\tau$ does not act freely and properly discontinuously on
$\partial\dom$.
By this result we deduce that \emph{the action of $\Gamma_\tau$ on $\Kn{M}{n+1}$ is
not free and properly discontinuous}. In particular we deduce that the domain
of dependence of a $\Gamma_\tau$-invariant spacelike hypersurface is a regular
domain (either future or past complete).
\begin{lemma}
Let $\Omega$ be a future complete regular domain.
Suppose that $\Sigma$ is closed in $\partial\Omega$. Then the retraction
$\Omega\rightarrow\Sigma$ extends uniquely to a  deformation retraction
$r:\overline\Omega\rightarrow\Sigma$
\end{lemma}
\Dim
Since $\Sigma$ is closed it is easy to see that for every point $p$ outside $\Sigma$
there exists a unique null ray $R$ in $\partial\Omega$ such that $p$ is contained
in the interior of $R$. Thus we can define the retraction on
$\partial\Omega$ by taking $r(p)$ the starting point of the ray $R$.
It is easy to show that this is a continuous extension of $r$.
\cvd
Let $\Omega$ be a future complete regular domain and suppose that $\Sigma$ is
closed. Let $X:=\partial\Omega$, we want to construct a boundary of $X$.
We know that $X-\Sigma$ is a $\mathrm C^1$-manifold foliated by rays with starting points in
$\Sigma$. 
Let $p\in X-\Sigma$, we denote by $R(p)$ the ray of the foliation which passes
through $p$.
The retraction on $X$ is so defined: $r(p)=p$ if $p\in\Sigma$
whereas $r(p)$ is the starting point of $R(p)$ if $p\in X-\Sigma$.\par 
\emph{The boundary of $X$ is the leaves space of the foliation:}
\[
 \partial X:=\{R|R \textrm{ is a
  ray of the foliation}\}. 
\]
We want to define a topology on $\overline X:=X\cup\partial X$
such that agrees with the natural topology on $X$ and makes $\overline X$ a
$n$-manifold with boundary equal to $\partial X$. Thus we have to define a fundamental
family of neighborhoods of a point $R\in\partial X$. Fix a $\mathrm C^1$-embedded closed
$(n-1)$-ball  $D$ which intersects transversally the foliation and passes
through $R$ and define
\begin{eqnarray*}
   U(R,D)=\{p\in X-\Sigma| R(p)\cap \mathrm{int}D\neq\varnothing \textrm{ and }
   \mathrm{int}D\cap R(p) \textrm{ separes } p \textrm{ from
   r(p)}\}\cup\\
   \cup\{S\in\partial X|S\cap\mathrm{int} D\neq\varnothing\}.
\end{eqnarray*}
Then we can consider the topology on $\overline X$ which agrees with the
natural topology on $X$ and such that for each $R\in\partial X$ the sets $U(R,D)$ are a fundamental
family of neighborhoods of $R$. It is easy to see that $\overline X$ is a
Hausdorff space.
In order to construct an atlas on $\overline X$ for $p\in X-\Sigma$ put $v(p)$
the future directed null vector tangent to $R(p)$ such that  $x_0(v(p))=1$.
For all $(n-1)$-ball $D$ as above, consider the maps 
$\mu_D:D\times(0,+\infty]\rightarrow U(R,D)$ defined by
the rule
\[
 \mu_D (x,t)=\left\{\begin{array}{ll} x+tv(x)&\textrm{if } t<+\infty\\
                                         R(x) &\textrm{ if }
                                         t=+\infty\end{array}\right.
\]
It is easy to see that these maps are local charts, 
so that $\overline X$ is a manifold with boundary.
Finally notice that the retraction $r$ extends uniquely to a retraction
$r:\overline X\rightarrow \Sigma$. Notice that this \emph{ retraction is a proper map}.\par   
Now suppose that $\Gamma_\tau$ acts freely and properly discontinuously on
$\overline{\Omega}$.
It is easy to see that the action of $\Gamma_\tau$ on $X$ extends uniquely to
an action on $\overline X$. Moreover the map $r:\overline X\rightarrow\Sigma$
is $\Gamma_\tau$-equivariant. By using this remark it follows that the action
of $\Gamma_\tau$ on $\overline X$ is free and properly discontinuous. Thus
$\overline X/\Gamma_\tau$ is a manifold with boundary.
Now we can state the annouced proposition.  
\begin{prop}
The action of $\Gamma_\tau$ on $\partial\dom$ is not free and properly discontinuous.
\end{prop}
\Dim
By contradiction suppose that the action is free and properly
discontinuous. Let $X:=\partial\dom$, and put $M':=X/\Gamma_\tau$ and
$K:=\Sigma/\Gamma_\tau$.
Let $\hat r:\dom/\Gamma_\tau\rightarrow K$ be the map which arises from the
retraction $r:\dom\rightarrow\Sigma$.
Notice that $K=\hat r(\widetilde S_1/\Gamma_\tau)$. Since $\widetilde
S_1/\Gamma_\tau$ is compact (in fact we have seen that it is homeomorphic to
$M$) we get that $K$ is compact. Since $M'$ is a
Hausdorff space  we have that $K$ is
closed in $M'$ and so $\Sigma$ is closed in $X=\partial\dom$.\par
Thus we can  construct as above the boundary $\partial X$ of $X$. Let  $\overline
M'=\overline X/\Gamma_\tau$, we know that $\overline M'$ is a manifold with
boundary and $M'$ is the interior of $\overline M'$.\par
 Consider the retraction $r:\overline X\rightarrow\Sigma$: since this map is
 $\Gamma_\tau$-equivariant and proper it induces to the quotient a proper map
$\overline r:\overline M'\rightarrow K$. Since $K$ is compact it follows that
$\overline M'$ is a compact manifold with boundary.\par
Since $\overline r:\overline M'\rightarrow K$ is a deformation retraction we
have that $\mathrm{H}_n(K)=\mathrm H_n(M')$ and by the Poicar\'e duality it
follows that
\[
    \mathrm H_n(K)=\mathrm H_n(M')=\mathrm H^0(\overline M',\partial\overline M')=0.
\]
On the other hand let $Y_\tau=\dom/\Gamma_\tau$ and $\overline
Y_\tau=\overline\dom/\Gamma_\tau$.
We know that the map $r:\overline{ \dom}\rightarrow \Sigma$ induces to the
quotient a deformation retraction $\overline Y_\tau\rightarrow K$. So that
\[
   \mathrm H_n(K)=\mathrm H_n(\overline Y_\tau)=\mathrm H_n(Y_\tau).
\]
We have $Y_\tau\cong\mathbb R\times M$. So $\mathrm H_n(Y_\tau)=\mathrm
H_n(M)=\mathbb Z$ and this is a contradiction. 
\cvd    
\begin{corollario}
$\Gamma_\tau$ does not act freely and properly discontinuously on the whole
$\Kn{M}{n+1}$.
Moreover let $\widetilde F$ be a $\Gamma_\tau$-invariant complete spacelike hypersurface
such that $\Gamma_\tau$ acts freely and properly discontinuously on it.
Then $D(\widetilde F)$ is a regular domain (either future or complete).
\end{corollario}
\Dim
The first statement follows from the last proposition.
In particular we have that $D(\widetilde F)$ is not the whole $\Kn{M}{n+1}$.
By corollary \ref{convexity of domain of dependence} we know that $D(\widetilde F)$ 
is either future or past complete and it is the intersection of the
future (resp. past) of null planes. By lemma \ref{support planes} we know
that $D(\widetilde F)$ has spacelike support planes and so it is a regular domain.
\cvd 
\section{Uniqueness of the Domain of Dependence}
Let us summarize what we have seen until now.
We have fixed $\Gamma$ a free-torsion co-compact discrete subgroup of
$\SOO^+(n,1)$ and $M=\Kn{H}{n}/\Gamma$. Then we have fixed a cocycle $\tau\in
Z^1(\Gamma,\Kn{R}{n+1})$ and we have studied the $\Gamma_\tau$-invariant 
domains of $\Kn{M}{n+1}$. In particular we have constructed a
$\Gamma_\tau$-invariant  future complete (resp. past complete)
regular domain $\mathcal D_\tau$ (resp.$\dom^-$) such that the action on $\Gamma_\tau$ on it
is free and properly discontinuous and the quotient $Y_\tau=\mathcal
D_\tau/\Gamma_\tau$ is a globally hyperbolic manifold homeomorphic to $\mathbb
R_+\times M$ with regular CT.
Moreover we have seen that if $\widetilde F$ is a $\Gamma_\tau$-invariant
complete spacelike hypersurface such that the action on it is free and
properly discontinuous then $D(\widetilde F)$ is a regular domain either future or
past complete.\par
In this section we want to show that $\dom$ (resp. $\dom^-$) is the unique
$\Gamma_\tau$-invariant future complete (resp. past complete) regular
domain. In particular we deduce that every $\widetilde F$ as above is contained in
$\dom$ or in $\dom^-$ and it is in fact a Cauchy surface of it.
\begin{teo}\label{uniqueness}
$\dom$ is the unique  $\Gamma_\tau$-invariant future complete regular domain.
\end{teo}
Let us give a scheme of the proof. Let $\Omega$ be a $\Gamma_\tau$-invariant
future complete regular domain  we have to show that $\Omega= \dom$.\par
Let $T_\Omega$  be the cosmological time on $\Omega$ (whereas $T$ is the
cosmological time on $\dom$). 
For every $a>0$ let  $\widetilde S_a^\Omega:=T_\Omega^{-1}(a)$ (whereas
$\widetilde S_a=T^{-1}(a)$ is the level surface of $\dom$). Since $\widetilde S_a^\Omega$
(resp. $\widetilde S_a$) is 
Cauchy surface in $\Omega$ (resp. $\dom$) we have that
$\Omega=D(\widetilde S_a^\Omega)$ (resp. $\dom=D(\widetilde S_a)$). Thus it is sufficient
to prove that $\widetilde S_a\subset\Omega$ and $\widetilde S_a^\Omega\subset \dom$ for
$a>>0$. 
Let us split the proof in some steps.\par
\emph{Step 1.} $\Omega\cap \dom\neq\varnothing$.\par
\emph{Step 2.} Fix $p\in\Omega\cap \dom$ and let $C$ 
be  \emph{the convex hull} of the $\Gamma_\tau$-orbit of $p$. Then $C$ is a
future complete
convex set.\par
\emph{Step 3.} Let  $\Delta=\partial C$ be the boundary of $C$ then
$\Delta/\Gamma_\tau$ is compact.\par
\emph{Step 4.} Let
$a>\sup_{q\in\Delta}T_\Omega(q)\vee\sup_{q\in\Delta}T(q)$ then 
     $\widetilde S_a^\Omega\subset \dom$ and $\widetilde S_a\subset\Omega$.\\
 
The \emph{first step} is quite evident. In fact let $p\in\Omega$ and $q\in \dom$.
Since $\Omega$ and $\dom$ are future complete then
$\fut(p)\cap\fut(q)\subset\Omega\cap \dom$. On the other hand the future sets of
two points in $\Kn{M}{n+1}$ are not disjoint.\par

The \emph{second step} is more difficult. We start with a lemma.
\begin{lemma}\label{convex sets}
Let $C$ be a closed convex set of $\Kn{M}{n+1}$ whose interior part is nonempty. Then
one of the following statement holds:
\begin{itemize}
\item[(1)]
there exists a non spacelike direction $v$ such that
$C=\{x\in\Kn{M}{n+1}|\alpha_1\leq\E{x}{v}\leq\alpha_2\}$;
\item[(2)]
there exists a timelike support plane for $C$;
\item[(3)]
every support plane is non-timelike and $C$ is a future or past convex set. 
\end{itemize}
\end{lemma}
\Dim
Suppose there exists $p,q\in \partial C$ such that $q\in\fut(p)$. 
Let us prove that $C$ verifies (1) or (2). More exactly suppose that $C$ has
not timelike support plane we want to show that $C$ satifies (1).\par
Let $P_p$ and $P_q$ be non-timelike support planes respectively in $p$ and in $q$. 
Since $q\in\fut(p)$ we argue that $C\subset\fut(P_p)\cap\pass(P_q)$. As in
the proof of corollary \ref{support planes}  if $P_p$ and $P_q$
are not parallel then we get that there exists a timelike support plane. Thus
$P_p$ and $P_q$ are parallel.\par 
Since $p$ and $q$ are generic it follows that
for all $p',q'\in\partial C$ such that $p'\in\fut(q')$  support planes in
$p'$ and in $q'$ are parallel. In particular there exists a nonspacelike
direction $v$ such that the unique support plane in $p'$ and the unique
support plane in $q'$ are orthogonal to $v$.\par
We want to show that $\partial C=P_p\cup P_q$. In fact it is sufficient to
prove the inclusion $P_p\cup P_q\subset\partial C$. Now let $e=q-p$ and define
the set 
\[
A=\{z\in P_q|z\in\partial C \textrm{ and } z-e\in\partial C\}.
\]
We have to show that $A=P_q$. Clearly  $A$ is closed and nonempty (in
fact $q\in A$). Thus it is
sufficient to show that $A$ is open. Now let $z\in A$ and $z'=z-e$: the set
$U:=\fut(z')\cap\partial C$ is a neighborhood of  $z$ and for
every $x\in U$ there is a unique support plane $P_x$ parallel to $P_q$ which passes
through $x$. Since $P_q$ is a support plane and $P_q\cap\partial C\neq\varnothing$
it follows that  $P_q=P_x$ so that $U\subset P_q$. In the same way we deduce
that $V=\pass(z)\cap\partial C$ is contained in $P_p$. Thus $U\cap(V+e)$ is a
neighborhood of $z$ in $\partial C$ which is contained in $A$. Thus we have
that $\partial C=P_p\cup P_q$ and it follows that $C$ verifies (1).\\

We have proved that if there exists $p,q\in\partial C$ such that $p-q$ is
timelike then $C$ verifies (1) or (2). Suppose now that for all
$p,q\in\partial C$ the vector $p-q$ is non-timelike, we want to show that $C$
verifies (3).\par
Suppose that there exists a timelike support plane. Then there exists a vector
$u$ and $K\in\mathbb R$ such that $\E{u}{u}=1$ and
\[
    \E{u}{p}\leq K\qquad \textrm{for all }p\in C.
\]
Take $p_0\in\mathrm{int}C$ and $v_+,v_-\in\Kn{H}{n}$ such that
\[
    \E{u}{v_+}>0\qquad\E{u}{v_-}<0.
\]
Consider for $t>0$
\[
      p_t=p_0+tv_+\qquad p_{-t}=x_0-tv_-.
\]
We have
\begin{eqnarray*}
\E{p_t}{u}=\E{p_0}{u}+t\E{v_+}{u}\rightarrow+\infty\quad\textrm{ for }t\rightarrow+\infty;\\
\E{p_{-t}}{u}=\E{p_0}{u}-t\E{v_-}{u}\rightarrow+\infty\quad\textrm{ for }t\rightarrow+\infty.
\end{eqnarray*}
Thus there exists $\alpha>0$ and $\beta<0$ such that $p_\alpha,p_\beta\notin
C$. Let $t_+=\sup\{t>0|x_t\in C\}$ and $t_-=\inf\{t<0|x_t\in C\}$.
Notice that $p_{t_+}$ and $p_{t_-}$ lies in $\partial C$. Furthermore since
$p_0$ is an interior point of $C$ we have that $t_+> 0$ and $t_-<0$. Thus
$p_{t_+}$ is in the future of $p_0$ which is in the future of $p_{t_-}$. But
then we have $p_{t_+}$ in the future of $p_{t_-}$ and this contradicts the our
assumption on $C$.\par
Hence all the support planes of $C$ are  non-timelike. Let $P$ be a support
plane. 
 We can suppose that $C\subset\overline{\fut(P)}$ (the other
case is analoguous). I claim that  for
every support plane $Q$ we have that $C\subset\overline{\fut(Q)}$. 
Otherwise there should exist $v_1,v_2$ future directed non-spacelike vectors
and $K\in\mathbb R$ such that
\[
    \E{p}{v_1}\leq K\qquad\E{p}{v_2}\geq K\qquad\textrm{for all } p\in C.
\]
Fix $p_0\in\mathrm{int}C$ and put $p_t=p_0+te$ where $e$ is a timelike vector.
Then it is easy to show that there exist $t_1<0<t_2$ such that
$p_{t_1},p_{t_2}\in\partial C$ and this contradicts the our assumption on $C$.
Thus we have
\[
   C=\bigcap_{P\textrm{ support plane}}\overline{\fut(P)}.
\]
It follows that $C$ is future complete.
\cvd
Now let us go back to the step 2. We have taken $p\in\Omega\cap \dom$ and we
have to show that the convex hull $C$ of the $\Gamma_\tau$-orbit of $p$ is
future complete. By last lemma it is sufficient to prove:\\
a) $\mathrm{int}C\neq\varnothing$;\\
b) $C$ is not of the form
$\{x\in\Kn{M}{n+1}|\alpha_1\leq\E{x}{v}\leq\alpha_2\}$;\\
c) $C$ has not timelike support plane.\\

a)  is quite evident. In fact if the interior of $C$ is empty then there exists
a unique $k$-plane $P$ with $0<k<n+1$ such that $C\subset P$ and
$\mathrm{int}_P(C)\neq\varnothing$.
But then since $C$ is $\Gamma_\tau$-invariant it follows that $P$ is
$\Gamma_\tau$-invariant
and so the tangent plane to $P$ is $\Gamma$-invariant. But we know that a
$\Gamma$ is co-compact and so it is irreducible.
In analogous way one can prove b).\par
\emph{It remains to prove that $C$ has not timelike support plane}. For this purpose
we introduce some notation. Fix a set of orthonormal affine coordinates $(y_0,\ldots,y_n)$.
For every $\gamma\in\Gamma$ we denote by
$x^+(\gamma)$ (resp. $x^-(\gamma)$) the attractor null eigenvector of $\gamma$
(resp.$\gamma^{-1}$)
 such that $y_0(x^+(\gamma))=1$ (resp. $y_0 (x^-(\gamma))=1$).
For every $z\in\Kn{M}{n+1}$ we can write
\[
     z= a_+(\gamma,z)x^+(\gamma)+ a_-(\gamma,z)x^-(\gamma)+x^0(\gamma,z)+x^1(\gamma,z)
\]
with $a_+,a_-\in\mathbb R$, $x^0\in\ker(\gamma-1)$ and
$x^1\in\ort{\ker(\gamma-1)}\cap\ort{\langle x^+(\gamma),x^-(\gamma)}$.
In order to prove c) we now have to show the following lemma.
\begin{lemma}\label{cocycle stime}
Fix $r\in\overline{\dom}$.
For every $\gamma\in\Gamma$ put $z(\gamma)=\gamma_\tau r -r$. Then we have
\begin{eqnarray*}
   a_+(\gamma,z(\gamma))\geq 0\,\,\\
   a_-(\gamma,z(\gamma))\leq 0.
\end{eqnarray*}
Furthermore if $ a_+(\gamma,z(\gamma)) a_-(\gamma,z(\gamma))=0$ then
$x^0(\gamma,z(\gamma))=0$.
\end{lemma}
\Dim
We have seen in \ref{domain of dependence is proper} that $z$ is a cocycle.
Thus it is easy to see that
\begin{eqnarray*} 
    z(\gamma^k)=\sum_{i=0}^{k-1}\gamma^iz(\gamma)\\
    z(\gamma^{-k})=-\sum_{i=1}^{k}\gamma^{-i}z(\gamma).
\end{eqnarray*}
Let $\lambda>1$ such that $\gamma x^+(\gamma)=\lambda x^+(\gamma)$ and put
$a_+=a_+(\gamma, z(\gamma))$, $a_-=a_-(\gamma, z(\gamma))$,
$x^0=x^0(\gamma, z(\gamma))$ and $x^1=x^1(\gamma, z(\gamma))$. By using 
last formulas we obtain
\begin{eqnarray}\label{calculation1-cocycle stime}
 z(\gamma^k)=\frac{\lambda^k-1}{\lambda-1} a_+x^+(\gamma)\,
+\,\frac{1}{\lambda^k}\frac{\lambda^k-1}{\lambda-1}a_-x^-(\gamma)\,+\,kx^0+(\gamma^{k-1}+\ldots+\gamma+1)x^1;\nonumber\\
  z(\gamma^{-k})=-\frac{1}{\lambda^k}\frac{\lambda^k-1}{\lambda-1}
    a_+x^+(\gamma)\,-\,\frac{\lambda^k-1}{\lambda-1}a_-x^-(\gamma)\,-\,kx^0-\gamma^{-k}(\gamma^{k-1}+\ldots+\gamma+1)x^1.
\end{eqnarray}
Now notice that $W=\ort{\ker(1-\gamma)}\cap\ort{\langle
  x^+(\gamma),x^-(\gamma)\rangle}$ is a spacelike $\gamma$-invariant subspace
  and the application $(1-\gamma)|_W$ is invertible. Let us denote
  by $B_\gamma$ the map $(1-\gamma)|_W^{-1}$. Now it easy to see that
$(\gamma^{k-1}+\ldots+\gamma+1)x^1=(\gamma^k-1)B_\gamma x^1$ and so we have
  that  the set $\{(\gamma^{k-1}+\ldots+\gamma+1)x^1\}_{k\in\mathbb N}$ is
  contained in a compact of $W$.\par
Fix a timelike future directed vector $e$. Since $r\in\overline{\dom}$ there exists
  $K\in\mathbb R$ such that $\E{\alpha r}{e}\leq K$ for all
  $\alpha\in\Gamma$ and thus $\E{z(\alpha)}{e}\leq 2K$ for all
  $\alpha\in\Gamma$.
Now let us impose $\E{z(\gamma^k)}{e}\leq 2K$ for every $k\in\mathbb N$. Since
 $\{(\gamma^{k-1}+\ldots+\gamma+1)x^1\}_{k\in\mathbb N}$ is
  contained in a compact set so that there exists $K'$ such that
\begin{equation}\label{calculation2-cocycle stime} 
    \frac{\lambda^k-1}{\lambda-1} a_+\E{x^+(\gamma)}{e}
    +\frac{1}{\lambda^k}\frac{\lambda^k-1}{\lambda-1}a_-\E{x^-(\gamma)}{e}+\,k\E{x^0}{e}\leq
K'.
\end{equation}
Suppose $a_+<0$: by passing to the limit we have that the left expression in
 (\ref{calculation2-cocycle stime})
 tends to $+\infty$  (in fact notice that
$\E{x^+(\gamma)}{e}<0$).But this contradicts (\ref{calculation2-cocycle stime}) and so we
 have $a_+\geq 0$. An analogous
argument shows that $a_-\leq 0$.\par
Now suppose that $a_+=0$ (the case $a_-=0$ is analogous). Suppose $x^0\neq
0$. Since $x^0$ is spacelike then we can choose the
vector $e$ so that $\E{x^0}{e}>0$ but then the expression on the left in
(\ref{calculation2-cocycle stime}) tends to $+\infty$ and this contradicts
 (\ref{calculation2-cocycle stime}).  
\cvd
Now we can prove that $C$ has not any timelike support plane. By contradiction
 suppose that there exists a spacelike vector $v$ and $K\in\mathbb R$ such that
\[
    \E{\gamma_\tau p}{v}\leq K\qquad\qquad\textrm{ for all }\gamma\in\Gamma.
\]
Let $z(\gamma)=\gamma_\tau p-p$, so we have $\E{z(\gamma)}{v}\leq 2K$ for all
$\gamma\in\Gamma$.\par
We can fix $\gamma\in\Gamma$ such that $\E{x^+(\gamma)}{v}\geq 0$ and
$\E{x^-(\gamma)}{v}\geq 0$ (a such $\gamma$ exists because the limit set of
$\Gamma$ is the whole $\partial\Kn{H}{n}$). Now put
$a_+=a_+(\gamma,z(\gamma))$. Notice that $a_+\neq 0$: in fact if $a_+=0$
we have that $x^0(\gamma,z(\gamma))=0$. Then from (\ref{calculation1-cocycle stime}) it follows
that $z(\gamma^k)$ runs in a compact set for $k\geq 0$. But we know that the
action of $\Gamma_\tau$ on $C$ is properly discontinuous (in fact $C\subset
\dom$) and this is a contradiction.
Thus $a_+$ is positive. By using (\ref{calculation1-cocycle stime}) we easily see
that
$\E{z(\gamma^k)}{v}\rightarrow +\infty$ and this is a contradiction.\par
Finally we have that $C$ is a future ore past convex set. Since $C\subset
\dom\cap\Omega$ and these are future convex then $C$ is future convex
set. This concludes th proof of the step 2.\\

Now  \emph{we prove the step 3}.
Let $\Delta=\partial C$: since $C$
is $\Gamma_\tau$-invariant it follows that $\Delta$ is
$\Gamma_\tau$-invariant. Furthermore since $C$ is a convex set with interior
non empty  then $\Delta$ is a topological $n$-manifold.\par
We have to show that \emph{$\Delta/\Gamma_\tau$ is compact}. Let $r:\dom\rightarrow\partial \dom$ be the
retraction, since $\Delta\subset\dom$ we can define
\[
   f: \Delta\ni p\mapsto r(p)+\frac{1}{T(p)}(p-r(p))\in \widetilde S_1.
\]
Notice that $f(p)$ is the intersection of the timelike line $p+\mathbb
R(p-r(p))$ with the surface $\widetilde S_1$.
Clearly we have that $f$ is $\Gamma_\tau$-equivariant so induces a map
$\overline f:\Delta/\Gamma_\tau\rightarrow \widetilde S_a/\Gamma_\tau$. We have seen that
$\widetilde S_a/\Gamma_\tau$ is homeomorphic to $M$ so that it is sufficient to show that
$\overline f$ is a homeomorphism.\par
Since $C$ is future convex it is easy to see that $f$ is injective. On
the other hand fix $q\in \widetilde S_1$, we have that $q+\lambda
(q-r(q))\in C$  for $\lambda>>1$ (in fact fix $p_0\in C$ then $q+\lambda
(q-r(q))\in\fut(p_0)\subset C$ for $\lambda>>1$) and $q+\lambda (q-r(q))\notin
C$ for $\lambda<<0$. Thus there exists $\lambda_0$ such that $p=q+\lambda_0
(q-r(q))\in \Delta$. Since $r\left(q+\lambda(q-r(q))\right)=r(q)$ for all
$\lambda\in(-1,+\infty)$,  we get $r(p )=r(q)$ so that the lines $p+\mathbb R
(p-r(p))$ and $q+\mathbb R (q-r(q))$ coincide. Thus $q$ is the intersection of
the line $p+\mathbb R(p-r(p))$ with $\widetilde S_1$ so that $f(p)=q$.
It follows that
the map $f$ is surjective. By theorem of the the invariance of domain  we get that
$f$ is a homemorphism and so $\overline f$. This concludes the step 3.\\

Now \emph{we prove the step 4}. Notice that the function $T:\Delta\rightarrow\mathbb R$
and $T_\Omega:\Delta\rightarrow\mathbb R$ are $\Gamma_\tau$-invariant. Since
$\Delta/\Gamma_\tau$ is compact these functions are bounded on $\Delta$ 
so that there exists $a>0$ such that $T(x)<a$ and $T_\Omega(x)<a$ for every
$x\in\Delta$.
We have to show  that $\widetilde S_a$ and $\widetilde S_a^\Omega$ are contained in $C$. Let  $y\in \dom$ and suppose
$y\notin C$. We have that $y\in\pass(\Delta)$:  in
fact if $e$ is a timelike future directed vector the ray $y+\mathbb R_+ e$
intersects $C$ and so it intersects $\Delta$. Thus there exists
$y'\in\Delta\cap\fut(y)$ so that  we have $T(y)<T(y')<a$. Thus $\widetilde S_a\subset C$. An
analogous argument shows that $\widetilde S_a^\Omega\subset C$. This concludes
the proof of step
4 and the proof of theorem \ref{uniqueness}.
\cvd   
 
Clearly an analogous theorem holds for  $\Gamma_\tau$-invariant past complete
regular domain. So that
$\dom^-$ is the unique  $\Gamma_\tau$-invariant past complete regular domain. 
\begin{corollario}
If $\tau$ and $\sigma$ differs by a traslation then $\dom$ and $\mathcal
D_\sigma$ differs by a traslation.
Moreover we have that $\mathcal D_{-\tau}$ coincides with $-(\mathcal
D_\tau^-)$.
\end{corollario}
\Dim
Suppose that $\sigma_\gamma-\tau_\gamma=\gamma(x)-x$ it is easy to see that
$\dom+x$ is a $\Gamma_\sigma$-invariant future complete regular domain.\par
On the other hand notice that $-(\dom^-)$ is a future complete regular domain
and it is invariant by the action of $\Gamma_{-\tau}$.
\cvd
Let $Y_\tau:=\dom/\Gamma_\tau$  (resp. $Y^-_\tau:=\dom^-/\Gamma_\tau$). We get that
$Y_\tau$ and $Y_\sigma$ are isometric if and only if $\tau$ and $\sigma$
differs by a coboundary (i.e. the isometric class of $Y_\tau$ depends only on
the cohomology class of $\tau$). Notice that for $\tau=0$ the domain $\mathcal
D_0$
(resp. $\mathcal D_0^-$) coincides with $\fut(0)$ (resp. $\pass(0)$) so that $Y_0$ is
the Minkowskian cone $\mathcal C^+(M)$.\par
On the other hand notice that there exists a time-orientation reversing
isometry between $Y_{-\tau}$ and $Y_\tau^-$.
\begin{corollario}
Let $\widetilde F$ be a $\Gamma_\tau$-invariant complete spacelike hypersurface such that
the action of $\Gamma_\tau$ is free and properly discontinuous. 
Then $\widetilde F$ is contained in $\dom$ or in
$\dom^-$. In particular every timelike coordinate is proper on $\widetilde F$ and the
Gauss map has degree $1$. Furthermore $\widetilde F/\Gamma_\tau$ is homeomorphic to $M$.
\end{corollario}
\Dim
We know that $D(\widetilde F)$ is a $\Gamma_\tau$-invariant future or past
complete regular domain.
By theorem \ref{uniqueness} we get that 
 either $D(\widetilde F)=\dom$ or $D(\widetilde F)=\dom^-$.  Thus $\widetilde F$ is contained either
in $\dom$ or in $\dom^-$. Notice that this
implies that every timelike coordinate on $\widetilde F$ is proper.\par
Suppose $\widetilde F\subset\dom$. Consider  the map
\[
   \varphi:\widetilde F\ni x\mapsto r(x)+\frac{1}{T(x)}\left(x-r(x)\right)\in \widetilde S_1.
\]
It is easy to see that this map is $\Gamma_\tau$-equivariant and
injective. Furthermore since $\widetilde F$ is a Cauchy surface for $\dom$ (in
fact $\dom=D(\widetilde F)$) one easily
see that $\varphi$ is surjective. Thus it is an homeomorphism
$\Gamma_\tau$-equivariant. It follows that it induces an homeomorphism
$\overline\varphi:\widetilde F/\Gamma_\tau\rightarrow \widetilde
S_1/\Gamma_\tau\cong M$.\par
Finally consider the Gauss map of $\widetilde F$. It is a $\Gamma$-equivariant map
$G:\widetilde F\rightarrow \Kn{H}{n}$ (i.e. $G(\gamma_\tau p)=\gamma G(p)$). It is
easy to see that this map induces to the quotient a map $\overline G:\widetilde
F/\Gamma_\tau\rightarrow M$ which is a homotopical equivalence. Thus it has
degree $1$.
\cvd
Now we want to prove that $Y_\tau$ and $Y_\tau^-$ are the unique maximal
globally hyperbolic spacetimes with a compact spacelike Cauchy surface such
that the holonomy group is $\Gamma_\tau$.
We need the following remark which 
was stated by Mess in \cite{Mess} for the case $n=2$. However his proof runs
in every dimension.
\begin{corollario}
For every $\tau\in Z^1(\Gamma, \Kn{R}{n+1})$ the intersection $\dom\cap\dom^-$ is empty.
\end{corollario} 
\Dim 
It is easy to see that $\dom\cap\dom^-$ is a $\Gamma_\tau$ invariant
compact set.
Thus if it is not empty its barycenter $p$ is a fix point of $\Gamma_\tau$.  
It is straightforward to recognize that $\fut(p)$ and $\pass(p)$ are
respectively a $\Gamma_\tau$-invariant future and past complete domain of
dependence (notice that the cohomology class of $\tau$ vanishes). Hence
$\dom=\fut(p)$ and $\dom^-=\pass(p)$ so that their intersection is
empty.
\cvd
\begin{corollario}
There exists only two maximal globally hyperbolic  flat spacetimes
with compact spacelike Cauchy surface such that the  holonomy group is
$\Gamma_\tau$.
\end{corollario}
\Dim
Let $Y$ be a maximal globally hyperbolic flat spacetime with compact spacelike Cauchy
surface $N$ and holonomy group equal to $\Gamma_\tau$. We have to show that $Y$ isometrically embbeds in $Y_\tau$ or in
$Y^-_\tau$. It is sufficient to show that the developing map $D:\widetilde
Y\rightarrow\Kn{M}{n+1}$ is an embedding such that the image is contained
either in $\dom$ or in $\dom^-$.\par
Let $N$ be the spacelike Cauchy surface of $Y$. We know that $D:\widetilde
N\rightarrow\Kn{M}{n+1}$ is an embedding and the image $D(\widetilde N)$ is a
$\Gamma_\tau$-invariant surface such that the $\Gamma_\tau$-action on it is
free and properly discontinuous. Thus $D(\widetilde N)$ is a Cauchy surface of  $\dom$ or
$\dom^-$. It follows that $N$ is homeomorphic to $M$. In \cite{Andersson2} it
is shown that $Y$ is foliated by spacelike hypersurfaces so that
$D(Y)\subset\dom\cup\dom^-$. Since these domain are disjoint it follows that
$D(Y)$ is contained in one of them, say $\dom$ (the other case is analogous).\par
Consider the map $T_D:=T\circ D$ where $T$ is the CT of $\dom$: we have that
$T_D$ is a $\pi_1(Y)$-invariant regular map such that the level surface $\widetilde
N_a$ are $\pi_1(Y)$-invariant spacelike Cauchy surface. Thus $\widetilde
N_a/\pi_1(Y)\cong N$ is compact. It follows that  $D|_{\widetilde N_a}$ is an embedding.
Moreover let $p\in \widetilde N_a$ and $q\in\widetilde N_b$ with $a\neq b$, since we have that $T(D(p))=a$
and $T(D(q))=b$ it follows $D(p)\neq D(q)$. Thus the map $D$ is an  embedding of $Y$
in $\dom$. This map induces to the quotient the embedding $Y\rightarrow Y_\tau$.   
\cvd

\section{Continuous family of  Domains of Dependence}
We use the notation introduced in the previous sections. In particular
$\Gamma$ is a free-torsion co-compact and discrete subgroup of $\SOO^+(n,1)$
and $M=\Kn{H}{n}/\Gamma$. We have seen that there exists a well defined
corrispondence
\[
  \coom1(\Gamma, \Kn{R}{n+1})\ni [\tau]\mapsto [Y_\tau]\in\mathcal T_{Lor}(M)
\]
where $Y_\tau$ is the quotient of the unique $\Gamma_\tau$-invariant future
complete regular domain $\mathcal D_\tau$ by the action of the deformed group
$\Gamma_\tau$ (we recall that $\mathcal T_{Lor}(M)$ is the Teichm\"uller
space of the globally hyperbolic flat Lorentzian structure on $\mathbb R\times
M$ with a spacelike Cauchy surface).
In this section we shall show that this correspondence is continuous.\par 
More precisely, we shall prove that for every boundend neighborhood $U$ of $0$
in $Z^1(\Gamma, \Kn{R}{n+1})$ there is a continuous map
\[
   dev:U\times\big(\mathbb R_+\times\widetilde M\big)\rightarrow\Kn{M}{n+1}
\]
such that for every $\sigma\in U$ the map $dev_\sigma=dev(\sigma,\cdot)$ is a
developing map with holonomy  equal to $\rho_\sigma$ (notice that
$\pi_1(\mathbb R_+\times M)=\pi_1(M)=\Gamma$) and it is
a homeomorphism onto $\mathcal D_\sigma$.\par
We start with the map
\begin{equation}\label{definition of dev0}
    dev^0:U\times\widetilde M\rightarrow\Kn{M}{n+1}
\end{equation}
constructed in theorem \ref{existence of structure}. We know that
$dev^0_\sigma$ is an embedding onto a $\Gamma_\sigma$-invariant strictly convex spacelike
hypersurface  and the map $dev^0_\sigma$ is
$\Gamma$-equivariant in the following sense
\[
    dev^0_\sigma(\gamma x)=\gamma_\sigma dev^0_\sigma(x).
\]
For every $\sigma\in U$ let $\widetilde F_\sigma=dev^0_\sigma(\widetilde M)$. Now fix
an orthonormal affine coordinates system $(y_0,\ldots,y_n)$. We know that for
every $\sigma\in U$ there exists a convex function
$\varphi_\sigma:\{y_0=0\}\rightarrow\mathbb R$ such that $\widetilde F_\sigma$ is
the graph of $\varphi_\sigma$.
The first remark is that $\varphi_\sigma$ is a continuous function of
$\sigma$. More exactly
let $(\sigma_k)_{k\in\mathbb N}$ be  a sequence in $U$ which converges to
$\sigma$ in $U$ then $\varphi_{\sigma_k}$ converges to $\varphi_{\sigma}$ in
compact open topology of the plane $\{y_0=0\}$.\par
We know that for every $\sigma\in U$ the domain $\mathcal D_\sigma$ is the domain
of dependence of $\widetilde F_\sigma$. Now let
$\psi_\sigma:\{y_0=0\}\rightarrow\mathbb R$ such that $\partial \mathcal D_\sigma$ is
the graph of $\psi_\sigma$. We want to prove that $\psi_\sigma$ is a
continuous function of $\sigma$
\begin{prop}\label{functions vary with continuity}
Let $(\tau_k)_{k\in\mathbb N}$ be a sequence in $U$ which converges to
$\tau\in U$. Then the maps $\psi_{\tau_k}$ converge to $\psi_\tau$ in
compact open topology.
\end{prop}
\Dim
First let us show that the family
$\{\psi_{\tau_k}:\{y_0=0\}\rightarrow\mathbb R\}_{k\in\mathbb N}$ 
is locally bounded and equicontinuous.
Since two points on $\partial \mathcal D_{\tau_k}$ are not chronological related it
follows that the maps $\psi_{\tau_k}$ are $1$-Lipschitzian so that they form
an equicontinuous family. We have to prove that they are locally bounded.
Now we know that $\varphi_{\tau_k}\rightarrow\varphi_{\tau}$ and since
$\widetilde F_\tau$ is contained in $\mathcal D_\tau$ we have  
$\psi_{\tau_k}\leq\varphi_{\tau_k}$. On the other hand we can construct
a family of future strictly convex $\Gamma_\tau$-invariant spacelike hypersurfaces $\left\{\widetilde
F^-_\sigma\right\}_{\sigma\in U}$ which ``vary'' continuously. So let
$\varphi^-_\sigma:\{y_0=0\}\rightarrow\mathbb R$ such that $\widetilde F^-_{\sigma}$
is the graph of $\varphi^-_\sigma$. We have that
$\varphi^-_{\tau_k}\rightarrow\varphi^-_\tau$.
But the domain of dependence of $F^-_{\tau_k}$ is $\mathcal D^-_{\tau_k}$ and since
it is disjoint from $\mathcal D_{\tau_k}$ we deduce that
\[
     \varphi^-_{\tau_k}\leq\psi_{\tau_k}\leq\varphi_{\tau_k}.
\]
Thus $\{\varphi^-_{\tau_k}\}_{k\in\mathbb N}$ and
$\{\varphi_{\tau_k}\}_{k\in\mathbb N}$ are convergent and hence locally
bounded. It follows that $\psi_{\tau_k}$ is locally bounded too.\par
Now it remains to prove that $\psi_{\tau_{k}}\rightarrow \psi_\infty$ then
$\psi_\infty=\psi_{\tau}$.
We have that $\psi_\infty$ is a convex function and the graph $S$ of $\psi_\infty$ is
$\Gamma_\tau$-invariant. Furthermore since $\psi_\infty$ is $1$-Lipshitz
function then $S$ has only non-timelike support plane. Hence $\fut(S)$ is the
future of the graph of $\psi_\infty$ and it is a future convex set. 
It is easy to see that $\fut(S)$ is a future complete regular domain. Since it
is $\Gamma_\tau$-invariant by theorem \ref{uniqueness} we get
$\fut(S)=\dom$. Thus $\graph\psi_\infty=\partial\dom$ and so
$\psi_\infty=\psi_\tau$.
\cvd
Fix $K$ a compact subset of $\mathcal D_\tau$. The last proposition implies that
there exists a neighborhood $V$ of $\tau$ (which depends on $K$) 
such that if $\tau\in V$ then $K\subset \mathcal D_\tau$. Thus for every $\tau\in V$
the cosmological time $T_\tau$ , the normal field $N_\tau$ and the retraction
$r_\tau$ of the domain $\mathcal D_\tau$ are maps defined over $K$. The following
propositions show that these maps change continuously on $K$ when $\tau$
varies in $V$.
\begin{prop}\label{CT is natural}
Let $\{\tau_k\}$ be a sequence of cocycles which belongs to the neighborhood
$V$ of $\tau$. Then the sequences of function $\{T_{\tau_k}|_K\}$ converges
uniformly to $T_\tau|_K$ on $K$.
\end{prop}
In order to prove the proposition we need the following technical lemma.
\begin{lemma}\label{technical}
For $C\in\mathbb R$ and for every cocycle $\sigma$ let
\[ 
   K_C(\sigma)=\{x\in\{y_0=0\}|\psi_\sigma(x)\leq C\}
\]
($\psi_\sigma$ is the function defined over the horizontal
plane such that $\partial \mathcal D_\sigma$ is the graph of such function)
Then for every $C\in\mathbb R$ and $\eps>0$ there exists $k_0\in\mathbb N$
such that
\[
     K_{C-\eps}(\tau)\subset K_C(\tau_k)\subset
     K_{C+\eps}(\tau)\qquad\textrm{for all }
     k\geq k_0.
\]
For every cocycles $\sigma$ put $M(\sigma)$ the minimum of the function 
$\psi_\sigma:\{y_0=0\}\rightarrow\mathbb R$. 
Then the sequence $M(\tau_k)$ converges to $M(\tau)$.
\end{lemma} 
\Dim
Notice that  $K_C(\sigma)$ is a convex compact set, moreover if $C>M(\sigma)$
then $K_C(\sigma)$ has non-empty interior part and $\partial K_C(\sigma)$ is the level
set $\{x|\psi_\sigma(x)=C\}$.\par
Now let $M=M(\tau)$. First let us show the first statement for $C>M$. 
Fix  $\eps>0$ and let $k_0\in\mathbb N$ such that
$||\psi_\tau-\psi_{\tau_k}||_{\infty,K_{C+\eps}(\tau)}<\frac{\eps}{2}$ for all
$k\geq k_0$ . 
Clearly $K_{C-\eps}(\tau)\subset K_C(\tau_k)$ for all $k\geq k_0$.
Now let $x\notin K_{C+\eps}(\tau)$. I claim that $\psi_{\tau_k}(x)\geq
C+\frac{\eps}{2}$ for all $k\geq k_0$ and this proves the other inclusion.\par
Let $k>k_0$ and fix  $x_0$  such that $\psi_\tau(x_0)=M$. Consider the map
$c(t)=\psi_{\tau_k}(x_0+t(x-x_0))$ for $t\in[0,1]$. Let $t_0$ such that
$y_0+t_0(x-x_0)\in\partial K_{C+\eps}(\tau)$. We have that
$c(0)\leq M+\frac{\eps}{2}$ and $c(t_0)\geq C+\frac{\eps}{2}$. By
imposing that $c(t_0)\leq (1-t_0)c(0)+t_0c(1)$ we get that
$\psi_{\tau_k}(x)=c(1)\geq C+\frac{\eps}{2}$.\par
Now suppose $C<M$. Fix $k_0$ such that $K_{M+1}(\tau_k)$ is contained in
$K_{M+2}(\tau)$  and
$||\psi_\tau-\psi_{\tau_k}||_{\infty,K_{M+2}(\tau)}<\frac{M-C}{2}$ for all
$k>k_0$. Then it turns out that $K_C(\tau_k)=\varnothing$ for all $k>k_0$.\par
By this fact it turns out that $M(\tau)\leq\underline{\lim}_{k\rightarrow+\infty} M(\tau_k)$.
On the other hand since $\psi_{\tau_k}$ converges to $\psi_{\tau}$ one easily see that
$M(\tau)\geq\overline{\lim}_{k\rightarrow+\infty} M(\tau_k)$.
\cvd
Now we can prove the proposition \ref{CT is natural}.\\
\Dim
Let $M$ be the minimum of $\psi_\tau$. By lemma \ref{technical}  there exists $k_0$
such that $\psi_{\tau_k}>M-1$ for all $k\geq k_0$. Now notice that the set
$\cpass(K)\cap\{y_0\geq M-1\}$ is  compact set and  let $H$ be the projection
  of it  on the horizontal plane $\{y_0=0\}$.
 Fix $\eps>0$ and let $k(\eps)$
  such that  $||\psi_\tau-\psi_{\tau_k}||_{\infty,H}<\frac{\eps}{2}$ for
  $k\geq k(\eps)$.\par
Let  $p\in K$ and $r$ be the projection of $p$ on $\partial
\mathcal D_\tau$: notice that $r\in \cpass(K)\cap\{y_0\geq M-1\}$. Now if
$k\geq k(\eps)$ then
$r+\eps\frac{\partial}{\partial y_0}$ belongs to $\mathcal D_{\tau_k}$ so that
\[
     T_{\tau_k}(p)>\sqrt{-\E{p-r+\eps \frac{\partial}{\partial y_0}}{p-r+\eps \frac{\partial}{\partial y_0}}}.
\]
Now notice that
$-\E{p-r+\eps \frac{\partial}{\partial y_0}}{p-r+\eps \frac{\partial}{\partial
    y_0}}=-T_{\tau}(p)^2-\eps^2+2\eps(p-r)_0$.
By the compactness of $\cpass(K)\cap\{y_0\geq M\}$ there exists a constant $C$
such that
\[
     T_{\tau_k}(p)>\sqrt{T_\tau(p)^2+\eps^2-2C\eps}
\]
Now fix $\eta>0$, we can choose $\eps>0$ so that $2\eps C-\eps^2\leq
\eta^2$. So we have that
$T_{\tau_k}(p)>T_\tau(p)-\eta$ for all $k\geq k(\eps)$ and $p\in K$.\par 
On the other hand we have that the projection $r_k(p)$ of $x$ on $\partial
\mathcal D_{\tau_k}$ lies in $\cpass(K)\cap\{y_0\geq M-1\}$ and so
the same argument shows that $T_\tau(p)>T_{\tau_k}(p)-\eta$ for all
$k>k(\eps)$ and $p\in K$.
\cvd  
Let $\psi_\tau^a:\{y_0=0\}\rightarrow\mathbb R$ such that
$\graph(\psi_\tau^a)$ is the CT level surface $\widetilde S_a(\tau):=T_\tau^{-1}(a)$. The proposition
\ref{CT is natural} implies that these maps are continuous functions of $\tau$.
\begin{corollario}\label{CT surface varies with continuity}
If $\tau_k\rightarrow\tau$ then $\psi_{\tau_k}^a\rightarrow\psi_\tau^a$
in the compact-open topology of $\{y_0=0\}$.
\end{corollario}
\Dim
Since the maps $\psi_{\tau_k}^a$ are $1$-Lipschitz 
$\{\psi_{\tau_k}^a\}_{k\in\mathbb N}$ is an equicontinuous family.
On the other hand notice that 
\[
    \psi_{\tau_k}<\psi_{\tau_k}^a<\psi_{\tau_k}+a.
\]
Since $\{\psi_{\tau_k}\}$ is locally bounded it follows that $\{\psi_{\tau_k}^a\}_{k\in\mathbb N}$
is bounded.
From proposition \ref{CT is natural} it follows that if $\psi_{\tau_k}^a$
converges then the limit is $\psi_\tau^a$.
\cvd
Now let us prove the retraction $r_\tau$ and the normal field $N_\tau$ are
continuous functions of $\tau$.
\begin{prop}\label{Gauss map is natural}
Let $\tau_k\rightarrow\tau$ be as
above. Let $r_k$ and $N_k$ be respectively the retraction and the
normal field of $\mathcal D_{\tau_k}$.
Now fix  a compact subset $K$ of $\mathcal D_\tau$ . The maps $r_k|_K$ and $N_k|_K$
converge in the compact open topology of $K$ respectively to the retraction $r$ and
to the normal field $N$ of the domain $\mathcal D_\tau$.
\end{prop}
\Dim
Let $M$ be the minimum of the map $\psi_\tau$ and fix $k_0$ such that 
$\psi_{\tau_k}\geq M+1$ for all $k\geq k_0$. In particular
$r_k(p)\in\cpass(K)\cap\{y_0\geq M-1\}$ for all $p\in K$ and $k>k_0$. 
Since $\cpass(K)\cap\{y_0\geq M-1\}$ is compact we can choose $C$ such that
$||p-r_k(p)||\leq C$  for all $p\in K$ and $k\geq k_1$ ($||\cdot||$ is the
euclidean norm).  
On the other hand because of the proposition \ref{CT is natural} we can choose
$k_1>k_0$ such that 
\[
           T_{\tau_k}(p)\geq\alpha>0\qquad\textrm{for all }p\in K\textrm{ and } k\geq
          k_1.
\]
Thus the image $N_k(K)$ is contained in the set $H=\{x\in\Kn{H}{n}|\,
||x||\leq \frac{C}{\alpha}\}$ for all $k\geq k_1$.
This is a compact set of $\Kn{H}{n}$ so that the family of functions
$\{N_k|_K\}$ is bounded. \par
In order to show that $N_k|_K\rightarrow N|_K$ it is sufficient to prove that 
$N_k(p_k)\rightarrow N(p)$ for all convergent sequences $p_k\rightarrow
p$. Since $N_k(p_k)$ runs in a compact set of $\Kn{H}{n}$ we can suppose
that $N_k(p_k)$ converges to a timelike vector $v$. Let $a=T(p)$, in order to
show that $N(p)=v$ it is sufficient to prove that $p+\ort v$ is a support
plane for the surface $\widetilde S_a=T^{-1}(a)$ i.e. we have to prove the following
inequality
\begin{equation}\label{fff}
    \E{q}{v}\leq\E{p}{v}\qquad\textrm{ for all } q\in \widetilde S_a.
\end{equation}
Fix $q\in \widetilde S_a$ and let $q=(\psi_\tau^a(y),y)$. 
Let $a_k=T_k(p_k)$ and we consider the sequences 
\begin{eqnarray*}
  q_k:=(\psi_{\tau_k}^{a_k}(y),y);\\
  q'_k:=(\psi_{\tau_k}^a(y),y).
\end{eqnarray*}
 By corollary \ref{CT surface varies with
  continuity} we have that $q'_k\rightarrow q$. On the other hand it turns out
that
$||q_k-q'_k||\leq |a_k-a|$ so that  $q_k\rightarrow q$.
We know that $\E{q_k}{N_k(p)}\leq\E{p_k}{N_k(p)}$: by passing to the limit the inequality
(\ref{fff}) follows.\par
Since $r_k+T_kN_k=id$ we get that $r_k|_K\rightarrow r|_K$ uniformly.
\cvd
\begin{corollario}
Let $K$ as above. Then the cosmological times $T_{\tau_k}$ converge to
$T_\tau$ in the $\mathrm C^1$-topology of $\mathrm C^1(K)$.
\end{corollario}
\cvd
Now let us go back to the original problem. 
\begin{teo}\label{continuity of dom}
For every boundend neighborhood $U$ of $0$ in
$Z^1(\Gamma, \Kn{R}{n+1})$ there exists a continuous map
\[
       dev:U\times\big(\mathbb R_+\times\widetilde M\big)\rightarrow\Kn{M}{n+1}
\]
such that for every $\sigma\in U$
\begin{enumerate}
\item
 $dev_\sigma$ is a developing map whose holonomy is $\rho_\sigma$;
\item
 $dev_\sigma$ is a homeomorphism with $\mathcal D_\sigma$.
\end{enumerate}
\end{teo}
\Dim
Let $dev^0:U\times\Kn{H}{n}\rightarrow\Kn{M}{n+1}$ be the map defined in
(\ref{definition of dev0}). Now fix $\sigma\in U$ , $x\in\Kn{H}{n}$ and $t>0$. 
Consider the timelike geodesic $\gamma$  
in $\mathcal D_\sigma$ which passes through $dev^0_\sigma(x)$ and has the direction of
the normal field in $dev^0_\sigma(x)$. Now put $dev(\sigma,t,x)$ the point on
$\gamma$ with CT equal to $t$:
\[
   dev(\sigma,t,x)=r_\sigma(dev^0_\sigma(x))\,+\,t N_\sigma(dev^0_\sigma(x)).
\]
Clearly $dev$ satisfies the three properties  required. On the other hand by
propositions \ref{CT is natural} and \ref{Gauss map is natural} one easily see
that it is continuous.
\cvd
\begin{oss}
\emph{
With the proof of theorem \ref{continuity of dom} the proof of
theorem \ref{main} is complete. In the following section we shall prove  theorem
\ref{main2}.
}\end{oss}
\begin{oss}
\emph{
Notice that the coordinate $t$ on $\mathbb R_+\times\widetilde M$ coincides with
the pull back of the cosmological times by the map $dev_\sigma$, i.e.}
\[
  T_\sigma(dev_\sigma(t,x))=t
\]
\emph{
On the other hand notice that $r_\sigma(dev_\sigma(t,x))$ and
$N_\sigma(dev_\sigma(t,x))$ depend only on the $x$ coordinate. Thus there are
well defined functions}
\begin{eqnarray*}
  \mathbf r_\sigma:\widetilde M\rightarrow\Sigma_\tau\\
  \mathbf N_\sigma:\widetilde M\rightarrow\Kn{H}{n}
\end{eqnarray*}
\emph{
such that $r_\sigma(dev_\sigma(t,x))=\mathbf r_\sigma(x)$ and
$N_\sigma(dev_\sigma(t,x))=\mathbf N_\sigma(x)$.
}\end{oss}

The map $dev$ is only continuous. But we can smooth this map to
obtain a $\mathbb C^\infty$ map $dev'$ which verifies the properties required
in the above theorem. In fact it is easy to see that  we can perturb the
normal field $N_\sigma$ on $\mathcal D_\sigma$ to obtain a $\Gamma_\tau$-invariant $\mathrm
C^\infty$ timelike vector field $V_\sigma$. By considering the restriction of the
flow of this vector field on the Cauchy surface $\widetilde F_\sigma$ we obtain a
$\mathrm C^\infty$ developing map $dev'_\sigma$ which verifies the properties
1. and 2. of theorem \ref{continuity of dom}.\par
We can construct the field $V_\sigma$
such that $V_0$ coincides with $N_0$ and $V_\sigma$ ``varies continuously'' 
with $\sigma$ in the following sense: for every convergent sequence $\tau_k\rightarrow\tau$ and for every
open set $K\subset \mathcal D_\tau$ the fields $V_{\tau_k}|_K$ converge to
$V_\tau|_K$ in the $\mathrm C^\infty$ topology of $K$.
In this way it is easy to see that the map $dev'(\sigma,x)=dev'_\sigma(x)$ is
$\mathrm{C}^\infty$ map which verifies the properties required in the theorem.

\section{Gromov Convergence of the CT-Level Surfaces}\label{Convergence-sec}

Let us summarize what we have seen until now. We have fixed a closed
hyperbolic $n$-manifold $M$ and we have identified $\pi_1(M)$ with 
a free-torsion co-compact discrete subgroup of $\SOO^+(n,1)$, say $\Gamma$, such that
$M=\Kn{H}{n}/\Gamma$. Given a cocycle
$\tau\in Z^1(\Gamma,\Kn{R}{n+1})$ we have considered the deformation
$\Gamma_\tau$ of $\Gamma$. We have proved that there exists a unique
$\Gamma_\tau$-invariant future complete regular
domain $\dom$  such that the action is free and properly discontinuous and the
quotient is a globally hyperbolic spacetime diffeomorphic to $\mathbb
R_+\times M$. The domain $\dom$ is provided with a $\Gamma_\tau$-invariant regular cosmological time
$T$ which is a $\mathrm C^1$-submersion. Moreover
there exists a retraction map $r:\dom\rightarrow\Sigma$ onto the singularity in
the past and a normal field $N:\dom\rightarrow\Kn{H}{n}$ 
which is up to the sign the Lorentzian gradient of the cosmological time $T$. 
The level surfaces $\widetilde S_a=T^{-1}(a)$ are $\mathrm C^1$ spacelike hypersurfaces, so that there is a
natural path distance $d_a$ on it. Since $\widetilde S_a$ is $\Gamma_\tau$-invariant and
$\widetilde S_a/\Gamma_\tau\cong M$ is compact by Hopf-Rinow theorem we get that $d_a$ is
a complete distance. In this section we study the metric properties of the
surface $\widetilde S_a$ and in particular the asymptotic behaviour for
$a\rightarrow+\infty$ and for $a\rightarrow 0$.\par
In the previous section we have constructed a developing map
\[
  dev_\tau:\mathbb R_+\times\widetilde M\rightarrow\Kn{M}{n}
\]
such that $dev_\tau(\{a\}\times\widetilde M)=\widetilde S_a$ and
there exist well defined maps
\begin{eqnarray*}
\mathbf r:\widetilde M\rightarrow\Sigma\\
\mathbf N:\widetilde M\rightarrow\Kn{H}{n}
\end{eqnarray*}
such that $r(dev_\tau(t,x))=\mathbf r(x)$ and $N(dev_\tau(t,x))=\mathbf N(x)$.
By taking the pull-back of the distance $d_a$ on $\widetilde S_a$  we get a
family of distance $\delta_a$ on $\widetilde M$ such that $\pi_1(M)(=\Gamma)$ acts by
isometry on $(\widetilde M, \delta_a)$.\par
The principal results of this section are the following propositions.
\begin{prop}\label{convergence in the future}
For all $x,y\in\widetilde M$
\[
   \lim_{a\rightarrow+\infty}\frac{\delta_a(x,y)}{a}=d_{\mathbb H}(\mathbf
   N(x), \mathbf N(y))
\]
where $d_{\mathbb H}$ is the distance of $\Kn{H}{n}$. Moreover 
the maps $a^{-1}\delta_a$ converge in the compact open topology of
$\mathrm C(\widetilde M\times\widetilde M)$ to the map $(x,y)\mapsto d_{\mathbb
  H}(\mathbf N(x),\mathbf N(y))$.
\end{prop} 
\cvd
\begin{prop}\label{convergence in the past}
There exists a natural distance $d_\Sigma$ on $\Sigma$ such that 
\[
\lim_{a\rightarrow 0}\delta_a(x,y)=d_\Sigma(\mathbf r(x),\mathbf
r(y))\qquad\textrm{ for all }x\in\widetilde M.
\]
Moreover the maps $\delta_a$ converge in the compact open topology of $\mathrm
C(\widetilde M\times\widetilde M)$ to the map $(x,y)\mapsto d_\Sigma(\mathbf r(x),
\mathbf r(y))$.
\end{prop}
\cvd
We shall see that proposition \ref{convergence in the future} implies that the
action of $\Gamma_\tau$ on $\widetilde S_a$ converge for $a\rightarrow +\infty$ in the Gromov sense to the action of
$\Gamma_\tau$ on $\Kn{H}{n}$. For the  $a\rightarrow 0$ we can deduce by
proposition \ref{convergence in the past} 
only the convergence of the spectrum of the action of $\Gamma_\tau$ on $\widetilde S_a$
to the spectrum of the action of $\Gamma_\tau$ on $\Sigma$.\\

We start showing that  $\{\delta_a\}_{a>0}$ and
$\{a^{-1}\delta_a\}_{a>0}$ are respectively increasing and decreasing functions
of $a$.
\begin{lemma}\label{normal field is lipschitz}
Fix a Lipschitz path $c:[0,1]\rightarrow \widetilde S_a$, then the paths 
$N(t)=N(c(t))$ and $r(t)=r(c(t))$
are differentiable almost everywhere and we have
\begin{eqnarray*}
  N(t)=N(0)+\int_0^t\dot N(s)\mathrm{d}
  s;\\
  r(c(t))=r(c(0))+\int_0^t\dot r(s)\mathrm{d}
  s.
\end{eqnarray*}
Moreover we have that $\dot N(t)$ and $\dot r(t)$ lie into $T_{c(t)}\widetilde S_a$ (so
they are spacelike) and $\E{\dot N(t)}{\dot r(t)}>0$ almost everywhere.
\end{lemma}
\Dim
In order to prove the first statement it is sufficient to show 
that the maps $N:\widetilde S_a\rightarrow\Kn{H}{n}\subset\Kn{M}{n+1}$ and
$r:\widetilde S_a\rightarrow \Sigma\subset\Kn{M}{n+1}$ are locally Lipschitz with respect
to the \emph{euclidean distance} $d_E$ of $\Kn{M}{n+1}$. Since $p=r(p)+aN(p)$
it is sufficient to show that $N$ is locally Lipschitz.\par
Fix a compact $K\subset \widetilde S_a$ and let $H=N(K)$: since $H$ is compact there
exists a constant $C$ such that 
\[
   d_E(x,y)=||x-y||\leq C (\E{x-y}{x-y})^{1/2}.
\] 
On the other hand by inequalities (\ref{fundamental inequality}) 
we have that
\[
(\E{N(p)-N(q)}{N(p)-N(q)})^{1/2}\leq\frac{1}{a}(\E{p-q}{p-q})^{1/2}.
\]
Since
$\E{p-q}{p-q}\leq||p-q||^2$ we deduce that $||N(p)-N(q)||\leq
\frac{C}{a}||p-q||$ for all $p,q\in K$.\par
Now notice that $N(t)$ is a path in $\Kn{H}{n}$ so that $\dot N(t)\in
T_{N(t)}\Kn{H}{n}=T_{c(t)}\widetilde S_a$.
Since $\dot r(t)=\dot c(t)-a\dot N(t)$ we have that $\dot r(t)\in T_{c(t)}\widetilde S_a$
almost everywhere. Finally by inequalities (\ref{fundamental inequality})
we have that $\E{N(t+h)-N(t)}{r(t+h)-r(t)}\geq 0$. Thus we easily deduce that
$\E{\dot N(t)}{\dot r(t)}\geq 0$.
\cvd
\begin{lemma}\label{distance on S_a}
For all $x,y\in\widetilde M$ and for all $a<b$ we have
\begin{eqnarray*}
\delta_a(x,y)\leq\delta_b(x,y)\\
d_{\mathbb H}(\mathbf N(x),\mathbf
N(y))\leq\frac{1}{b}\delta_b(x,y)\leq\frac{1}{a}\delta_a(x,y).
\end{eqnarray*}
\end{lemma}
\Dim
For $t>0$ let $p_t=dev_\tau(t,x)$ and $q_t=dev_\tau(t,y)$. Let
$c_b:[0,1]\rightarrow \widetilde S_b$ be a minimizing-lenght geodesic path between $p_b$
and $q_b$ 
Consider $r(t)= r(c_b(t))$ and $N(t)=N(c_b(t))$ and let $c:[0,1]\rightarrow \widetilde S_a$ be the path defined by
the rule $c(t):=r(t)+aN(t)$. We have that $c$ is a rectifiable arc between
$p_a$ and $q_a$ so that the lenght of $c$ is greater than the distance
$\delta_a(x,y)$.
Now notice that $\dot c_b(t)=\dot c(t)+(b-a)\dot N(t)$. By lemma \ref{normal
  field is lipschitz} we get that $\E{\dot c_b(t)}{\dot c_b(t)}\geq\E{\dot
  c(t)}{\dot c(t)}$. This proves that the lenght of $c_b$ is greater than the
lenght of $c$ and so the first inequality follows. Now we shall prove the
second one.  \par
Let $c:[0,1]\rightarrow \widetilde S_a$ be a Lipschitz path and let $N(t)=N(c(t))$. By lemma \ref{normal field
  is lipschitz} we have that $a^2\E{\dot N(t)}{\dot N(t)}\leq\E{\dot c(t)}{\dot
    c(t)}$. By this inequality it follows that $d_{\mathbb
    H}(N(p),N(q))\leq\frac{1}{a} d_a(p,q)$ for all $p,q\in \widetilde S_a$.
On the other hand let $x,y\in\widetilde M$ and $c_a:[0,1]\rightarrow \widetilde S_a$ the $d_a$-minimizing geodesic
between $p_a=dev_\tau(a,x)$ and $q_a=dev_\tau(a,y)$. Let
$c(t)=c_a(t)+(b-a)N(t)$ (where $N(t)=N(c_a(t))$): the endpoints of this path are
$p_b$ and $q_b$ so that
\[
   \frac{\delta_b(x,y)}{b}=\frac{d_b(p_b,q_b)}{b}\leq\frac{1}{b}\int_0^1 \E{\dot
   c(t)}{\dot c(t)}^{1/2}\mathrm dt.
\]
By lemma \ref{normal field is lipschitz} we know that $a^2\E{\dot N(t)}{\dot
  N(t)}\leq\E{\dot c_a(t)}{\dot c_a(t)}$ so that
\begin{eqnarray*}
  \E{\dot c(t)}{\dot c(t)}^{1/2}\leq\E{\dot c_a(t)}{\dot
  c_a(t)}^{1/2}+(b-a)\E{\dot N(t)}{\dot N(t)}^{1/2}\leq\\
  \E{\dot c_a(t)}{\dot c_a(t)}^{1/2}+\frac{(b-a)}{a}\E{\dot c_a(t)}{\dot
  c_a(t)}^{1/2}= \frac{b}{a}\E{\dot c_a(t)}{\dot c_a(t)}^{1/2}.
\end{eqnarray*}
Thus we  have
\[
\frac{\delta_b(x,y)}{b}\leq\frac{1}{a}\int_0^1(\E{\dot
   c_a(t)}{\dot c_a(t)})^{1/2}\mathrm dt=\frac{\delta_a(x,y)}{a}.
\]
\cvd
Now we can prove proposition \ref{convergence in the future}.\\

\emph{Proof of proposition \ref{convergence in the future}:} 
Let $x,y\in\widetilde M$:
by lemma \ref{distance on S_a} we have that $\frac{1}{a}\delta_a(x,y)$
is decreasing with respect $a$ so that there exists
\[
   \delta_\infty(x,y)=\lim_{a\rightarrow +\infty}\frac{1}{a}\delta_a(x,y).
\]
Let us show that $a^{-1}\delta_a$ converges to $\delta_\infty$ in the
compact open topology of $\widetilde M$. Since $a^{-1}\delta_a\leq\delta_1$  the family
$\{a^{-1}\delta_a|_K\}_{a>1}$ is locally bounded. On the other hand by
triangular inequality we have for $a>1$
\[
  |a^{-1}\delta_a(x,y)-a^{-1}\delta_a(x',y')|\leq a^{-1}\delta_a(x,x')+
   a^{-1}\delta_a(y,y')\leq
   \delta_1(x,x')+\delta_1(y,y').
\]
Thus the family $\{a^{-1}\delta_a|_K\}_{a>1}$ is equicontinuous. By these
remarks we easily get that $a^{-1}\delta_a\rightarrow\delta_\infty$ in the
compact open topology of $\widetilde M\times \widetilde M$.\par
Clearly $\delta_\infty$ is a pseudo-distance on $\widetilde M$. We claim that
$\delta_\infty(x,y)=0$ if and only if $\mathbf N(x)=\mathbf N(y)$. In fact
from lemma \ref{distance on S_a} it follows that
\[
    d_{\mathbb H}(\mathbf N(x),\mathbf N(y))\leq\frac{\delta_a(x,y)}{a}
\]
so that $d_{\mathbb H}(\mathbf N(x),\mathbf N(y))\leq \delta_\infty(x,y)$. Thus if
$d_\infty(x,y)=0$ then $\mathbf N(x)=\mathbf N(y)$. On the other hand if
$\mathbf N(x)=\mathbf N(y)$ the
segment $[\mathbf r(x), \mathbf r(y)]$ is contained in $\Sigma$ and one easily
see that $\frac{\delta_a(x,y)}{a}=\frac{1}{a}(\E{\mathbf r(y)-\mathbf r(x)}{\mathbf r(y)-\mathbf
  r(x)})^{1/2}$. By passing to the limit we obtain that
$\delta_\infty(x,y)=0$.\par
It follows that there exists a distance $d$ on $\Kn{H}{n}$ such that
\[
   \delta_\infty(x,y)=d(\mathbf N(x),\mathbf N(y))
\]
In the last part of this proof we shall show that $d=d_{\mathbb H}$. We
already know that $d_{\mathbb H}\leq d$.\par 
By using theorem \ref{uniqueness} it is easy to see that 
$\mathcal D_{\frac{\tau}{a}}=\frac{1}{a}\dom$. Consider the map
\[
  f:\dom\ni p\mapsto \frac{p}{a}\in\mathcal D_{\frac{\tau}{a}}.
\]
Let $p\in\dom$ and $c_p$ be the Lorentz-lenght maximizing timelike geodesic of
$\dom$ with future-endpoint equal to $p$.
Then $f(c_p)$ is a Lorentz-lenght maximizing timelike geodesic of 
$\mathcal D_{\frac{\tau}{a}}$. Since the lenght of $f(c_p)$ is $a^{-1}T_\tau(p)$
we get
\[
  T_{\frac{\tau}{a}}(\frac{p}{a})=\frac{T_\tau(p)}{a}.
\]
Thus $\frac{1}{a}\widetilde S_a$ is the CT-level
surface $\widetilde S_1(\frac{\tau}{a})=T_{\frac{\tau}{a}}^{-1}(1)$. Moreover the distance $a^{-1}\delta_a$ is the
pull-back of the natural path-distance on $\widetilde S_1(\frac{\tau}{a})$.\par
First we claim that $\lim_{a\rightarrow+\infty}\frac{p_a}{a}=N(x)$ and
$\lim_{a\rightarrow+\infty}\frac{q_a}{a}=N(y)$ (recall that $p_a=dev_\tau(a,x)$ and $q_a=dev_\tau(a,y)$).
Fix a set of orthonormal affine coordinates $(y_0,\ldots, y_n)$ on
$\Kn{M}{n+1}$ such that $\mathbf N(x)=\der{0}{}$ and let
\[
   \psi_a:\{y_0=0\}\rightarrow \mathbb R
\]
such that $\widetilde S_1(\frac{\tau}{a})=\graph\psi_a$ and let $N_a$ the normal field of
$\mathcal D_{a^{-1}\tau}$. Since 
$N_a(a^{-1}p_a)=N(p_a)=\mathbf N(x)$ by using lemma \ref{technical} we have
that $\{a^{-1}p_a\}_{a>1}$ is a bounded set. Suppose that $a^{-1}p_a\rightarrow
p$ we have to show that $p=\mathbf N(x)$. By corollary \ref{CT surface varies
  with continuity} we have that $p\in\Kn{H}{n}$ and
$N_0(p)=\lim_{a\rightarrow+\infty} N_a(a^{-1}p_a)=\mathbf N(x)$. Since the
normal field on $\Kn{H}{n}$ is the identity we get that $p=\mathbf N(x)$.\par
Now let $c:[0,1]\rightarrow\Kn{H}{n}$ be a geodesic path between $\mathbf
N(x)$ and $\mathbf N(y)$
\[
    c(t)=(\psi_0(u(t)),u(t))
\]
where $\psi_0:\{y_0=0\}\rightarrow\mathbb R$ is the fuction such that
$\Kn{H}{n}=\graph\psi_0$.
Let $c_a(t)=(\psi_a(u(t)),u(t))$ be the corrisponding path on the surface
$\widetilde S_1(\frac{\tau}{a})$ and let $p'_a=c_a(0)$ and $q'_a=c_a(1)$. Since
$\psi_a\rightarrow\psi_0$ in $\mathrm C^1$-topology we have that
\[
\int_0^1(\E{\dot c_a(t)}{\dot c_a(t)})^{1/2}\mathrm d t\rightarrow\int_0^1(\E{\dot
  c(t)}{\dot c(t)})^{1/2}\mathrm d t=d_{\mathbb H}(\mathbf N(x),\mathbf N(y)).
\]
Let $a^{-1}p_a=(\psi_a(v_a),v_a)$ and $a^{-1}q_a=(\psi_a(w_a),w_a)$: it is easy to see that
\[
\frac{\delta_a(x,y)}{a}\leq ||v_a-u(0)||+||w_a-u(1)||+\int_0^1(\E{\dot c_a(t)}{\dot c_a(t)})^{1/2}
\]
Since $v_a\rightarrow u(0)$ and $w_a\rightarrow u(1)$  by passing to the
limit we get $d(\mathbf N(x),\mathbf N(y))\leq d_{\mathbb H}(\mathbf
N(x),\mathbf N(y))$.
\cvd

We want to show that the action of
$\Gamma_\tau$ on $(\widetilde S_a, a^{-1}d_a)$ converges in the Gromov sense to the action of
$\Gamma$ on  $\Kn{H}{n}$ for $a\rightarrow +\infty$.
For a complete definition of converegence in the Gromov sense of a sequence of
isometric actions on metric spaces see for istance \cite{Paulin}. However we
need only the following statement which is an immediate corollary of the
definition.\\

\emph{
Suppose that $\Gamma$ acts by isometries on a sequence of metric spaces $(X_i,d_i)$ and on a
metric space $(X_\infty, d_\infty)$. Suppose that 
there exists a sequence of  $\Gamma$-equivariant maps $\pi_i:X_i\rightarrow
X_\infty$ which verifies the following property: 
for every $K_\infty$ compact subset of $X_\infty$ and $\eps>0$ for $i>>0$ 
there exists a compact $K_i$ such that
$\pi_i(K_i)=K_\infty$ and $|d_\infty(\pi_i(x),\pi_i(y))-d_i(x,y)|<\eps$ for
all $x,y\in K_i$. Then the action of $\Gamma$ on $X_i$ converge in the sense
of Gromov to the action of $\Gamma$ on $X_\infty$.
}\\

\begin{corollario}\label{Gromov convergence in the future}
The action of $\Gamma_\tau$ on the rescaled surface $(\widetilde S_a, a^{-1}d_a)$
converges in the Gromov sense to the action of $\Gamma$ on $\Kn{H}{n}$.
\end{corollario}
\Dim
We want to see that the maps $N:\widetilde S_a\rightarrow\Kn{H}{n}$ satisfy the condition
above.
Fix a compact $K\subset\Kn{H}{n}$ and let $H=\mathbf N^{-1}(K)\subset\widetilde
M$. Since $\mathbf N$ is a proper map we have that $H$ is compact. Let
$K_a= dev_\tau(a,H)\subset \widetilde S_a$. By  proposition \ref{convergence in the future} 
we have that for all $\eps>0$ there exists $a_0$ such that for all $a>a_0$
\[
   \left|\frac{d_a(x,y)}{a}-d_{\Kn{H}{}}(N(x),N(y))\right|\leq\eps\qquad\textrm{for all }x,y\in K_a
\]
\cvd
Now we are interested in the asymptotic behaviour of the metrics $\delta_a$
when $a\rightarrow 0$. 
The results  are very similar to those that we have proved in the
previous case.
However in this case we shall see that there exist some
technical problems in the proofs.
In particular since $\partial\dom$ is an achronal set it is
defined a notion of lenght of a curve: however the lenght of a curve can be
$0$. By taking the infimum of the lenghts of the curves with fixed endpoints
we get a \emph{pseudo-distance} on $\partial\dom$. The first problem arises
when one tries to prove that this pseudo-distance restricted to $\Sigma$ is in
fact a distance. It seems to be more convenient to change viewpoint.
First we show that the distances
$\delta_a$ converges to a pseudo-distance $\delta_0$ on $\widetilde M$ such that
\[
     \delta_0(x,y)=0\quad\Leftrightarrow \quad \mathbf r(x)=\mathbf r(y).
\]
By this implication it follows that there exists a distance $d_\Sigma$ on
$\Sigma$ such that $\delta_0(x,y)=d_\Sigma(\mathbf r(x), \mathbf r(y))$.
After we shall prove that this distance coincides with the natural path-distance on
$\Sigma$.
\begin{prop}
There exsists a pseudo-distance $\delta_0$ on $\widetilde M$ such that
$\delta_a\rightarrow\delta_0$ in the compact-open topology of $\widetilde
M\times\widetilde M$. Moreover $\delta_0(x,y)=0$ if and only if $\mathbf r(x)=\mathbf r(y)$.
\end{prop}
\Dim
By using  lemma \ref{distance on S_a} one easily argues the first statement
as in the proof of proposition \ref{convergence in the future}.\par
The proof of the second statement is more difficult. We need the following
technical lemma.
\begin{lemma}\label{geodesics on S_a}
Let $\varphi:\mathbb R^n\rightarrow\mathbb R$ be a $\mathrm C^1$ convex
function such that $||\nabla\varphi(x)||<1$ for all $x\in\mathbb R^n$. Let 
$S=\{(x_0,\ldots, x_n)\in\Kn{M}{n+1}|x_0=\varphi(x_1,\ldots, x_n)\}$ be the
  corresponding spacelike surface in the Minkowski space and suppose $S$
  complete. Let $0$ be a minimum point of $\varphi$, then for every
  $y\in\Kn{R}{n}$ there exists a distance-minimizing geodesic arc 
  $c(t)=(\varphi (x(t)), x(t))$ 
  with starting point equal to $(\varphi(0),0)$ and ending point equal to
  $(\varphi(y), y)$ such that
   the functions $t\mapsto ||x(t)||$ and $t\mapsto\varphi(x(t))$ are 
  increasing ( $||\cdot||$
  is the euclidean norm of $\Kn{R}{n}$).
\end{lemma}
\emph{Proof of lemma:} 
First suppose that $\varphi$ is $\mathrm C^\infty$. By imposing that $c$ is a
geodesic we deduce that the path $x(t)$ satisfies the following equation
\[
   \ddot x(t)=\frac{\dot x\mathrm H\varphi(x)\dot
   x}{\left(1-||\nabla\varphi(x)||^2\right)^{3/2}}\nabla\varphi(x)
\]  
where $\mathrm H\varphi(x)$ is the Hessian matrix of $\varphi$ in $x$.\par
Let $f(t)=||x(t)||^2$, we have that 
\begin{eqnarray*}
\dot f (t)=2 x(t)\cdot\dot x(t);\\
\ddot f(t)=2 \Big(\dot x(t)\cdot \dot
x(t)+x(t)\cdot\ddot x(t)\Big).
\end{eqnarray*}
Now we have that $x(0)=0$ so that $\dot f(0)=0$. Hence it is sufficient to
prove that $\ddot f(t)\geq 0$ for $t\geq 0$. By looking at the last expression
it follows that it is sufficient to show that $x(t)\cdot\ddot x(t)\geq
0$. Since $\varphi$ is convex function we have that $\dot x\mathrm
H\varphi(x)\dot x\geq 0$ so that $x(t)\cdot\ddot x(t)\geq 0$ if and only if
$x\cdot\nabla\varphi(x)\geq 0$. 
On the other hand by using that $\varphi(0)$ is the minimum of
$\varphi$ and by imposing the convexity on the rays starting from $0$ one
easily deduces that this inequality holds for all $x\in\Kn{R}{n}$.
An analogous calculation shows that $t\mapsto \varphi(x(t))$ is increasing.\\

Now suppose that $\varphi$ is only $\mathrm C^1$. Let $\{\rho_\eps\}$ be a family
of $\mathrm C^\infty$ positive functions on $\Kn{R}{n}$ such that:\par
1. $\mathrm{supp}\rho_\eps=\{x\in\Kn{R}{n}|\,||x||\leq\eps\}$;\par
2. $\int_{\Kn{R}{n}}\rho_\eps=1$.\\
Let $\varphi_\eps$ the convolution $\varphi*\rho_\eps$
\[
    \varphi_\eps(x)=\int_{\Kn{R}{n}}\varphi(x-y)\rho_\eps(y)\mathrm dy.
\]
We know that $\varphi_\eps$ is $\mathrm C^\infty$ and 
$\varphi_\eps\rightarrow\varphi$ in $\mathrm C^1$-topology.
Moreover it is easy to see that $\varphi_\eps$ is a convex function so that
$S_\eps:=\graph\varphi_\eps$ is a $\mathrm C^\infty$ future convex spacelike surface.\par 
Fix $y\in\Kn{R}{n}$. By using the completeness of $S$ we have that
 for $\eps<<1$ there exists a path
\[
  x_\eps:[0,L_\eps]\rightarrow\Kn{R}{n}
\]
such that 
\begin{enumerate}
\item
$c_\eps(t)=(\varphi_\eps(x_\eps(t)),x_\eps(t))$ is a parametrization
of a distance-minimizing geodesic arc 
of the surface $S_\eps$;
\item
$x_\eps(0)=x$ and $x_\eps(L_\eps)=y$;
\item
$||\dot x_\eps(t)||=1$ and $L_\eps$ is bounded.
\end{enumerate}
Thus $x_\eps$ converges to a Lipschitz arc $x(t)$ and it is easy to see that
 the path $t\mapsto(\varphi(x(t)),x(t))$ is a distance-minimizing
geodesic between $(\varphi(0),0)$ and $(\varphi(y),y)$.\par
Let $p_\eps(t)$ be the orthogonal projection of $c_\eps(t)$ onto $T_{(\varphi_\eps(x),x)}S_\eps$ 
\[
 p_\eps(t)=
 c_\eps(t)+\frac{\E{c_\eps(t)}{(1,\nabla\varphi_\eps(x))}}{1-||\nabla\varphi_\eps(x)||^2}
   (1,\nabla\varphi_\eps(x)).
\]
We know that $\E{p_\eps(t)}{p_\eps(t)}$ is an increasing function of $t$. On
the other hand since $\nabla\varphi_\eps(x)\rightarrow 0$ for $\eps\rightarrow
0$ we get that $p_\eps(t)\rightarrow x(t)$. Thus $||x(t)||$ is an increasing
function of $t$. 
\cvd
Let us go back to the proof of proposition. We have to show that for all
$x,y\in\widetilde M$ such that
$\delta_0(x,y)=0$ we have  $\mathbf r(x)=\mathbf r(y)$.
By contradiction suppose that there exist $x,y\in\widetilde M$ such that
$\delta_0(x,y)=0$ and $\mathbf r(x)\neq\mathbf r(y)$.
Fix a set of affine coordinates $(y_0,\ldots, y_n)$ in such way that
$\der{0}{}=\mathbf N(x)$ and $\mathbf r(x)=0$. 
Let $\widetilde S_a=\mathrm{graph}\varphi_a$ and $\partial\dom=\graph\varphi$.
Moreover let $p_a=dev(a,x)=(a,0)$  and 
$q_a=dev(a,y)=(\varphi(z_a),z_a)$. 
Now for all $a>0$ fix a distance-minimizing
geodesic path 
\[
 c_a(t)=(\varphi_a(x_a(t)),x_a(t))\qquad \textrm{for }
  t\in[0,L_a]
\]
between $p_a$ and $q_a$ such that $||x_a(t)||$ is increasing.
Since $z_a\rightarrow z_0$ for $a\rightarrow 0$ there exists a constant
$K$ such that
\[
   ||x_a(t)||\leq K\qquad\textrm{ for all }a\leq 1.
\]
First suppose that there exists a sequence $a_k$ such that $L_{a_k}$ is 
bounded. Then up to passing to a subsequence we have that $x_{a_k}$
converges to a $1$-Lipschitz path $x:[0,L]\rightarrow\Kn{R}{n}$ such that
$x(0)=0$ and $x(L)=z_0$. Let $\varphi_a(t):=\varphi_a(x(t))$ then by the
hypothesis on $\delta_0$ we have that
\[
\lim_{k\rightarrow+\infty}
\int_0^{L_{a_k}}\sqrt{1-\dot\varphi_{a_k}(t)^2}\mathrm dt=0
\]
Thus $|\dot\varphi_{a_k}(t)|\rightarrow 1$ for almost all $t\in[0,L]$.
By lemma \ref{geodesics on S_a}  we know that $\varphi_{a_k}(x_a(t))$ are increasing functions of $t$
so that $\dot\varphi_{a_k}(t)\rightarrow 1$. Thus we have that
$\varphi_{a_k}(t)-\varphi_{a_k}(s)\rightarrow t-s$. On the other hand we
have that
$\varphi_{a_k}(t)-\varphi_{a_k}(s)\rightarrow\varphi(x(t))-\varphi(x(s))$.
So we obtain $\varphi(x(t))=t$. Thus we have that the path
$t\mapsto(\varphi(x(t),x(t))$ is a null path contained in $\partial\dom$
between $\mathbf r(x)$ and $\mathbf r(y)$. But this is a contradiction (in
fact we know that no point in $\Sigma$ lies in the interior of any null ray
contained in $\partial\dom$).\par
Hence suppose that $L_a\rightarrow +\infty$. Then we have that there exists
a sequence $a_k\rightarrow 0$ and a Lipschitz path
\[
  x:[0,+\infty)\rightarrow\Kn{R}{n}
\]
such that $x_{a_k}\rightarrow x$ in the compact open topology of 
$\Kn{R}{n+1}$. Since $||x_{a_k}(t)||\leq K$ we have that $||x(t)||\leq
K$. On the other hand the same argument used above shows that
$\varphi(x(t))=t$. Since $\varphi$ is $1$-Lipschitz we get $||x(t)||\geq
t$ and this gives a contradiction.
\cvd
From this proposition it follows that there exists a distance $d$ on $\Sigma$
such that
\[
   d(\mathbf r(x), \mathbf r(y))=\lim_{a\rightarrow 0}
   \delta_a(x,y)\qquad\textrm{for all }x,y\in\widetilde M.
\]
We have to see that $d$ coincides with the natural path-distance $d_\Sigma$.
Fix $r,s\in\Sigma$ and let $C(r,s)$ be the set of Lipschitzian path (with
respect the euclidean distance on $\Kn{M}{n+1}$)  in $\partial\dom$ between
$r$ and $s$  then $d_\Sigma(r,s)$ is defined by the rule
\[
   d_\Sigma(r,s):=\inf_{c\in C(r,s)}\int\sqrt{\E{\dot c(s)}{\dot c(s)}}\mathrm d s.
\]
\begin{prop}
For all $r,s\in\Sigma$ we have $d_\Sigma(r,s)=d(r,s)$.
\end{prop}
\Dim
It is easy to see that if $c:[0,1]\rightarrow \widetilde S_a$ is a rectifiable path then
$r\circ c$ is a rectifiable path with lenght lesser than the lenght of $c$. It
follows that $d_\Sigma(\mathbf r(x), \mathbf r(y))\leq \delta_a(x,y)$. Thus 
$d_\Sigma(r,s)\leq d(r,s)$.\par
Let us show the other inequality. Let $r,s\in\Sigma$ and $x,y\in\widetilde M$ such
that $\mathbf r(x)=r$ and $\mathbf r(y)=s$. Moreover let $p_a=dev(a,x)$ and
$q_a=dev(a,y)$.
Fix a set of affine orthonormal coordinates $(y_0,\ldots,y_n)$ and let
$\varphi:\Kn{R}{n}\rightarrow\mathbb R$
(resp. $\varphi_a:\Kn{R}{n}\rightarrow\mathbb R$) such that
$\partial\dom=\graph\varphi$ (resp. $\widetilde S_a=\graph\varphi_a$). Now we have
$r=(\varphi(u),u)$, $s=(\varphi(v),v)$, $p_a=(\varphi_a(u_a),u_a)$ and
$q_a=(\varphi_a(v_a),v_a)$. 
Finally let $q'_a=(\varphi_a(u),u)$ and $p'_a=(\varphi_a(u),u)$.\par
Now let $E\subset\Kn{R}{n}$ the set of points where $\varphi$ is not
differentiable. There exists a sequence of $1$-Lipshitz path
$x_k:[0,L_k]\rightarrow\Kn{R}{n}$ between $u$ and $v$ such
that $x_k^{-1}(E)$ has null Lebesgue measure on $[0,L_k]$ and 
\[
  \lim_{k\rightarrow+\infty}\int_0^{L_k}\sqrt{1-(\nabla\varphi(x_k(t))\cdot\dot
  x_k(t))^2}\mathrm d t=d_\Sigma(r,s).
\]
Consider the path $c_a^k(t)=(\varphi_a(x_k(t)),x_k(t))$. It is a path in $\widetilde S_a$
between $q'_a$ and $p'_a$ so that
\[
d_a(p_a,q_a)\leq
d_a(p_a,p'_a)\,+\,\int_0^{L_k}\sqrt{1-(\nabla\varphi_a(x_k(t))\cdot\dot
  x_k(t))^2}\mathrm d t\,+\,d_a(q_a,q'_a)
\]
Notice $(\varphi_a(x),x)+\ort{(1,\nabla\varphi_a(x))}$ is a support plane for $\widetilde S_a$ at
$(\varphi_a(x),x)$. So that the sequence of plane
$(\varphi_a(x),x)+\ort{(1,\nabla\varphi_a(x))}$ converges to a support plane
for $\partial\dom$ in $(\varphi(x),x)$. Thus it is easy to see that
$\nabla\varphi_a(x)\rightarrow\nabla\varphi(x)$ for $a\rightarrow 0$ and for
all $x\in\Kn{R}{n}-E$.
It follows that
\[
\lim_{a\rightarrow 0}\int_0^{L_k}\sqrt{1-(\nabla\varphi_a(x_k(t))\cdot\dot
  x_k(t))^2}\mathrm d t=\int_0^{L_k}\sqrt{1-(\nabla\varphi(x_k(t))\cdot\dot
  x_k(t))^2}\mathrm d t.
\] 
Now we have that $d_a(p_a,p'_a)\leq||u-u_a||$ (resp. $d_a(q_a, q'_a)\leq
||v-v_a||)$) so that by passing to the limit for $a\rightarrow 0$ we get
\[
  d(r,s)\leq\int_0^{L_k}\sqrt{1-(\nabla\varphi(x_k(t))\cdot\dot
  x_k(t))^2}\mathrm d t\qquad\textrm{ for all }k.
\] Now by passing to the limit for
$k\rightarrow+\infty$ we get $d(r,s)\leq d_\Sigma(r,s)$.
\cvd
In order to show the Gromov convergence of $\widetilde S_a$ to $\Sigma$
notice that we cannot us the argument of corollary
\ref{Gromov convergence in the future}. In fact there exists compact sets of
$(\Sigma, d_\Sigma)$ such that for all $a>0$ there exists no compact in $\widetilde S_a$
which projects on them. For istance consider the case $n=2$ and let $\tau\in
Z^1(\Gamma,\Kn{R}{2+1})$ such that the lamination associated with $\dom$ is
simplicial.
In this case the singularity is a simplicial tree such that every vertex is
the endpoints of a numerable set of edges. Fix a vertex $r_0$ and consider a
numeration $(e_k)_{k\in\mathbb N}$ of the edges with a endpoint equal to $r_0$.
Let 
\[
   K=\bigcup_{k\in\mathbb N}\{r\in e_k|d_\Sigma(r,r_0)\leq C/k\}
\]
where $C$ is the minimum of the lenghts of the edges of $\Sigma$.
It is easy to see that $K$ is compact. By contradiction suppose that for some
$a>0$ there exists a compact  $K_a$ such that $r(K_a)=K$. Now let
$\mathcal F(r_0)=N(r^{-1}(r_0))$: it is a complementary region of the lamination
and $\mathcal F(r)$ is a component of the boundary of $\mathcal F(r_0)$ for all $r\in K$.
Moreover $\mathcal F(r)$ depends only on the edge which contains $r$. Let $F_k$ be the
leaf corripsonding to $e_k$. Now fix $p_0\in K_a$ such that $r(p_0)=r_0$ and
for all $k$ let $p_k\in K_a$ such that $r(p_k)\in e_k$. We have that
$d_a(p_k, p_0)\geq a d_{\Kn{H}{}}(F_k, N(p_0))$. On the other hand
$d_{\Kn{H}{}}(F_k, N(p_0))\rightarrow +\infty$ for $k\rightarrow+\infty$ and
this contradicts the compactness of $K_a$.\par
In what follows we shall prove the convergence of the spectra of the
$\Gamma_\tau$-action  on $\widetilde S_a$ to the spectrum of the $\Gamma_\tau$-action on
$\Sigma$.
Generally let $(X,d)$ be a metric space provided with an action of
$\Gamma$. 
For every $\gamma\in\Gamma$ we can define the \emph{traslation lenght} of $\gamma$ as
$\ell_X(\gamma)=\inf_{x\in X} d(x, \gamma x)$. Clearly $\ell_X(\gamma)$
depends only on the conjugation class of $\Gamma$ and so it is defined a
function $\ell_x:\mathcal C\rightarrow[0,+\infty)$ where $\mathcal C$ be the
set of conjugation classes of $\Gamma-\{1\}$. This function is called
\emph{the marked lenght spectrum} of the action.
For semplicity we denote by $\ell_a$ ($a>0$) the marked lenght spectrum of the
$\Gamma_\tau$-action on the CT level surface $\widetilde S_a$, by $\ell_0$ the marked
lenght spectrum of the $\Gamma_\tau$-action on $\Sigma$ and by
$\ell_{\Kn{H}{}}$ the spectrum of the action on $\Kn{H}{n}$ .
\begin{corollario}   
With the above notation we have that for all $\gamma\in\Gamma$:
\begin{eqnarray*}
 \lim_{a\rightarrow +\infty}\ell_a(\gamma_\tau)/a=\ell_{\Kn{H}{}}(\gamma);\\
 \lim_{a\rightarrow 0}\ell_a(\gamma_\tau)=\ell_0(\gamma_\tau).
\end{eqnarray*}
\end{corollario}
\Dim
The first limit is consequence of the Gromov convergence. 
For the second limit notice that
$\ell_a(\gamma_\tau)\geq\ell_0(\gamma_\tau)$. On the other hand let
$x\in\widetilde M$ : then we have
\[
  \ell_a(\gamma)\leq\delta_a(x, \gamma x)
\]
so that we get $\limsup_{a\rightarrow 0}\ell_a(\gamma)\leq d_\Sigma(\mathbf
r(x),\gamma_\tau\mathbf r(x))$. Thus $\limsup_{a\rightarrow
  0}\ell_a(\gamma_\tau)\leq\ell_0(\gamma_\tau)$.
\cvd
\section{Measured Geodesic Stratification}\label{measure-section}
In section \ref{CT-section} we have associated  a geodesic stratification of $\Kn{H}{n}$ 
 with every future complete regular domain with surjective normal field. 
In this section we define the \emph{transverse measure} on a stratification. 
We have seen that in dimension $n=2$ the geodesic stratifications are in fact
 geodesic laminations.
 We see that for $n=2$   transverse measures on  geodesic stratifications are
 equivalent to  transverse measures on the corresponding geodesic laminations (in the classical
 sense).
Since
the behaviour of a stratification is quite more complicated than the behaviour
of a lamination, the general definition of transverse measure on a stratification is more involved. \par
We see that every measured geodesic stratification gives a future complete regular domain.
In his work Mess exposed a technique to associate a future complete regular
 domain of $\Kn{M}{2+1}$ to a measured geodesic lamination. This
construction is a generalization of that technique to any dimension,
$n\geq 2$.\\

Fix a \emph{complete weakly continuous  geodesic stratification} $\mathcal C$.
For $p\in\Kn{H}{n}$ let us indicate with $C(p)$ the piece in $\mathcal C$ which contains $p$ and
has minimum dimension.\par
The first notion that we need is the \emph{ transverse measure} on a piece-wise geodesic path.
Let $c:[0,1]\rightarrow\Kn{H}{n}$ be a \emph{piece-wise geodesic  path}. 
A \emph{transverse measure} on it is a \emph{$\Kn{R}{n+1}$-valued measure} $\mu_c$ on $[0,1]$ such that
\begin{enumerate}
\item
  there exists a finite positive measure $|\mu_c|$ such that  $\mu_c$ is
  $|\mu_c|$-absolutely continuous 
  and $\supp|\mu_c|$ it the topological closure of the set $\{t\in(0,1)|\,\dot c(t)\notin T_{c(t)}C(c(t))\}$;
\item
  let $v_c=\frac{\mathrm{d}\mu_c}{\mathrm{d}|\mu_c|}$ be the $|\mu_c|$-density
  of $\mu_c$, then 
  \begin{equation}\label{density property}
  \begin{array}{ll}
  v_c(t)\in T_{c(t)}\Kn{H}{n}\cap\ort{T_{c(t)}C(c(t))},&\\
  \E{v_c(t)}{v_c(t)}=1, & \\
  \E{v_c(t)}{\dot c(t)}>0 & |\mu_c|-\textrm{ a.e.}
  \end{array}
   \end{equation}
\item
  the endpoints of $c$ are not atoms of the measure $|\mu_c|$.
\end{enumerate}
Let us point out an useful property of a transverse measure on a geodesic path.
\begin{lemma}\label{density separes endpoints}
Let $c:[0,1]\rightarrow\Kn{H}{n}$ be  a \emph{geodesic} path and $\mu_c$ be a
transverse measure on it. Then for $|\mu_c|$-almost all $t$  we have
\[
     \E{c(0)}{v_c(t)}<0\qquad\E{c(1)}{v_c(t)}>0.
\]
Thus $\ort{v_c(t)}$ separes $c(0)$ from $c(1)$. 
\end{lemma}
\Dim
Since $c$ is a geodesic path there exists $v\in T_{c(0)}\Kn{H}{n}$ such that
\[
      c(t)=\cosh(s(t))c(0)+\sinh(s(t))v
\]
with $s(t)$ an increasing function.
By (\ref{density property}) we have  $|\mu_c|$-almost everywhere
\begin{eqnarray*}
   0=\E{c(t)}{v_c(t)}=\cosh(s(t))\E{c(0)}{v_c(t)}+\sinh(s(t))\E{v}{v_c(t)};\\
   0<\E{\dot c(t)}{v_c(t)}=\dot
   s(t)\big(\ \sinh (s(t))\E{c(0)}{v_c(t)}\ +\ \cosh(s(t))\E{v}{v_c(t)}\ \big).
\end{eqnarray*}
Looking at this expressions we easily get that
$\E{c(0)}{v_c(t)}<0$ . An analogous calculation shows the other inequality.
\cvd
A first consequence of this lemma is that the measure $\mu_c$ determines the
positive measure $|\mu_c|$. 
\begin{corollario}
Let $c$ be a piece-wise geodesic path and $\mu_c$ be a transverse measure on
it. Suppose that $\lambda$ is a positive measure such that $\mu_c$ is
$\lambda$-absolutely continuous and the density  $u=\frac{\mathrm
  d\mu_c}{\mathrm d\lambda}$ verifies (\ref{density property}). Then
$\lambda=|\mu_c|$.
\end{corollario}
\Dim
First let us show  that $\lambda$ is $|\mu_c|$-absolutely
continuous. Let $E\subset [0,1]$ such that $|\mu_c|(E)=0$. We can suppose that
$E$ is contained in an interval $I=[t_0,t_1]$ such that $c|_I$ is a geodesic
path.  Since $\mu_c(E)=0$
we deduce that
\[
    0=\E{c(t_0)}{\int_E u(t)\mathrm d\lambda}=\int_E\E{c(t_0)}{u(t)}\mathrm
    d\lambda.
\]
The same argument of lemma \ref{density separes endpoints} shows that
$\E{c(t_0)}{u(t)}<0$ for $\lambda$-almost all points in $I$. Thus $\lambda(E)=0$.\par
Now let $a=\frac{\mathrm d\lambda}{\mathrm d|\mu_c|}$. We recall that
\[
    v_c(t)=\frac{\mathrm d\mu_c}{\mathrm d|\mu_c|}=\frac{\mathrm
      d\mu_c}{\mathrm d\lambda}\frac{\mathrm d\lambda}{\mathrm d|\mu_c|}=a(t)u(t)
    \qquad |\mu_c|-\textrm{a.e.}
\]
Since $\E{v_c(t)}{v_c(t)}=\E{u(t)}{u(t)}=1$ we deduce that $a(t)=1$.
\cvd

In order to define a transverse measure on a geodesic stratification we need
the following definition.
\begin{defi}
Let $\varphi_s:[0,1]\rightarrow\Kn{H}{n}$ be a homotopy between $\varphi_0$
and $\varphi_1$. We say that $\varphi$ is
$\mathcal C$-\emph{preserving} if $C(\varphi_s(t))=C(\varphi_0(t))$ for all
$(t,s)\in[0,1]\times[0,1]$ (recall that $C(x)$ is the piece of $\mathcal C$
which contains $x$ and has minimum dimension).
\end{defi} 
   
Now we can give the definition of trasverse measure on a geodesic
stratification.
\begin{defi}\label{transverse measure-definition}
Let $\mathcal C$ be a complete weakly continuous stratification and fix a
subset $Y\subset\Kn{H}{n}$ which is union of pieces of $\mathcal C$ such that
the Lebesgue measure of $Y$ is $0$. 
We mean by $(\mathcal C, Y)$-\emph{admissible path} (or simply
\emph{admissible path}) any
piece-wise geodesic path $c:[0,1]\rightarrow\Kn{H}{n}$ such that every maximal
geodesic subsegment has no endpoint in $Y$.\par
A \textbf{transverse measure} on $(\mathcal C, Y)$  is the assignment of a
transverse measure $\mu_c$ to every admissible path
$c:[0,1]\rightarrow\Kn{H}{n}$ such that
\begin{enumerate}
\item
  if  there exists a $\mathcal C$-preserving homotopy between two paths $c$ and
  $d$ then $\mu_c=\mu_d$;
\item
  for every admissible path  $c$ and every  parametrization
  $s:[0,1]\rightarrow[0,1]$ of an admissible
  subarc of $c$  we have
  that $\mu_{c\circ s}=s^*(\mu_c)$; 
\item 
 the atoms of $|\mu_c|$ are contained in $c^{-1}(Y)$, and for every $y\in Y$
 there exists an admissible path $c$ such that $|\mu_c|$ has some atoms on $c^{-1}(y)$;
\item
  $\mu_c(c)=0$ for every closed admissible path $c$;
\item
  for all sequences $(x_k)_{k\in\mathbb N}$ such that $x_k\in\Kn{H}{n}-Y$ and
  $x=\lim_{k\rightarrow+\infty}x_k\in\Kn{H}{n}-Y$ we have that
  $\mu_{c_k}(c_k)\rightarrow 0$ where $c_k$ is the admissible arc $[x_k, x]$. 
\end{enumerate}
A \textbf{measured geodesic stratification} is given by a weakly continuous
geodesic stratification, a subset $Y$ as above and a measure on $\mu$ on
$(\mathcal C,Y)$.\par
A measured geodesic stratification $(\mathcal C, Y, \mu)$ is
$\Gamma$-\emph{invariant} if $\mathcal C$ is $\Gamma$-invariant, $Y$ is
$\Gamma$-invariant and for all $c:[0,1]\rightarrow\Kn{H}{n}$ admissible path,
$E\subset[0,1]$ borelian set, and $\gamma\in\Gamma$ we have
\[
    \mu_{\gamma\circ c}(E)=\gamma(\mu_c(E)).
\]
\end{defi}
\begin{figure}[h]
\begin{center}
\input{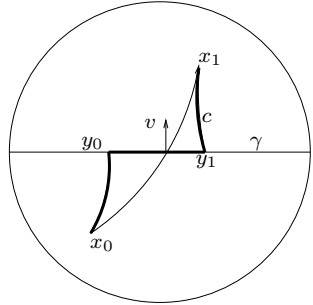}
\caption{{\small The figure shows a non-admissible arc.}}\label{controesempio-fig}
\end{center}
\end{figure}
\begin{oss}\emph{
The restriction to admissible subarcs is necessary for the foundation of
definition. For istance consider the stratification of $\Kn{H}{2}$ with one
geodesic $\gamma$. Fix a geodesic arc $c$: if $c$ intersects transversally
$\gamma$ in a point $t_0$ put $\mu_c=v\delta_{t_0}$ where $v$ is the normal to
$\gamma$ such that $\E{v}{\dot c(t_0)}>0$. If $c$ is contained in $\gamma$ or
does not intersect it put $\mu_c=0$. It is not possible to extend this
definition to all piece-wise geodesic paths: in fact consider the segment $c$ as
in fig.\ref{controesempio-fig}. Applying property 4. of definition of measure
to the closed admissible arc $c*[x_1,x_0]$ we get $\mu_c(c)=v$. On the other hand we
have that $\mu_c([x_0,y_0))=0$, $\mu_c([y_0,y_1])=0$ and $\mu_c((y_1,x_1])=0$
and this is a contradiction.}\par
\emph{
Consider a measured geodesic stratification $(\mathcal C,
Y,\mu)$, notice that condition 3. in definition
\ref{transverse measure-definition} impose a minimality property of $Y$.  
}\end{oss}
\begin{oss}\emph{
Consider the case  $n=2$. Let $\Gamma$ be a co-compact Fuchsian group. Let
$\mathcal C$ be a $\Gamma$-invariant geodesic stratification. The $1$-stratum
$L$ of $\mathcal C$ is a $\Gamma$ -invariant geodesic lamination of
$\Kn{H}{2}$. We know that $L=S\cup L_1$ where $S$ is a simplicial lamination
and $L_1$ is a lamination with no closed leaf (see \cite{Casson} for further
details about geodesic laminations).}\par
\emph{
We want to see that maximal measured geodesic stratification  $(\mathcal C, Y,\mu)$ are naturally
identified with the transverse measures (in the classical sense)  on the lamination $L=X_{(1)}$ (notice that
these concepts are quite different, in fact one is a $\Kn{R}{n+1}$ valued
measure and the other is a positive measure).}\par
\emph{
Fix a transverse measure $\mu$ of $(\mathcal C,Y)$: we want to see that there exists
a unique transverse measure $\lambda$ on $L$ such that $\lambda_c=|\mu_c|$ for
all admissible paths
$c$. First notice that $Y$ is union of geodesics of $L$ (in fact the interior
of $L$ is empty). It follows that every transverse path $d$ is composition of
paths $d_i$ such that  there exists a $\mathcal C$-preserving
homotopy between $d_i$ and a suitable parametrization $c_i$ of the admissible geodesic
segment $[d_i(0),d_i(1)]$. Thus one define $\lambda_d$ such that the
restriction on $d_i$ is $|\mu_{c_i}|$. By using property 1. and 2. of
definition \ref{transverse measure-definition} it is easy
to see that this definition does not depend on the choice of the
decomposition and in fact it is the only possible one. Finally since
$\mu_c$ is $\Gamma$-invariant we have that $\lambda$ is $\Gamma$-invariant
too.}\par
\emph{
By general facts about measured geodesic laminations (see \cite{Casson}) it
follows that for every admissible path $c$ the atoms of $\mu_c$ are exactly
$c^{-1}(S)$. It follows that $S=Y$.}\par
\emph{
Conversely let $\lambda$ be a $\Gamma$-invariant transverse measure on $L$.  Put $Y=S$. We want to construct a
$(\mathcal C,Y)$-measure on $\Kn{H}{n}$.
Fix an admissible path
$c$ and let $v_c:[0,1]\rightarrow\Kn{R}{n+1}$ the function so defined:
$v_c(t)=0$ if $c(t)\notin L$, otherwise $v_c(t)$ is the normal vector to the
leaf $C(c(t))$ such that $\E{v_c(t)}{\dot c(t)}>0$. Then we can define $\mu_c$
the $\Kn{R}{n+1}$-measure on $[0,1]$ which is $\lambda_c$ absolutely continuous
and has $\lambda_c$ density equal to $v_c$. It is easy to see that in this way
$\mu_c$ is a transverse measure on $c$. Furthermore by definiton of the assignment
$c\mapsto\mu_c$ verifies condition 1. and 2. of definition \ref{transverse measure-definition}. An
easy  analysis of the geometry of a lamination show that the condition 4. and
5. are satisfied. 
}\par
\emph{
It is easy to see that this correspondence gives an identification between
$\Gamma$-invariant transverse measure on $\mathcal C$ and $\Gamma$-invariant transverse measure on $L$.
}\par
\emph{
Notice that in dimension $n=2$ the condition 4. of definition
\ref{transverse measure-definition} is  assured by
the geometry of the stratification. Furthermore in this case the set $Y$ is
determined by the lamination (i.e. it does not depend on the measure). 
}\end{oss}
Before constructing a future complete regular domain with a given geodesic
stratification, let us point out an easy property of measured geodesic
stratifications.
\begin{lemma}\label{uniqueness of maximal piece}
Let $\mu$ be a transverse measure on $(\mathcal C,Y)$. Then for all
$x\in\Kn{H}{n}-Y$ we have that there exists a unique maximal piece of $\mathcal
C$ which contains $x$ (maximal with respect the inclusion).
\end{lemma}
\Dim
Suppose that there exist two pieces $C_1,C_2$ which contain $x$. We want to show
that there exists a piece $C$ which contains $C_1\cup C_2$.\par
Let $x_i$ be a point on $C_i$ such that $C(x_i)=C_i$. Notice that $x_i$ does
not lie
into $Y$ (in fact $Y$ is union of pieces so that if $y\in Y$ then $C(y)\subset
Y$). Consider the piece-wise geodesic arc $c =[x_1,x]\cup[x,x_2]\cup[x_2,x_1]$. It
is closed and admissible so that
\[
        \mu_c([x_2,x_1])=-\mu_c([x_1,x])-\mu_c([x,x_2])
\]
(notice that we are using that $x,x_1$ and $x_2$ are not atoms). Since $(x_1,x)$ and
$(x,x_2)$ are contained in $C$ one easily see that $\mu_c([x_1,x_2])=0$.
On the other hand  by lemma \ref{density separes endpoints} we have that
$\E{v_c(t)}{x_2-x_1}>0$ for $|\mu_c|$-almost all $t$. Since 
$\E{\mu_c([x_1,x_2])}{x_2-x_1}=\int_{[x_1,x_2]}\E{v_c(t)}{x_2-x_1}$ we have
that $|\mu_c|([x_1,x_2])=0$. Thus the segment $[x_1,x_2]$ is contained in a
piece $C$. It follows that $C\supset C(x_i)=C_i$.
\cvd
 
Now given a measured geodesic stratification $(\mathcal C, Y,\mu)$ we 
construct a regular domain  with stratification equal to $\mathcal C$.  
Fix a base point $x_0\in\Kn{H}{n}-Y$ and define for $x\notin Y$
\[
      \rho(x)=\mu_{c_x}(c_x) 
\]
where $c_x$ is a admissible path between $x_0$ and $x$. It is quite evident that this
definition does not depend on the choice of the path. Furthermore notice that
\[
     \rho(y) = \rho(x) +\mu_{c_{x,y}}(c_{x,y})
\]
where $c_{x,y}$ is the geodesic arc between $x$ and $y$.
By using property 5. it follows that the map
$\rho:\Kn{H}{n}-Y\rightarrow\Kn{M}{n+1}$ is continuous.
\par 
Let us denote $M(x)$ the maximal piece of $\mathcal C$ which contains $x$ for
$x\notin Y$. By lemma \ref{uniqueness of maximal piece} this piece is unique.
By using lemma \ref{density separes endpoints} it is easy to see that
\begin{equation}\label{dominio-formula1}  
\E{\rho(y)-\rho(x)}{y}\geq 0\qquad\E{\rho(y)-\rho(x)}{x}\leq 0.
\end{equation}
Furthermore arguing as in lemma \ref{uniqueness of maximal piece} one easily
see that the equality holds if and only if $M(x)=M(y)$. Thus $\rho(x)=\rho(y)$
if and only if $M(x)=M(y)$.\par
Let us define a convex set
\[
    \Omega=\bigcap_{x\in\Kn{H}{n}-Y}\fut(\rho(x)+\ort x).
\]
\begin{teo}\label{stratification gives  domain}
We have that $\Omega$ is a future complete regular domain. Furthermore
$\rho(x)\in\Sigma$ for all $x\in\Kn{H}{n}-Y$ and $\Omega$ is the convex hull
of the points $\{\rho(x)|x\in\Kn{H}{n}-Y\}$.
\end{teo}
\Dim
First let us show that $\rho(x)\in\partial\Omega$. Clearly $\rho(x)\notin\Omega$. Now
let $v\in\fut(0)$: we have to show that $\rho(x)+v\in\Omega$ i.e. 
$\E{\rho(x)+v-\rho(y)}{y}< 0$ for all $y\in\Kn{H}{n}$.
By inequalities (\ref{dominio-formula1}) $\E{\rho(x)-\rho(y)}{y}\leq 0$ 
and since $v$ is future directed $\E{v}{y}<0$ so that  $\fut(\rho(x))\subset\Omega$.\par
Now notice that $\rho(x)=\rho(y)$ for every $y\in M(x)$ so that the plane
$\rho(x)+\ort{y}$ is a support plane for $\Omega$ for all $y\in M(x)$.
Let $v$ be a null direction such that $[v]$ is in the boundary of $M(x)$. By
taking a sequence $(y_k)\in M(x)$ such that $y_k\rightarrow [v]$ it is easy to
see that the plane $\rho(x)+\ort{v}$ is a support plane for $\Omega$.
So that we have shown
\[
   \Omega\subset\bigcap_{\begin{array}{c} 
       {\scriptstyle x\in\Kn{H}{n}-Y\textrm{ and}}\\
     {\scriptstyle [v]\in M(x)\cap\partial\Kn{H}{n}}\end{array}}
     \fut(\rho(x)+\ort{v}).
\]
If we prove the other inclusion we obtain that $\Omega$ is a future complete
 regular domain.
 Fix $p\in\Kn{M}{n+1}$ and suppose that $\E{p-\rho(x)}{v}<0$ for
all $x\in\Kn{H}{n}$ and $[v]\in M(x)\cap\partial\Kn{H}{n}$. We have to show
that $p\in\Omega$. Notice that every $x\in\Kn{H}{n}-Y$ is a convex combination of
a collection $v_1,\ldots,v_n$ such that $[v_i]\in M(x)\cap\partial\Kn{H}{n}$.
It follows that $\E{p-\rho(x)}{x}<0$ and so that $p\in\Omega$.
Since $\rho(x)+\ort{x}$ is a spacelike support plane, we argue that
$\rho(x)\in\Sigma$ and  $\mathcal F(\rho(x))\supset M(x)$.\par 
Now let us prove
 that $\Omega$ is the convex hull of the points $\rho(x)$. Notice that it is
 sufficient to show that the future $\fut(\widetilde S_a)$ of the CT level surface $\widetilde S_a$
 is the convex hull of the set  $\mathcal S_a=\{ax+\rho(x)|x\in\Kn{H}{n}-Y\}$. 
Since $\rho(x)+\ort{x}$ is a support plane of $\Omega$ we have that
 $ax+\rho(x)\in \widetilde S_a$ so that the future of $\widetilde S_a$ contains the convex hull of
 $\mathcal S_a$. On the other hand it is easy to see that the spacelike
 support planes of the convex hull of $\mathcal S_a$ are support planes of
 $\fut(\widetilde S_a)$. Thus in order to prove the statement it is sufficient to show
 that the convex hull of $\mathcal S_a$ has not timelike support plane.
By contradiction suppose that there exists a vector $v$ such that $\E{v}{v}=1$
 and $\E{ax+\rho(x)}{v}<C$. Up to translating $\Omega$ we
 can suppose that the base point $x_0$ is orthogonal to $v$. Consider the
 geodesic $\gamma$ such that $\gamma(0)=x_0$ and $\dot\gamma(0)=v$. We can
 suppose that there exists a sequence $t_k\rightarrow+\infty$ such that
 $\gamma(t_k)\notin Y$. Now let $x_k=\gamma(t_k)$, by using that
 $\E{\rho(x_k)}{x_k}\geq 0$ and $\E{\rho(x_k)}{x_0}\leq 0$ we deduce that
 $\E{\rho(x_k)}{v}\geq 0$. Since $\E{x_k}{v}\rightarrow+\infty$ we have that
$\E{ax_k+\rho(x_k)}{v}\rightarrow+\infty$ and this gives a contradiction.\par
\cvd
Now we want to show that the stratification associated with $\Omega$ coincides
with $\mathcal C$ at least on $\Kn{H}{n}-Y$.
First we give a technical result.
\begin{lemma}\label{counterimage of N}
Let $\Omega$ be the regular domain constructed in theorem \ref{stratification
  gives domain} and suppose that the map $N:\widetilde
  S_1\rightarrow\Kn{H}{n}$ is proper. 
Then for all $x\in\Kn{H}{n}$ we have
  that $N^{-1}(x)\cap\widetilde S_1$ is the convex hull  of the limit points
  of the sequences $(x_k+\rho(x_k))_{k\in\mathbb N}$ such that 
$x_k\in\Kn{H}{n}-Y$ and $Nx_k\rightarrow x$.
\end{lemma}
\Dim
Let $p\in\widetilde S_1$ and suppose that $N(p)=x$. We have to show that for
all $v\in\ort{N(p)}$ there exists a sequence $x_k\in\Kn{H}{n}-Y$ , 
such that $x_k+\rho(x_k)\rightarrow p_\infty$ with $N(p_\infty)=x$ and 
$\E{p_\infty-p}{v}\geq 0$.
We know that there exists a sequence $x_k\in \Kn{H}{n}-Y$ such that
$x_k=\cosh d_k x+\sinh d_k v_k$ such that $d_k\rightarrow 0$ and
$v_k\rightarrow v$. Let $p_k=x_k+\rho(x_k)$. Since $N$ is a proper map, up to
passing to a subsequence we get that $p_k\rightarrow p_\infty$ such that
$N(p_\infty)=x$.
Since $\E{p_k-p}{x}\leq 0$ and $\E{p_k-p}{x_k}\geq 0$ we have that
$\E{p_k-p}{v_k}\geq 0$ and passing to the limit we get $\E{p_\infty-p}{v}\geq
0$.
\cvd
\begin{prop}\label{stratification11}
Let $\Omega$ be the regular domain associated with the measured stratification
$(\mathcal C, Y, \mu)$, and suppose that the restriction of the normal field
$N|_{\widetilde S_1}$ is a proper map. Then for all $x\in\Kn{H}{n}-Y$ we have that
$rN^{-1}(x)=\{\rho(x)\}$.
Moreover $\mathcal F(\rho(x))$ is the maximal
piece $M(x)$. (We recall that $\mathcal F(r_0)=N(r^{-1}(r_0))$ where
$r_0\in\Sigma$ and $r$ is the retraction onto the singularity
$\Sigma$). 
\end{prop}
\Dim
Let $x\in\Kn{H}{n}-Y$.
Since $\rho:\Kn{H}{n}-Y\rightarrow\Kn{M}{n+1}$ is a continuous map  lemma
\ref{counterimage of N} implies that $N^{-1}(x)\cap\widetilde
S_1=\{x+\rho(x)\}$.  Thus $rN^{-1}(x)=\{\rho(x)\}$.
Now let us prove that $\mathcal F(\rho(x))=M(x)$. Clearly $M(y)\subset\mathcal
F(\rho(x))$, thus we have to prove the other inclusion.  
Let $y\in \mathcal F(\rho(x))$, first suppose that $y\notin Y$. 
We have $\E{y}{\rho(y)-\rho(x)}\leq 0$. On the other hand by the (\ref{dominio-formula1}) we know that the other
inequality holds so that $\E{y}{\rho(y)-\rho(x)}=0$, but then
$\rho(x)=\rho(y)$ and so $M(x)=M(y)$.\par
Suppose now that $y\in Y\cap\mathcal F(\rho(x))$. Let us prove that $y$ lies
in the boundary  $b\mathcal F(\rho(x))$ (see section \ref{CT-section} for the definition of the
boundary $bK$ of a convex $K$). Since $y\in Y$ we have that $N^{-1}(y)\cap
\widetilde S_1$ is not a point. So that there exists a spacelike vector $v$
orthogonal to $y$ such that the segment $[\rho(x),\rho(x)+\eps v]\subset\Sigma$.
We have that $\ort v$ is a support plane of $\mathcal F(\rho(x))$ and contains
$y$. So that if $y\notin b\mathcal F(\rho(x))$ we get $x\in\mathcal
F(\rho(x)+\eps v)$ i.e. $\rho(x)+\eps v\in rN^{-1}(x)$. Since $x\notin Y$ we
get a contradiction. Thus $\mathcal F(\rho(x))-b\mathcal
F(\rho(x))\subset\Kn{H}{n}-Y$ so that $\mathcal F(\rho(x))-b\mathcal
F(\rho(x))\subset M(x)$. It follows that $\mathcal F(\rho(x))\subset M(x)$.
\cvd
\begin{oss}\emph{
Notice that the stratification induced by the domain $\Omega$ coincides with
$\mathcal C$ on $\Kn{H}{n}-Y$, so that
$Y$ is the union of pieces of the stratification associated with
$\Omega$. Moreover by property 3. in definition \ref{transverse
  measure-definition} it turns out that $Y=\{y\in\Kn{H}{n}|\#N^{-1}(y)>1\}$.
}\end{oss}
\begin{corollario}
Let $(\mathcal C, Y,\mu)$ be a $\Gamma$-invariant measured geodesic
stratification of $\Kn{H}{n}$ (where $\Gamma$ is a co-compact free-torsion
discrete subgroup of $\SOO^+(n,1)$). Fix a
base point $x_0\notin Y$ and let $\tau_\gamma=\rho(\gamma(x_0))$
then we have that $\tau\in Z^1(\Gamma,\Kn{R}{n+1})$.
Let $\Omega$ be a domain associated with $(\mathcal C, Y,\mu)$ . 
Then  $\Omega=\dom$ and $\mathcal F(\rho(x))=M(x)$ for all $x\notin Y$.
\end{corollario}
\Dim
Since $\mu$ is $\Gamma$-invariant we have
\[
       \rho(\gamma(x))=\gamma \rho(x)+\rho(\gamma(x_0))
\]
Thus $\tau_{\alpha\beta}=\alpha\tau_\beta+\tau_\alpha$, so that $\tau$ is a
cocycle. The same equality shows that $\Omega$ is a $\Gamma_\tau$-invariant
regular domain. Thus   theorem \ref{uniqueness} implies that it is equal to
$\dom$. Since in this case we know that $N:\widetilde S_1\rightarrow\Kn{H}{n}$
is a proper map proposition \ref{stratification11} implies that $\mathcal F(\rho(x))=M(x)$ for all $x\notin Y$.
\cvd
\section{Simplicial Stratifications}\label{simplicial-section}
In this section we restrict to the case $n+1=4$.
In particular we study the future complete regular domain  of $\Kn{M}{3+1}$ with simplicial geodesic
stratification.  
\begin{defi}
We say that a geodesic stratification $\mathcal C$ is \emph{simplicial} if for every
$p\in\Kn{H}{n}$ there exists a neighborhood $U$ which intersects only a finite
number of pieces of $\mathcal C$.
\end{defi}
\begin{figure}
\begin{center}
\input{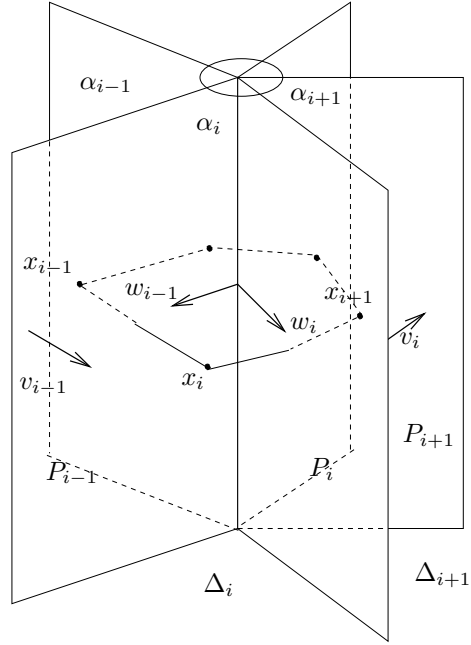}
\caption{{\small A neighborhood of a point on a $1$- piece of a simplicial
  stratification of $\Kn{H}{3}$}}\label{simplicialstrat-fig}
\end{center}
\end{figure}
We shall see that the correspondence between measured geodesic stratifications
and regular domains induce an identification between measured simplicial
stratifications and regular domain \textbf{ with simplicial
singularity}. Finally we shall recover the duality between stratification
and singularity.\par
Notice that in dimension $n=2$  simplicial stratifications correspond to  simplicial
laminations.
Moreover in all the dimensions a simplicial stratification has closed
strata: it is in fact a tasselation of $\Kn{H}{n}$ by locally finite ideal
convex polyhedra.
In fig.\ref{simplicialstrat-fig}  we show the local behaviour of a simplicial stratification of $\Kn{H}{3}$.\\

First we describe the measures on a simplicial stratification.
Let $(\mathcal C, Y,\mu)$ be a measured geodesic lamination with simplicial
support $\mathcal C$. Let $X$ be the $2$-stratum of $\mathcal C$. We want to
show that $X=Y$. Since $Y$ has not interior we have $Y\subset X$. On the other
hand let $c$ be a geodesic path with no endpoint in $X$ (a such path is
admissible). Notice that $\supp|\mu_c|$ is  $c^{-1}(X)$, but this set is
finite so that the measure $|\mu_c|$ has an atom on every point of
$c^{-1}(X)$. Thus $X\subset Y$.\par
Fix a $2$-piece $P$ and let $\Delta_1$ and $\Delta_2$ be 
the $3$-pieces which incide on $P$. Let $c$ an admissible  geodesic path
which start in $\Delta_1$ and termines in $\Delta_2$. Clearly $\mu_c=a
v_{1,2}\delta_{c^{-1}(P)}$ where $a$ is positive constant $v_{1,2}$ is the
normal vector to $P$ which points toward $\Delta_2$ 
and $\delta_x$ is the Dirac measure centered in $x$. By using  property 1. and 2. of definition
\ref{transverse measure-definition} one easily see that the constant $a$ does not depends on the
path.  By imposing that $\mu_c(c)+\mu_{c^{-1}}(c^{-1})=0$ one deduce that the
measure of a geodesic path $c'$ which starts in $\Delta_2$ and termines in
$\Delta_1$ is $\mu_\delta=av_{2,1}\delta_{c^{-1}(P)}$. It follows that the
constant $a$ depends only on the piece $P$: we call it the \emph{weight} of $P$
and we denote it by $a_\mu(P)$.\par
We want to show that the set of the weight $\{a_\mu(P)|P\textrm{ is a }2\textrm{
  piece}\}$ satisfies a certain set of equations and determines the measure $\mu$.\par
Fix a geodesic $l$ in $\mathcal C$ and let $P_1,\ldots P_k$  and
$\Delta_1,\ldots,\Delta_k$ be respectively the $2$-pieces and the $3$-pieces
which incide on $l$. We can suppose that the numeration is such that $\Delta_i$
incides on $P_{i-1}$ and $P_i$ (the index $i-1$ and $i$ are considered $\mod
k$, see fig.\ref{simplicialstrat-fig}).
Fix $x_i$ in $\mathrm{int}(\Delta_i)$ and consider the admissible closed path 
$c= [x_1, x_2]\cup[x_2, x_3]\cup\ldots\cup[x_{k-1}, x_k]\cup[x_k, x_1]$. By
imposing that $\mu_c(c)=0$ we deduce 
\begin{equation}\label{condition}
    \sum_{i=1}^k a_\mu(P_i)v_i=0
\end{equation}
where $v_i$ is  the normal vector to $P_i$ which points outward
$\Delta_i$. Notice that $v_i$ lies in the linear subspace of $\mathbb M^{3+1}$
which is orthogonal to the space generated by $l$ (which we denote $\ort
l$). Fix a point $x\in l$, notice that $\ort l$ is identified with the
subspace of $T_x\Kn{H}{3}$ which is orthogonal to $l$. 
By performing a $\frac{\pi}{2}$-rotation on $\ort l$ one see that equation
(\ref{condition}) is equivalent to the equation
\begin{equation}\label{condition-bis}
   p_l(a_\mu)=\sum_{l\subset P} a_\mu(P)w(P)=0
\end{equation}
where $w$ is the unitary vector of $\ort l$ which is tangent to $P$ and
points inward (see fig.\ref{simplicialstrat-fig}).
\begin{defi}
A family of positive constants $a=\{a(P)\}$ paramatrized by the set of the
$2$-pieces of $\mathcal C$ is called a family of \emph{weights} for the
stratification if the equation $p_l(a)=0$ is satisfied for every $1$-piece $l$.
\end{defi}
We have shown that there is a family of weights associated with every
transverse measure $\mu$  on $\mathcal C$. 
Now we want to prove that this
corrispondence is bijective.
\begin{prop}
For every family of weights $\{a(P)\}$ there exists a unique transverse
measure $\mu$ such that $a(P)=a_\mu(P)$.
\end{prop}
\Dim
First let us prove the uniqueness. Let $\mu$ and $\nu$ measure such that  
$a_\mu(P)=a_\nu(P)$ for all $2$-pieces $P$. If $c$ is an admissible arc which
intersects only $2$-pieces it follows that $\mu(c)=\nu(c)$. Suppose now that
$c\cap X$ is a point $p$ which lies on the $1$-piece $l$. Consider an arc $c'$ which has the
same endpoints of $c$ and does not intersect the $1$-stratum. We have 
$\mu_c(c)=\mu_{c'}(c')=\nu_{c'}(c')=\nu_c(c)$. Since
$\supp\mu_c=\supp\nu_c=c^{-1}(p)$ we have $\mu_c=\nu_c$.
\begin{figure}[!h]
\begin{center}
\input{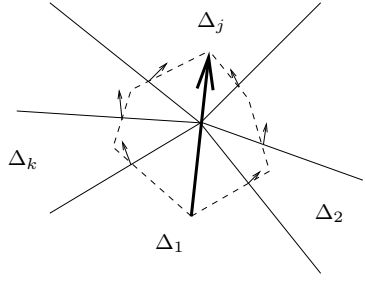}
\caption{{\small Definition of the measure on a geodesic which pass through
    the $1$-stratum.}}\label{somma-figure}
\end{center}
\end{figure}
Notice that every admissible path is
a composition of paths $c_1*\ldots*c_n$ such that every $c_i$ either does not
intersect
the $1$-stratum or intersects only one geodesic. It follows that $\mu_c=\nu_c$.\par
Now let us prove the existence. Let $c$ be an admissible  geodesic path:
notice that $c^{-1}(P)$ is at most a point for every $2$-piece $P$. 
Suppose that $c$ does not intersect any geodesics of the stratification. Then we can define
$\mu_c=\sum_Pa(P)v(P)\delta_{c^{-1}(P)}$ where $v(P)$ is the normal vector to
$P$  such that $\E{v(P)}{\dot c (c^{-1}P)}>0$ (notice that this sum is finite).
Suppose now that $c$ intersects only a geodesic $l$. Let $P_1,\ldots, P_k$ and
$\Delta_1,\ldots,\Delta_k$ be respectively the $2$-pieces and the $3$-pieces
which incide on $l$. We choose the numeration as above and suppose
that $c$ come from $\Delta_1$ and goes into $\Delta_j$. Thus we can define
\[
   \mu_c=\left(\sum_{i=1}^{j-1} a(P_i)u(P_i)\right)\delta_{c^{-1}(l)}
\]
where $u(P_i)$ is the normal vector to $P_i$ which points towards
$\Delta_j$. Since
$\sum_{i=1}^{j-1}a(P_i)u(P_i)=\sum_{i=j}^k a(P_i)u(P_i)$ we have that this
definition does not depend on the numeration (see fig. \ref{somma-figure}).
Now let $c$ an admissible path. Consider a decomposition of $c$ in geodesic
admissible paths $c_1*\ldots*c_k$ such that every $c_i$ intersects only a
geodesic of the stratification or a $2$-piece. Thus we can define $\mu_c$ such that
$\mu_c|{c_i}$ is $\mu_{c_i}$. Notice that this definition does no depends on
the decomposition of $c$.\par
\begin{figure}[!h]
\begin{center}
\input{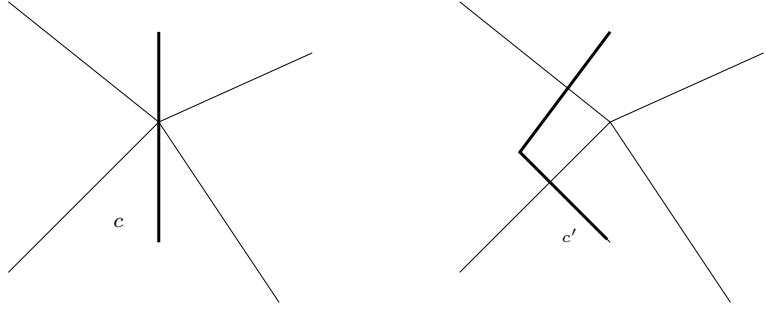}
\caption{{\small Performing the moss in the figure we obtain an arc which does
  not intersect $X_{(1)}$.}}\label{somma3-figure}
\end{center}
\end{figure}

Clearly this measure satisfies the property 1. 2. and 3. of definition
\ref{transverse measure-definition}. Let us show that if $c$ is a closed
admissible path then $\mu_c(c)=0$. First notice that we can assume that $c$
does not intersect any geodesic of $\mathcal C$. In fact if $c$ intersect $l$
we can perform the moss in figure \ref{somma3-figure}. Notice that we obtain a closed
admissible arc $c'$ such that $\mu_{c'}(c')=\mu_c(c)$ and $card(c'\cap
X_{(1)})<card (c\cap X_{(1)})$ (here $X_{(1)}$ is the $1$-statum). Since these
intersections are finite we can reduce to the case $c\cap X_{(1)}=\varnothing$.\par
\begin{figure}[!h]
\begin{center}
\input{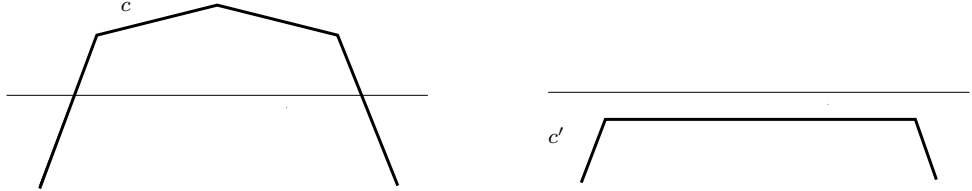}
\caption{{\small Performing the moss in the figure we obtain an arc tight.}}\label{somma2-fig}
\end{center}
\end{figure}

Now $c-X$ is union of component $m_1,\ldots, m_N$, where $m_i$ is an oriented arc whith
endpoints on the faces $P_i^-$ and $P_i^+$. Notice that we can suppose that
the faces $P_i^-$ and $P_i^+$ are different. In fact otherwise we can perform
the moss in figure \ref{somma2-fig}. We call a such path \emph{tight}.\par
Now let $n(P,c)$ the cardinality of the intersection of $c$ with the face $P$. It is clear that
if $n(P,c)=n(P, c')$ for every $P$ then $\mu_c(c)=\mu_{c'}(c')$. On the other
hand let $c$ and $c'$ be tight paths such that there exists a homotopy between
them in $\Kn{H}{n}-X_{(1)}$: it is easy to see that $n(P,c)=n(P,c')$ for every $P$
(in fact for every tight path $c$ there exists $\eps>0$ such that if $c'$ is
in a $\eps$-neighborhood of $c$ then $n(P,c)=n(P,c')$).\par
Thus we have that $\mu_c(c)$ depends only on the homotopy class of $c$ in
$\Kn{H}{n}-X_{(1)}$. Now fix a base point $x_0\in\Kn{H}{n}-X$. 
Since every $\alpha\in\pi_1(\Kn{H}{n}-X^1,x_0)$  is represented by an
admissible path we have a map $\pi_1(\Kn{H}{n}-X^1)\ni
[c]\rightarrow\mu_c(c)\in\Kn{M}{3+1}$. It is quite evident that this map is an
homomorphism.\par
For every geodesic $l\subset X_{(1)}$ consider an admissible path $s_l$ which winds aroud
$l$. Now we can join $s_l$ to $x_0$ with an admissible arc $d$. Consider the
loop $c_l=d*s_l*d^{-1}$.  
Notice that the family $\{c_l\}$ generates $\pi_1(\Kn{H}{n}-X_{(1)},x_0)$. On the
other hand we have $\mu_{c_l}(c_l)=\mu_{s_l}(s_l)$ and since the weights
verify the equation $p_l$ we get $\mu_{c_l}(c_l)=0$. It follows that
$\mu_c(c)=0$.
Thus $\mu$ is a measure on $(\mathcal C, X)$.
\cvd
Now consider a family of weights $a=\{a(P)\}$. We have seen that this family
induces a measure $\mu$ on $\mathcal C$. By  theorem \ref{stratification gives
  domain} we know that there is a domain which corresponds to $\mu$. Let us
denote this domain by $\Omega$: we want to describe the CT-level surface $S=T^{-1}(1)$
and the singularity $\Sigma$.
\begin{figure}[h]
\begin{center}
\input{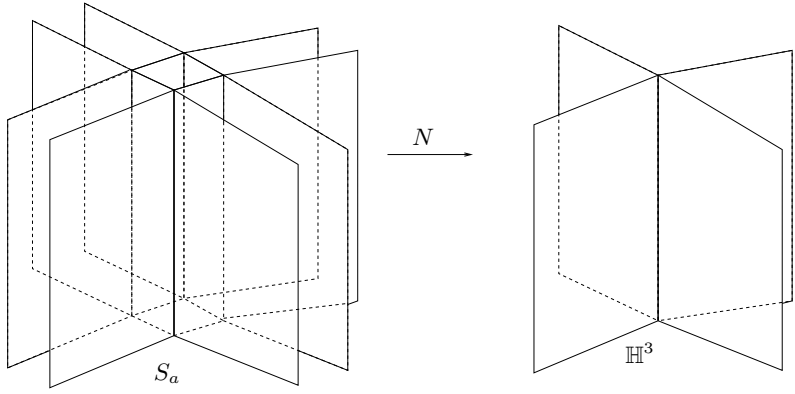}
\caption{{\small The inverse image of a neighborhood of a point in $\Kn{H}{3}$ on the
  CT-level-surface of a regular domain with simplicial singularity.}} 
\end{center}
\end{figure}
\begin{prop}
Consider the decomposition of $S$:  $S_{(i)}:=\{x\in S|\dim N^{-1}N(x)\cap S=i\}$ (for i=0,1,2). Then we have
\[
 \overline S_{(i)}=\sqcup\{N^{-1}(C)\cap S|C\textrm{ is a }3-i\textrm{ piece}\}.
\] 
If $\Delta$ is a $3$-piece $N^{-1}(\Delta)$ is obtained from
$\Delta$ by a traslation and in particular the normal field $N$ restricted to
every component of $S_{0}$ is an isometry.\par
If $P$ is a $2$-piece then $N^{-1}(P)$ is isometric to $P\times[0,a(P)]$ and
the normal field coincides with the projection on the first factor.\par
Finally if $\,l$ is a geodesic piece then $N^{-1}(l)$ is isometric to $l\times
F_l$ where $F_l$ is an euclidean $k$-gon. More precisely, let $P_1,\ldots,P_k$
and $\Delta_1,\ldots,\Delta_k$ be respectively the $2$-piece and $3$-piece which
incide on $l$ ( the numeration is choosen as above).Then there is a
numeration of the $F_l$ edges  $e_1,\ldots,e_k$ such that the lenght of $e_i$ is $a(P_i)$ and the
angle between $e_{i-1}$ and $e_i$ is $\pi-\alpha_i$ where $\alpha_i$
is the dihedral angle of $\Delta_i$ along $l$.
\end{prop}  
\Dim
Fix a base point $x_0$ and for every $3$-piece $\Delta$ let
$\rho_\Delta=\mu_c(c)$ where $c$ is an admissible path which join $x_0$ to
$\Delta$. Clearly $\rho_\Delta$ does not depend on the path $c$ and so it is
well defined. Moreover we have
\[
  \Omega=\bigcap_\Delta\bigcap_{x\in\Delta}\{p\in\Kn{M}{3+1}|\E{p-\rho_\Delta}{x}\}<
  0\}.
\]
Now for every $x\in \Kn{H}{n}$ there exists a unique support plane $P_x$ of
$\Omega$ which is orthogonal to $x$ and such that
$P_x\cap\overline\Omega\neq\varnothing$.
On the other hand we have that $r(N^{-1}(x))=P_x\cap\overline\Omega$.\par
Suppose now that $x\in\Delta$: by  definition of $\Omega$  we have that
 $P_x$ is the plane which passes through $\rho_\Delta$ and is orthogonal to
 $x$.\par
By using inequality (\ref{dominio-formula1}) one can see that 
$\Omega\cap P_x=\{\rho_\Delta\}$ if $x\in\mathrm{int}\Delta$. Let now $x\in P$ where $P$
is a $2$-piece. The same inequality gives that $\Omega\cap
P_x=[\rho_{\Delta_-},\rho_{\Delta_+}]$ where $\Delta_-$ and $\Delta_+$ are the
$3$-pieces which incide on $P$. Finally suppose that $x\in l$ for some
geodesic pieces $l$. We have that $\Omega\cap P_x$ is the convex hull of
$\rho(\Delta_i)$ where $\Delta_1,\ldots,\Delta_k$ are the $3$-pieces which
incide on $l$. Let $P_i$ the $2$-piece which separes $\Delta_i$ from
$\Delta_{i+1}$. Notice that $P_x\cap\overline\Omega$ is a $k$-gon with
vertices $p_i=\rho_{\Delta_i}$. Moreover we have
$\rho_{\Delta_{i+1}}-\rho_{\Delta_i}=a(P_i)v_i$ (where $v_i$ it the normal vector to $P_i$ which point towards
$\Delta_{i+1}$). It is easy to see that the edges of $P_x\cap\overline\Omega$
are $e_i=[p_i,p_{i+1}]$ so that the lenght of $e_i$ is $a(P_i)$.  Moreover the
angle between $e_{i-1}$ and $e_i$ is the equal to the angle between $-v_{i-1}$
and $v_i$. Since the angle between $v_{i-1}$ and $v_i$ is equal to the
dihedral angle of $\Delta_i$ along $l$ we get that $P_x\cap\overline\Omega$ is
isometric to the $k$-gon  $F_l$ defined in the proposition.\par
From this analysis it follows that $N^{-1}(\Delta)\cap S=\Delta+\rho_\Delta$ and so
the normal field is an isometry on $N^{-1}(\Delta)$.\par
Consider now a $2$-piece $P$. We have seen that
$r(N^{-1}(x))=[\rho_{\Delta_-},\rho_{\Delta_+}]$ for all $x\in P$. Thus 
$[\rho_{\Delta_-},\rho_{\Delta_+}]\times P\ni (p,x)\rightarrow p+x\in
N^{-1}(P)\cap S$ is a parametrization: notice that this map is an isometry (in fact
the segment $[\rho_{\Delta_-},\rho_{\Delta_+}]$ is orthogonal to $P$), and the
normal map coincides with the projection on the second factor.\par
An analogous argumentation shows the last statement of the proposition.
\cvd
For $i=0,1,2$ let $\Sigma_i=\{p\in\Sigma|\dim \mathcal F(p)=3-i\}$ (recall that $\mathcal F(p)=N(r^{-1}(p))$).
From the last proposition we imediately get the following corollary:
\begin{corollario}\label{singularity is simplicial}
If $\Omega$ is as above then $\Sigma$ is naturally a cellular complex in the
following sense. $\Sigma_0$ is a numerable set; every component $s$ of $\Sigma_1$ is
an open  segment moreover the closure of $s$ is a closed segment and $\partial
s=\overline s-s$ is contained in $\Sigma_0$; every component $\sigma$ of $\Sigma_2$ is
an open $2$-cell, moreover the closure  of $\sigma$ is a closed $2$-cell and
in fact it is a finite union of
$\sigma$ with segments in $\Sigma_1$ and points in $\Sigma_0$.
\end{corollario}
\cvd
\begin{oss}\emph{
Notice that we have not done any hypotesis of local finiteness of
the cells.}
\end{oss}
\begin{oss}\emph{
If $\mathcal C$ is a simplicial stratification of $\Kn{H}{3}$ we can construct
the dual  complex. For every $3$-piece $\Delta$ there is a vertice
$v_{\Delta}$, for every $2$-piece $P$ there is the segment
$[v_\Delta,v_\Delta']$ where $\Delta$ and $\Delta'$ are the $3$-pieces which
incide on $P$, and for every $1$-piece $l$ there is the polygon with vertices 
$v_{\Delta_1},\ldots, v_{\Delta_n}$ where $\Delta_1,\ldots,\Delta_n$ are the
pieces which incide on $l$.}\par
\emph{
Notice that $\Sigma$ is combinatorially equivalent to the dual complex of
$\mathcal C$. 
Thus the combinatorial structure of $\Sigma$ depends only on the
combiatorial structure of $\mathcal C$. This remark points out the duality
between stratifications and singularities.}
\end{oss}
We want to describe a class of regular domains whose stratification is
simplicial. As we are going to see, there are regular domains with
simplicial stratification which do not belong to this class. However 
the domain $\dom$ has simplicial stratification if
and only if it belongs to this class.
\begin{defi}
The singularity $\Sigma$ of a future complete regular domain  is \emph{simplicial} if there
exists a cellularization $\Sigma_0\subset\Sigma_1\subset\Sigma_2$ such that
$\Sigma_0$ is a discrete set,
every component of $\Sigma_1-\Sigma_0$ is a straight segment and every
component of $\Sigma_2-\Sigma_1$ is a finite-sided spacelike polygon with
vertices  in $\Sigma_0$ and edges in $\Sigma_1$.
\end{defi}
\begin{figure}
\begin{center}
\input{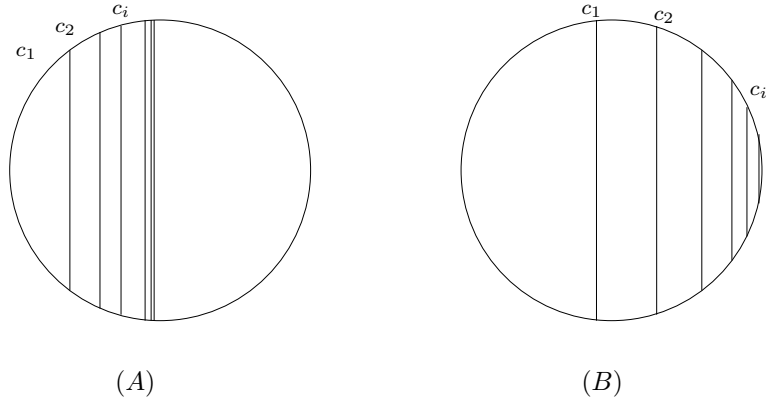}
\caption{On the left a non-simplicial stratification and on the right a
  simplicial stratification.}\label{lam-fig}
\end{center}
\end{figure}
\begin{oss}\emph{
We shall  prove that a regular domain with simplicial singularity has simplicial
stratification. If we do not require that $\Sigma_0$ is discrete this result
is no more true. Consider for istance the stratification of $\Kn{H}{2}$ given
in fig.\ref{lam-fig} (A) (here the example is given for $n=2$, but a very
analogous example holds for $n=3$). It is easy to construct a regular domain 
with a such stratification and singularity with a cell decomposition.}\par
\emph{
On the other han there exists regular domain with simplicial stratification
which has not simplicial singularity. For istance consider the stratification
of $\Kn{H}{2}$ given in fig.\ref{lam-fig} (B). It is easy to construct a regular
domain with a such stratification whose singularity $\Sigma$ is compact in
$\Kn{M}{2+1}$. It follows that $\Sigma_0$ is not discrete.
}
\end{oss} 
\begin{oss}\label{singularity distance}\emph{
Let $\Sigma$ be a simplicial singularity of a future complete regular domain
$\Omega$. Notice that we can consider on $\Sigma$ a weak topology induced by
the cellularization ($A\subset\Sigma$ is open in this topology if and only if
the intersection of $A$ with every open cell is open). 
Since we do not require the local finitness of the cells, generally this
topology is finer than the topology induced by $\Kn{M}{3+1}$.}\par 
\emph{
Notice that  every open cell
has a natural distance. Thus  if $c$ is a path in $\Sigma$ we can define the
lenght of $c$ as the sum of the lenghts of the intersections of $c$ with the
cells of $\Sigma$. Finally we can define a path-distance on $\Sigma$ such that 
$d_\Sigma(r,r')$ is the infimum of the lenghts of the paths from $r$ to $r'$.
It is easy to see that this distance agrees with the natural distance on
$\Sigma$ described in section \ref{Convergence-sec}. Thus  the
topology induced by $\Kn{M}{3+1}$ on $\Sigma$ generally is finer than the topology induced by $d_\Sigma$.
}
\end{oss}
Now we can prove that at least the regular domains with a proper and surjective
normal field and simplicial singularity  are
given by a measured simplicial stratification of $\Kn{H}{3}$.
We start with a technical lemma.
\begin{lemma}\label{simplicial singularity-lemma}
Let $\Omega$ be a regular domain with proper normal field $N$ and simplicial singularity. We have
$\Sigma_i=\{p\in\Sigma|\dim \mathcal F(p)\geq 3-i\}$.
\end{lemma}
\Dim
By proposition \ref{domain gives stratification} we have that if $p$ and $q$
lie in the same component of $\Sigma_1-\Sigma_0$ or $\Sigma_2-\Sigma_1$
then $\mathcal F(p)=\mathcal F(q)$. Let $p\in\Sigma_1-\Sigma_0$ and let $[r_0,r_1]$ be a
component of $p$ in  $\Sigma_1-\Sigma_0$. Take $q\in (r_0,r_1)$, since $\ort{p-q}$ separes
$\mathcal F(p)$ from $\mathcal F(q)$ it follows that this plane contains  $\mathcal F(p)$, then $\dim \mathcal F(p)<3$.
In analogous way we can prove that if  $p\in\Sigma_2-\Sigma_1$ then $\dim
\mathcal F(p)<2$. This proves an inclusion.\par
Fix now a component $[r_0,r_1]$ of $\Sigma_1-\Sigma_0$. We want to prove that
$\dim \mathcal F(p)>1$ for $p\in [r_0,r_1]$. First suppose that there not
exists any cell
in $\Sigma_2-\Sigma_1$ which incides on $p$. Since $\Sigma$ is contractile it
follows that $\Sigma-(r_0,r_1)$ is not connected. Thus $\Omega-r^{-1}(r_0,r_1)$ is
not connected and one easily see that it is possible only if $\mathcal F(p)$ is a plane
for $p\in(r_0,r_1)$. Suppose now that there exists a
component $T$ of $\Sigma_2-\Sigma_1$ which incides on $[r_0,r_1]$. We want to
see that this component is not unique. Let $l$ be the geodesic which
corresponds to $T$ (i.e. the plane which contains $T$ is orthogonal to $l$). 
Since $[r_0,r_1]$ is an edge of $l$ there is a vector $v\in\ort l$ such that $v$ is orthogonal
to $[r_0,r_1]$ and $\E{x-q}{v}<0$ for all $x\in T-[r_0,r_1]$ and all $q\in
[r_0,r_1]$. Fix now a point $x\in l$ and consider a geodesic segment $c$ starting from
$x$ with direction $v$. Since $N$ is a proper map we have that $N^{-1}(c)$ is
compact. So that there exists a sequence $p_k\in N^{-1}(c)$ such that $p_k\rightarrow
p$ with $N(p)=x$. Put $\rho_k=r(p_k)$ we have $\E{\rho_k-r_0}{N(p_k)-x}\geq
0$. On the other  hand the direction of $N(p_k)-x$ converges to $v$ so that we
have
$\E{r(p)-r_0}{v}=0$. Since $N(p)\in l$ we have $\rho_\infty=r(p)\in[r_0,r_1]$. Since
$rN^{-1}(c)$ is compact and since $\Sigma_0$ is discrete it is easy to see
that there is only a finite number of cells which intersects $rN^{-1}(c)$. So
that there exists a cell $T'$ such that $\rho_k\in\overline T'$ for
$k>>0$ so that $\rho_\infty\in T\cap T'$.
Since there exists only a finite number of pieces which correspond to the points in
$T'$ it follows that we can suppose that $\mathcal F(\rho_k)$ does not depend on
$k$. We know that $N(p_k)\in c$ so that the segment $c$ is contained in
$\mathcal F(\rho_k)$. 
On the other hand it is easy to see that $\mathcal F(\rho_k)\subset\mathcal F(\rho_\infty)$ so that
$C=\mathcal F(\rho_\infty)\cap\ort{(r_1-r_0)}$ is a face of $\mathcal
F(\rho_\infty)$ with dimension equal to $2$ (in fact it contains $l$ and the segment $c$). 
Since $\rho_\infty\in[r_0,r_1]$  we have that  $C\subset \mathcal F(r)$ for
all $r\in [r_0,r_1]$. 
An analogous argument shows that $\dim \mathcal F(p)=3$ for all vertices $p$.
\cvd

\begin{prop}
Let $\Omega$ be a future complete regular domain  with simplicial singularity. Suppose
that the normal field $N$ is a proper map. Then the stratification associated
with $\Omega$ is simplicial. Moreover there exists a unique measure $\mu$ on
$\mathcal C$ such that $\Omega$ is equal (up to traslations) to the domain associated with $\mu$.
\end{prop}
\Dim
Let $\mathcal C$ be the stratification associated with $\dom$. We show that
$\mathcal C$ is simplicial. 
Fix a point $x$ in $\Kn{H}{3}$ we have to construct a neighborhood of $x$ which
intersects only a finite number of pieces. If $x$ does not lie in  the
$2$-stratum $X$ then it in the interior of a $3$-piece $\Delta$ and this neighborhood
intersects only $\Delta$.
Suppose now that $x\in X-X_{(1)}$ where $X_{(1)}$ is the $1$-stratum. There is a
$2$-piece $P$ such that $x\in P-\partial P$.
Now by the hypotesis $rN^{-1}(P)$ is a segment $[r_0,r_1]$ and
moreover there exists $3$-pieces $\Delta_0$ and $\Delta_1$ such that
$\{r_i\}=rN^{-1}(\Delta_i)$. By lemma \ref{domain gives stratification} we have that
$P$ is a face of $\Delta_i$ and moreover $P=\Delta_0\cap\Delta_1$.
Thus $x\in\mathrm{int}(\Delta_0\cup\Delta_1)$. Moreover 
$\mathrm{int}(\Delta_0\cup\Delta_1)$ intersects only $P,\Delta_0$ and
$\Delta_1$.\par
Finally suppose $x\in l$ where $l$ is $1$-piece of $\mathcal C$ . 
Let $F=rN^{-1}(l)$: it is a compact polygon with
vertices $p_1,\ldots,p_k$. For every $p_i$ let $\Delta_i$ be the piece such
that $\{p_i\}=rN^{-1}(\Delta_i)$. We have that $l$ is equal to the intesection of
$\Delta_i$. Furthermore $\Delta_i\cap\Delta_{i+1}$ is a face $P_i$ and
$\Delta_i\cap\Delta_j=\{l\}$ if the vertices $p_i,p_j$ are not adjacent.
Notice that the dihedral angle of $\Delta_i$ along $l$ is equal to $\pi-\alpha_i$
where $\alpha_i$ is the angle of $F$ at $p_i$: in fact $p_{i+1}-p_i$ is a
orthogonal vector to $P_i$ which points toward $\Delta_{i+1}$. It follows that the
sum of dihedral angle of $\Delta_i$ along $l$ is $2\pi$ and so
$\mathrm{int}(\bigcup\Delta_i)$ is a neighborhood of $x$. Notice that this neghborhood intersects only
$l,\Delta_1,\ldots,\Delta_n$ and $P_1,\ldots, P_n$. Thus $\mathcal C$ is a
simplicial stratification.\par
Now we can define a family of weights $\{a(P)\}$. In fact if $P$ is a
$2$-piece we know that $rN^{-1}(P)$ is a spacelike segment $[r_0,r_1]$, thus we can
define $a(P)=\left(\E{r_1-r_0}{r_1-r_0}\right)^{1/2}$. Let us show that this
is a family of weights on $\mathcal C$. Fix a geodesic $l$ and consider the
$2$-pieces $P_1,\ldots, P_k$ and the $3$-pieces $\Delta_1,\ldots,\Delta_k$
which incides on $l$. We suppose $P_i=\Delta_{i-1}\cap\Delta_i$ and take $i\mod k$. Let
$\{p_i\}=rN^{-1}(\Delta_i)$ we know that $p_{i+1}-p_i$ is an orthogonal vector
to $P_i$ which points toward $\Delta_{i+1}$. So that if $v_i$ is the normal
vector to $P_i$ which points toward $\Delta_{i+1}$ we have
$p_{i+1}-p_i=a(P_i)v_i$. Thus 
\[
  \sum_i a(P_i)v_i=\sum_i p_{i+1}-p_i=0.
\]
Let $\mu$ be the measure which corresponds to the family $\{a(P)\}$, we have
to show that up to traslation $\Omega$ is the domain which corresponds to the
measure $\mu$.\par
Fix a base point $x_0\in\Kn{H}{n}-X^2$. Up to traslation we can suppose that $rN^{-1}(x_0)=0$. 
Now let $p$ a vertex of $\Sigma$, by construction is evident that $p=\mu_c(c)$
where $c$ is an admissible path which starts in $x_0$ and termines in the
piece which corresponds to $p$. It follows that $\Omega$ is the regular domain 
 which corresponds to the measure $\mu$.
\cvd

In the last part of this section we want to study the $\Gamma$-invariant
simplicial geodesic laminations where $\Gamma$ is a free-torsion discrete
co-compact subgroup of $\SOO^+(3,1)$.
We see that for a $\Gamma$-invariant simplicial lamination the set of measures
on it is parametrized by a a finite number of positive number which satisfies
a finite set of linear equations.\par
We start with some remarks about $\Gamma$-invariant simplicial
stratifications.
\begin{prop}\label{sss}
Let $\mathcal C$ be a $\Gamma$-invariant simplicial stratification of
 $\Kn{H}{3}$. Then $\pi(C)$ is compact for all $C\in\mathcal C$ . In particular the projection
of a $1$-piece is a simple geodesic whereas the projection of a $2$-piece is
either a closed hyperbolic surface or a hyperbolic surface with geodesic
boundary. 
Finally  there are only a finite number of pieces up to the action of $\Gamma$.
\end{prop}
\Dim
Since $\mathcal C$ is simplicial one easily see that $\Gamma\cdot C$ is closed
and the projection of $C$ is compact. Thus the projection of a $1$-pieces $l$ is a
compact complete geodesic and so it is closed. Since the orbit of $l$ is
formed by the disjoint union of geodesics it follows that the projection is a
simple geodesic.\par
An analogous argument shows that the projection of a $2$-piece $P$ is a
hyperbolic surface. Notice that it is closed if $P$ is a plane otherwise it has
totally geodesic boundary. \par
Let $K$ be a compact fundamental region for $\Gamma$: since $K$ intersects
only a finite number of pieces we get the last statement.
\cvd
Fix a $\Gamma$-invariant simplicial stratification $\mathcal C$. We denote by
$T_{\mathcal C}$ the projection of the $2$-stratum $X$ onto the quotient
$M=\Kn{H}{3}/\Gamma$. Notice that there exists a finite set of simple
geodesics $\{c_1,\ldots,c_N\}$  of $M$ such that $c_i\subset T_\mathcal C$ and
$T_{\mathcal C}-\bigcup c_i$ is a finite union of totally geodesic submanifold
$F_1,\ldots, F_L$ such that $\overline F_i=F_i\cup c_{i_1}\cup\ldots\cup
c_{i_k}$. 
The geodesics $c_i$ are called the \emph{edges} of the
surface whereas the surfaces $F_i$ are  the \emph{faces}. 
A subset $X\subset M$ which has a such decomposition is called
\emph{piece-wise geodesic surface}.
Notice that piece-wise geodesic  surfaces correspond bijectively with
$\Gamma$-invariant geodesic stratification.\par
No let $\mu$ be a $\Gamma$-invariant measure on $\mathcal C$. Notice that the
family of weights $\{a(P)\}$ associated with it satisfies $a(\gamma P)=a(P)$
for all $2$-pieces $P$ and all $\gamma\in\Gamma$. Thus there exists a family
 of positive constants $\{\alpha(F)\}$ parametrized by the face of
$T_{\mathcal C}$ such that $a(P)=\alpha(\pi P)$.\par 
This remark suggests the following definition. Let $T$ be a piece-wise geodesic
surface in $M$ and let $\mathcal C$ be the stratification
associated. 
A\emph{ family of weights} on $T$ is a family $\{\alpha(F)\}$ parametrized by the faces of $T$ such that 
 $\{a(P):=\alpha(\pi(P))\}$ is a family of weights on $\mathcal C$. Notice
 that for every $1$-piece $l$ we have that the equation associated with $l$
 and equation associated with $\gamma(l)$ are realted by the identity
\[
    p_{\gamma (l)}(a)= \gamma p_l(a)
\]
so that the solutions of the equation $p_l$ coincide with the solutions of
$p_{\gamma l}$. Thus the conditions we have to impose to ensure
$\{\alpha(F)\}$ is a family of weights can be parametrized by the edges of $T$.\par
Finally we have that the families of weights on $T$
correspond bijectively with $\Gamma$-invariant measures on $\mathcal C$.
Notice that $T$ is a piece-wise geodesic  surface with $f$ faces and $e$
edges then the weights on $T$ corresponds to a subset of $\mathbb R_+^f$ defined by
$2e$ linear equations (in fact every $p_l$ is equivalent to $2$ linear equations).
Thus if there exists positive solutions they forms a convex cone of dimension
greater than $f-2e$.\par
\begin{es}\emph{
Now we exhibit some examples of piece-wise geodesic surfaces. Fix an
hyperbolic $3$-manifolds with totally geodesic boundary $N$ and cosider the
canonical decomposition of $N$ in truncated polyedra (see \cite{Kojima} for
the definition). Let $M$ be the double of $N$, notice that the double of the  
$2$-skeleton of the decomposition of $N$ give a piece-wise geodesic surface
$T$.}\par
\emph{
Suppose that every polyedron of the decomposition is a tetrahedron. 
We want to exstimate the number of the edges and the faces of $T$. On
the boundary of $N$ the decomposition gives a triangulation. Let $v,l,t$ be
respectively the number of vertices, edges and faces of this triangulation. We
have that $v-l+t=2-2g$ where $g$ is the genus of $\partial N$ ($g\geq 2$). On
the other hand we have that $t= 2/3 l$ so that we get that
$v-1/3 l = 2-2g$. On the other hand let $e,f$ be respectively the number of
edges and faces of $T$. We have that $v=2e$ (in fact every edge of $T$
intersects in two vertces $\partial N$) and $l=3e$.
It follows that $2e-f=2-2g<0$.}
\end{es}  
We conclude this section proving that regular domains which are invariant for
some affine deformation of $\Gamma$ correspond bijectively to
$\Gamma$-invariant measured simplicial stratifications.
\begin{corollario}
There exists a bijective corrispondence between $\Gamma$-invariant measured
simplicial stratifications of $\Kn{H}{3}$ and future complete regular domains
which are invariant for some affine deformation $\Gamma_\tau$ of $\Gamma$ and
has simplicial singularity.
\end{corollario}
\Dim
It is sufficient to show that $\Gamma$-invariant measured simplicial
stratifications give domains with simplicial singularity. Now fix a
$\Gamma$-invariant measured
simplicial stratification $(\mathcal C, X,\mu)$ and let $\{a(P)\}$ be the
family of weights associated with it. By proposition \ref{sss} we get that there exists
$a>0$ such that $a\leq a(P)$ for all $2$-faces $P$. Now for a $3$-piece
$\Delta$ let $\rho_\Delta$ be the corresponding point on $\Sigma$. We have
that $\Sigma_0=\{\rho_\Delta|\Delta\textrm{ is a }3\textrm{-piece}\}$. On the other hand we know that
$d_\Sigma(\rho_\Delta, \rho_\Delta')\geq a$ so that $\Sigma_0$ is a discrete
set.
\cvd

\section{Conclusions}\label{conclusions-section}
Let $\Gamma$ be  a free-torsion co-compact discrete subgroup of
$\SOO^+(n,1)$. We have seen that for every cocycle $\tau\in
Z^1(\Gamma, \Kn{R}{n+1})$ 
there exists a unique future complete future complete regular domain $\dom$ which is
$\Gamma_\tau$-invariant. So  the future complete regular domains arise
naturally in the study of the Lorentz flat structures on $\mathbb R\times M$.\par
In section \ref{CT-section}  we have seen that every future complete future complete regular domain is
associated with a $\Gamma$-invariant geodesic stratification of $\Kn{H}{n}$. On the other hand in
section \ref{measure-section}  we have seen that given a $\Gamma$-invariant measured geodesic stratification,
we get a future complete regular domain which is invariant for an affine
deformation of $\Gamma$. Moreover in dimension $n=2$ this
corrispondence agree with the Mess identification between measured geodesic
lamination and future complete regular domain.\par
We can ask if this corrispondence is in fact an identification in all dimensions. 
We have seen in section \ref{simplicial-section} that this corrispondence induces an identification between
simplicial stratification and future complete regular domain with simplicial
singularity.\par
The general case seems more difficult. Given a future complete regular domain $\Omega$
we should construct a \textbf{measured geodesic stratification} $(\mathcal C,
Y,\mu)$ which gives $\Omega$. By looking at the construction of a domain
$\Omega$ we have $Y=\{x|\#(N^{-1}(x))>1\}$ and in fact it is easy to see
that this set has null Lebesgue measure in $\Kn{H}{n}$. Now 
suppose that for all admissible paths $c:[0,1]\rightarrow\Kn{H}{n}$ there
exists a \emph{Lipschitz path} $\tilde c:[0,1]\rightarrow \widetilde S_1$ such that
$N(\tilde c([0,1]))=c([0,1])$. Under this assumption we could define a measure $\mu$
on an admissible path in this way.   
Since the retraction is locally Lipshitz, the map $r(t)=r(\tilde c(t))$ is Lipshitz
so that it is differentiable almost everywhere. Consider the  $\mathbb
R^{n+1}$-valued measure $\tilde\mu$  on $[0,1]$  defined by the identity
\[
   \tilde\mu(E)=\int_Er'(t)\mathrm dt.
\]
Then we could define the transverse measure $\mu_c$ as the image of the measure
$\tilde \mu$:
\[
   \mu_c=N_*(\tilde\mu).
\]
Notice that the assumption is always verified for the regular domains $\Omega$
which arise from a measured geodesic stratifications. 
In fact  given an admissible path $c$ the
Lipschitz path $\tilde c$ always exists. In fact we
can define 
\[ 
\tilde c(t)=c(t)+\int_0^t\mu_c.
\]
Thus the problem is: given an admissible path $c:[0,1]\rightarrow\Kn{H}{n}$   
we have to find a rectifiable curve $\tilde c\subset \widetilde S_1$ such that $N(\tilde
c)$ is the curve $c$ (notice that if $\widetilde S_a$ is strictly convex there exists a
unique curve such that $N(\tilde c)=c$ but we do not know if such arc is
rectifiable).\par
In dimension $n=2$ one easily see that this problem has always solution
because if $v_1$ and $v_2$ are the orthogonal vector to two leaves of the
lamination then they generate a timelike vector space so that
\[
   |\E{v_1}{v_2}|\geq\E{v_1}{v_1}^{1/2}\E{v_2}{v_2}^{1/2}.
\]
By using this inequality one can see that the lenght of $\tilde c$ is lesser
than 
\[
\ell(c)+\E{r(\tilde c(1))-r(\tilde c(0))}{r(\tilde c(1))-r(\tilde c(0))}
\]
where $\ell(c)$ is the lenght of $\Kn{H}{n}$.
Unfortunately in dimension $n\geq 3$ this argument fails. 

\section*{Acknowledgements}
The author gratefully thanks Riccardo Benedetti for many helpful conversations.

\end{document}